\documentclass[twoside,reqno]{amsart}
\usepackage[square,comma]{natbib}
\usepackage[english]{babel} 
\usepackage{exscale}   
\usepackage{amsmath}
\usepackage{amsfonts}
\usepackage{amssymb}
\usepackage{amsxtra}
\usepackage{stmaryrd}
\usepackage{xspace}
\usepackage{epsfig}
\usepackage{mathrsfs}
\usepackage{pstricks}

\numberwithin{equation}{section}
\newtheorem{Theorem}{Theorem}[section]
\newtheorem{Lemma}[Theorem]{Lemma}
\newtheorem{Corollary}[Theorem]{Corollary}
\newtheorem{Proposition}[Theorem]{Proposition}
\newtheorem{Definition}[Theorem]{Definition}
\newtheorem{Remark}[Theorem]{Remark}

\newtheorem{Example}[Theorem]{Example}
\newtheorem{Notation}[Theorem]{Notation}

\newtheorem{Conjecture}[Theorem]{Conjecture}

\newcommand{\Bset}{\mathbb{B}}
\newcommand{\Cset}{\mathbb{C}}

\newcommand{\Nset}{\mathbb{N}}
\newcommand{\Fset}{\mathbb{F}}

\newcommand{\Rset}{\mathbb{R}}
\newcommand{\Sset}{\mathbb{S}}

\newcommand{\Zset}{\mathbb{Z}} 
\newcommand{\cA}{\ensuremath{{\mathcal A}}\xspace}         
\newcommand{\cB}{\ensuremath{{\mathcal B}}\xspace}         
\newcommand{\cC}{\ensuremath{{\mathcal C}}\xspace}         
\newcommand{\cM}{\ensuremath{{\mathcal M}}\xspace}         
\newcommand{\cN}{\ensuremath{{\mathcal N}}\xspace}         
\newcommand{\cZ}{\ensuremath{{\mathcal Z}}\xspace}         

\newcommand{\scrI}{\ensuremath{{\mathscr{I}}}\xspace} %
\newcommand{\scrtI}{\ensuremath{{\tilde{\mathscr{I}}}}\xspace} %
\newcommand{\ii}{\mathbf{i}}
\newcommand{\jj}{\mathbf{j}}

\newcommand{\tgamma}{\ensuremath{\tilde{\gamma}}\xspace} %
\newcommand{\tsigma}{\ensuremath{{\widetilde{\sigma}}}} %
\newcommand{\1}{\ensuremath{{\rm 1\kern-.25em l}}\xspace}  
\newcommand{\alg}{\operatorname{alg}} 
\newcommand{\Ad}{\mathop{\mathrm{Ad}}}                     
\newcommand{\distr}{\operatorname{distr}}                  
\renewcommand{\i}{\mathrm{i}}                              
\newcommand{\id}{\operatorname{id}}                        
\newcommand{\inv}{\operatorname{inv}}
\newcommand{\sh}{\operatorname{sh}}           
\newcommand{\tw}{\operatorname{tw}}           
\newcommand{\trace}{\operatorname{tr}}        
\newcommand{\tail}{\ensuremath{\mathrm{tail}}}                 
\newcommand{\Aut}[1]{\ensuremath{\operatorname{Aut}(#1)}\xspace}  

\newcommand{\set}[2]{\mathopen{\{}#1\mathop{|}#2\mathclose{\}}}
\newcommand{\linh}{\ensuremath{\operatorname{lin}}\xspace}  

\newcommand{\sotlim}{\ensuremath{\textsc{sot-}\lim}\xspace}   
\newcommand{\sot}{\ensuremath{\textsc{sot}}\xspace}           

\newsavebox{\artin}
\newcommand{\masterartin}{
\linethickness{1pt}
\qbezier(10,20)(0,10)(-10,0)
\qbezier(-10,20)(-7,17)(-4,14)
\qbezier(10,0)(6,4)(4,6)
}
\newsavebox{\artininv}
\newcommand{\masterartininv}{
\linethickness{1pt}
\qbezier(-10,20)(0,10)(10,0)
\qbezier(10,20)(6,16)(4,14) 
\qbezier(-10,0)(-5,5)(-4,6)
}
\newsavebox{\strandr}
\newcommand{\masterstrandr}{
\linethickness{1pt}
\qbezier(10,20)(10,10)(10,0)
}
\newsavebox{\strandl}
\newcommand{\masterstrandl}{
\linethickness{1pt}
\qbezier(-10,20)(-10,10)(-10,0)
}
\newsavebox{\horizontaldots}
\newcommand{\masterhorizontaldots}{
\linethickness{1pt}
\put(5,0){\circle*{2}}
\put(11,0){\circle*{2}}
\put(17,0){\circle*{2}}
}

\begin{document}
\title[Noncommutative independence from the braid group $\Bset_\infty$]
{Noncommutative independence\\ from the braid group $\Bset_\infty$}
\author{Rolf Gohm}
\author{Claus K\"ostler}
\address{Institute of Mathematics and Physics\\ 
         Aberystwyth University\\
         Aberystwyth, SY23 3BZ, UK}
\address{Department of Mathematics\\
         University of Illinois at Urbana-Champaign\\          
         Altgeld Hall, 1409 West Green Street\\
         Urbana, 61801, USA}
\email{rog@aber.ac.uk}      
\email{koestler@math.uiuc.edu}
\subjclass[2000]{Primary 46L53; Secondary 20F36}
\keywords{Braid groups, braid group von Neumann algebras, 
fixed point algebras, distributional symmetries, 
exchangeability, spreadability, noncommutative independence, 
noncommutative extended de Finetti theorem, braided Hewitt-Savage 
Zero-One Law, noncommutative random sequences, noncommutative 
stationary processes, noncommutative Bernoulli shifts} 
\begin{abstract}
We introduce `braidability' as a new symmetry for (infinite) sequences 
of noncommutative random variables related to representations of the
braid group $\Bset_\infty$. It provides an extension of exchangeability 
which is tied to the symmetric group $\Sset_\infty$. Our key result is 
that braidability implies spreadability and thus conditional  
independence, according to the noncommutative extended de Finetti theorem 
\cite{Koes08aPP}. This endows the braid groups $\Bset_n$ with a new 
intrinsic (quantum) probabilistic interpretation. We underline this 
interpretation by a braided extension of the Hewitt-Savage Zero-One Law. 

Furthermore we use the concept of product representations of endomorphisms 
\cite{Gohm04a} with respect to certain Galois type towers of fixed point 
algebras to show that braidability produces triangular towers of commuting 
squares and noncommutative Bernoulli shifts. As a specific case 
we study the left regular representation of $\Bset_\infty$ and the 
irreducible subfactor with infinite Jones index in the non-hyperfinite 
$II_1$-factor $L(\Bset_\infty)$ related to it. Our investigations reveal 
a new presentation of the braid group $\Bset_\infty$, the `square root of 
free generator presentation' $\Fset_\infty^{1/2}$. These new generators give 
rise to braidability while the squares of them yield a free family. Hence 
our results provide another facet of the strong connection between 
subfactors and free probability theory \cite{GJS07aPP}; and we speculate 
about braidability as an extension of (amalgamated) freeness on the 
combinatorial level. 
\end{abstract}
\maketitle 
{%
\def\widedotfill{\leaders\hbox to 10pt{\hfil.\hfil}\hfill}
\def\pg#1{\widedotfill\rlap{\hbox to 15pt{\hfill{\small#1}}}\par}
\rightskip=15pt\leftskip=10pt%
\newcommand{\sct}[2]{\noindent\llap{\hbox to%
    10pt{{#1}\hfill}}~\mbox{#2~}}
\newcommand{\separ}{\vspace{.2em}}

\hrule
\section*{Contents}
\sct{}{Introduction and main results}%
\pg{\pageref{section:intro}}%
\separ
\sct{\ref{section:symmetries}}{Distributional symmetries}%
\pg{\pageref{section:symmetries}}%
\separ
\sct{\ref{section:sequences-braid}}{Random variables generated by the braid group $\Bset_\infty$}%
\pg{\pageref{section:sequences-braid}}%
\separ
\sct{\ref{section:endomorphisms-braid}}{Endomorphisms generated by the braid group $\Bset_\infty$}%
\pg{\pageref{section:endomorphisms-braid}}%
\separ
\sct{\ref{section:presentation}}{Another braid group presentation, $k$-shifts and braid handles}%
\pg{\pageref{section:presentation}}%
\separ
\sct{\ref{section:left-regular-rep}}{An application to the group von Neumann algebra $L(\Bset_\infty)$}%
\pg{\pageref{section:left-regular-rep}}%
\separ
\sct{\ref{section:examples}}{Some concrete examples}%
\pg{\pageref{section:examples}}%
\separ
\sct{\ref{section:appendix}}{Appendix: Operator algebraic noncommutative probability}%
\pg{\pageref{section:appendix}}%
\separ
\sct{}{References}%
\pg{\pageref{section:bibliography}}
}
\vspace{12pt}
\hrule

\section*{Introduction and main results}
\label{section:intro}
\begin{figure}[h]
\setlength{\unitlength}{0.3mm}
\begin{picture}(360,35)
\savebox{\artin}(20,20)[1]{\masterartin} 
\savebox{\artininv}(20,20)[1]{\masterartininv} 
\savebox{\strandr}(20,20)[1]{\masterstrandr} 
\savebox{\strandl}(20,20)[1]{\masterstrandl} 
\savebox{\horizontaldots}(20,20)[1]{\masterhorizontaldots}
\put(0,00){\usebox{\strandl}}  
\put(0,00){\usebox{\strandr}}
\put(20,00){\usebox{\horizontaldots}}
\put(60,0){\usebox{\strandl}}
\put(80,0){\usebox{\artin}}
\put(100,0){\usebox{\strandr}}
\put(120,0){\usebox{\horizontaldots}}
\put(-1,25){\footnotesize{$0$}}
\put(19,25){\footnotesize{$1$}}
\put(72,25){\footnotesize{$i-1$}}
\put(100,25){\footnotesize{$i$}}
\put(200,00){\usebox{\strandl}}  
\put(200,00){\usebox{\strandr}}
\put(220,00){\usebox{\horizontaldots}}
\put(260,00){\usebox{\strandl}}
\put(280,00){\usebox{\artininv}}
\put(300,00){\usebox{\strandr}}
\put(320,00){\usebox{\horizontaldots}}
\put(199,25){\footnotesize{$0$}}
\put(219,25){\footnotesize{$1$}}
\put(272,25){\footnotesize{$i-1$}}
\put(300,25){\footnotesize{$i$}}
\end{picture}
\caption{Artin generators $\sigma_i$ (left) and $\sigma_i^{-1}$ (right)}
\label{figure:artin}
\end{figure}
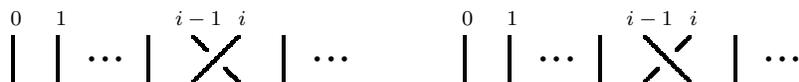

The braid groups $\Bset_n$ were introduced by Artin in \cite{Arti1925a} 
where it is shown that, for $n \ge 2$, $\Bset_n$ is presented by $n-1$ 
generators $\sigma_1,\ldots,\sigma_{n-1}$ satisfying the relations
\begin{align}
&&\sigma_i \sigma_{j} \sigma_i 
&= \sigma_{j} \sigma_i \sigma_{j} 
&\text{if $ \; \mid i-j \mid\, = 1 $;}&& \tag{B1} \label{eq:B1}\\
&&\sigma_i \sigma_j 
&= \sigma_j \sigma_i  
&\text{if $ \; \mid i-j \mid\, > 1 $.}&& \tag{B2} \label{eq:B2}
\end{align}
One has the inclusions $ \Bset_2 \subset \Bset_3 \subset \cdots \subset 
\Bset_{\infty}$, where  $\Bset_{\infty}$ denotes the inductive limit. 
For notational convenience, $\sigma_0$ will denote the unit element in 
$\Bset_\infty$ and $\Bset_1$ is the subgroup $\langle \sigma_0 \rangle$. 
The Artin generator $\sigma_i$ and its inverse $\sigma_i^{-1}$ will be
presented as geometric braids according to Figure \ref{figure:artin}.

Due to their rich algebraic and topological properties, braid groups 
are a key structure in mathematics and their better understanding is 
crucial for many applications, for example entanglement in quantum 
information theory \cite{Kauf07a}. Of special interest for us will be 
that the braid group $\Bset_n$ is an extension of the symmetric group 
$\Sset_{n}$ and contains the free group $\Fset_{n-1}$ as a 
subgroup \cite{Bir75a,BiBr05a}.

Aside from their algebraic and topological properties, $\Sset_n$ and $\Fset_n$ have also intrinsic probabilistic interpretations which are connected to independence structures. This was revealed for the symmetric groups $\Sset_n$ already in the 1930s by the celebrated work of de Finetti on exchangeability [Fin31]. Here the groups $\Sset_n$ are represented by automorphisms of the underlying probability space and we will generalize this idea to braid groups in this paper. In the case of the free groups $\Fset_n$ a breakthrough result had to wait until the 1980s when Voiculescu discovered freeness during his investigations of free group von Neumann algebras [Voi85] and soon after established its intimate connection to random matrix theory [Voi91]. Even more directly, the $\Fset_n$'s serve as noncommutative models for the underlying probability spaces themselves.

On the other hand, it is known that braid groups carry a probabilistic 
interpretation through quantum symmetries \cite{EvKa98a}, in particular 
quantum groups \cite{Maji95a,FSS03a}. Their statistical and entropy 
properties are examined in \cite{DeNe97a,DeNe98a,VeNeBi00a,NeVo05a,MaMa07a} and a physical interpretation in terms of quantum coin tosses is given in 
\cite{KuMa98a}. 
Moreover, a probabilistic facet of braid groups is apparent in subfactor 
theory, for example from Markov traces or commuting squares 
\cite{Jone83a,Popa83a,GHJ89a,Popa90a,Jone91a,Popa93a,JoSu97a}.   

In this paper we take a new way towards a probabilistic interpretation 
of braid groups and look at braided structures from the perspective of 
distributional symmetries and invariance principles \cite{Aldo85a,
Kalle05a}. The guiding idea is that, as representations of $\Sset_\infty$ are connected to 
exchangeability of infinite random sequences, the representations of 
$\Bset_\infty$ should be connected to a new symmetry called 
`braidability' such that a noncommutative notion of conditional 
independence appears. Two important pillars for the realization of this 
idea are the noncommutative extended de Finetti theorem obtained by one 
of the authors \cite{Koes08aPP} and product representations of 
endomorphisms studied by the other author \cite{Gohm04a}.

Our main result is that `braidability' of infinite random sequences 
implies conditional independence and provides an interesting symmetry 
between exchangeability and spreadability; of course these notions 
are meant in a noncommutative sense. At first sight this might be 
surprising for a probabilist, since Ryll-Nardzewski showed that 
exchangeability and spreadability are equivalent for infinite random 
sequences \cite{RyNa57a}. But we will see that this equivalence fails 
in the noncommutative realm.

Hence we need to consider noncommutative random variables and 
stochastic processes and, because topological and analytical 
arguments are vital together with algebraic ones, we use a framework 
of von Neumann algebras. In particular this allows a rich theory of 
conditioning and independence. To get more directly to the heart of the matter without many preliminaries, we have collected in Appendix \ref{section:appendix} what we need about an operator algebraic noncommutative probability theory; and the reader may find it necessary to consult this appendix from time to time. Note however that for a first reading it is fine to
concentrate on tracial states and to avoid the Tomita-Takesaki theory;
in fact most of our examples in this paper are tracial. 

As an additional help for the reader and to streamline the flow of arguments around our main results some important information is demoted to the level of remarks. These remarks serve a number of purposes, from bringing together the different areas touched upon in this paper to providing background information for readers interested in future developments or open problems. 

In the following we comment on the most significant issues and describe briefly the contents of the paper.

Throughout this paper a \emph{probability space} $(\cA,\varphi)$ consists of 
a von Neumann algebra $\cA$ with separable predual and a faithful 
normal state $\varphi$ on $\cA$. A \emph{random variable} $\iota$ from
$(\cA_0,\varphi_0)$ to $(\cA,\varphi)$ is an injective *-homomorphism 
$\iota\colon \cA_0 \to \cA$ such that $\varphi_0 = \varphi \circ \iota$ 
and $\iota(\cA_0)$ embeds as a $\varphi$-conditioned von Neumann subalgebra 
of $\cA$ (see Definition \ref{def:nc-prob-space}). A \emph{random sequence} 
$\scrI$ is an infinite sequence of (identically distributed) random 
variables $\iota\equiv (\iota_n)_{n\in \Nset_0}$ from $(\cA_0, \varphi_0)$ 
to $(\cA,\varphi)$. We may assume $\cA_0 = \iota_0(\cA_0) \subset \cA$ and 
$\varphi_0 = \varphi|_{\cA_0}$ whenever it is convenient. If we restrict or 
enlarge (where possible) the domain $\cA_0$ of the random variables $\iota$ 
to another von Neumann algebra $\cC_0$ (with $\varphi$-conditioned embedding), 
then to simplify notation we will just write $\scrI_{\cC_0}$ instead of $\scrI$.

Our notion of (noncommutative) conditional independence actually emerges 
from the noncommutative extended de Finetti theorem \cite{Koes08aPP}. 
Consider the random sequence $\scrI_{\cC_0}$ which generates the von 
Neumann subalgebras $\cC_I$ with $I \subset \Nset_0$ and the tail algebra 
$\cC^\tail$: 
\[
\cC_I := \bigvee_{i \in I} \iota_i(\cC_0), \qquad \qquad
\cC^\tail:= \bigcap_{n \in \Nset_0}  \bigvee_{k \ge n}\iota_k(\cC_0).
\]
We say that $\scrI_{\cC_0}$ is (\emph{full/order}) \emph{$\cC^\tail$-independent} if the family $(\cC_I)_{I \subset \Nset_0}$ is (full/order) $\cC^\tail$-independent (in the sense of Definition \ref{def:order-independence}). Note that we do not require $\cC^\tail$ to be contained in $\cC_I$ since this would be far too restrictive in the context of conditioning and distributional symmetries. 

We next introduce several distributional symmetries on an intuitive 
level (see Section \ref{section:symmetries} for equivalent definitions 
which are less intuitive but more efficient in proofs).  Given the two 
random sequences $\scrI$ and $\scrtI$ with random variables 
$(\iota_n)_{n\ge 0}$ resp.~$(\tilde{\iota}_n)_{n\ge 0}$ 
from $(\cA_0,\varphi_0)$ to $(\cA,\varphi)$, we write
\[
(\iota_0, \iota_1, \iota_2, \ldots ) \stackrel{\distr}{=} 
(\tilde{\iota}_0, \tilde{\iota}_1, \tilde{\iota}_2, \ldots )
\]
if all their multilinear functionals coincide:
\[
\varphi\big(\iota_{\ii(1)}(a_1) \iota_{\ii(2)}(a_2) 
           \cdots \iota_{\ii(n)}(a_n)\big) 
= \varphi\big(\tilde{\iota}_{\ii(1)}(a_1) \tilde{\iota}_{\ii(2)}(a_2) 
              \cdots  \tilde{\iota}_{\ii(n)}(a_n)\big) 
\]
for all $n$-tuples $\ii\colon \{1, 2, \ldots, n\} \to \Nset_0$, 
$(a_1, \ldots, a_n) \in \cA_0^n$ and $n \in \Nset$. Now a random 
sequence $\scrI$ is said to be \emph{exchangeable} if its 
multilinear functionals are invariant under permutations: 
\[
 (\iota_0,\iota_1, \iota_2,  \ldots ) \stackrel{\distr}{=} 
 (\iota_{\pi(0)},\iota_{\pi(1)}, \iota_{\pi(2)},  \ldots ) 
\]
for any finite permutation $\pi \in \Sset_\infty$ of $\Nset_0$.  
We say that the random sequence $\scrI$ is \emph{spreadable} 
if every subsequence has the same multilinear functionals:
\[
(\iota_0,\iota_1, \iota_2, \ldots )  \stackrel{\distr}{=}  
(\iota_{n_0},\iota_{n_1}, \iota_{n_2}, \ldots ) 
\] 
for any subsequence $(n_0, n_1, n_2,\ldots)$ of $(0,1,2,\ldots)$. 
Finally, $\scrI$ is \emph{stationary} if the multilinear 
functionals are shift-invariant:
\[
 (\iota_0,\iota_1, \iota_2,  \ldots )  \stackrel{\distr}{=} 
 (\iota_{k},\iota_{k+1}, \iota_{k+2},\ldots ) 
\]  
for all $k\in \Nset$. 

The key definition of `braidability' is motivated from the 
characterization of exchangeability in Theorem \ref{thm:exchange-character}.  
\begin{Definition}\normalfont \label{def:braidability}
A random sequence $\scrI$ consisting of random variables 
$(\iota_n)_{n\in \Nset_0}$ from $(\cA_0, \varphi_0)$ to $(\cA,\varphi)$ 
is said to be \emph{braidable} if there exists a representation of the 
braid group, $\rho \colon \Bset_\infty \to \Aut{\cA,\varphi}$, such that 
the properties {\normalfont(PR)} and \normalfont{(L)} are satisfied:
\begin{align}
&&&&\iota_n &= \rho(\sigma_n \sigma_{n-1} \cdots \sigma_1)\iota_0 
&&\text{for all $n \ge 1$};&&&&\tag{PR} \label{eq:PR-braidable}\\
&&&& \iota_0    & = \rho(\sigma_n) \iota_0&& \text{if $n \ge 2$.}   
&&&&  \tag{L} \label{eq:L-braidable}
\end{align}
Here the $\sigma_i$'s are the Artin generators and $\Aut{\cA,\varphi}$ 
denotes the $\varphi$-preserving automorphisms of $\cA$. 
\end{Definition}
Throughout this paper we will make use of fixed point algebras of the braid group representation $\rho$. We denote by $\cA^{\rho(\sigma_n)}$ the fixed 
point algebra of $\rho(\sigma_n)$ in $\cA$, and by $\cA^{\rho(\Bset_\infty)}$ the fixed point algebra of $\rho(\Bset_\infty)$. 
 
Given $(\cA,\varphi)$ and a braid group representation $\rho \colon \Bset_\infty \to \Aut{\cA,\varphi}$ then a $\varphi$-conditioned 
von Neumann subalgebra $\cC_0$ of $\cA$ with the localization property 
\[
\cC_0 \subset \cA_0^{\rho}:=\bigcap_{n \ge 2} \cA^{\rho(\sigma_n)} 
\]
induces canonically a braidable random sequence $\scrI_{\cC_0}$ by $\iota_0 = \id|_{\cC_0}$ and \eqref{eq:PR-braidable}. The maximal choice is $\cC_0 = \cA_0^{\rho}$.

Our main result is a refinement of the noncommutative extended de Finetti theorem by inserting braidability into the scheme of distributional symmetries.
\begin{Theorem}\label{thm:main-1}
Consider the following assertions for the random sequence 
$\scrI_{\cC_0}$:
\begin{enumerate}
\item[(a)]
$\scrI_{\cC_0}$ is exchangeable;
\item[(b)]
$\scrI_{\cC_0}$ is braidable;
\item[(c)]
 $\scrI_{\cC_0}$ is spreadable; 
\item[(d)]
$\scrI_{\cC_0}$ is stationary and full $\cC^\tail$-independent;
\item[(d$_{\text{o}}$)]
$\scrI_{\cC_0}$ is stationary and order $\cC^\tail$-independent.
\end{enumerate}
Then we have the implications:
(a) $\Rightarrow$ (b) $\Rightarrow$ (c) $\Rightarrow$ (d) $\Rightarrow$ (d$_\text{o}$).
\end{Theorem}
Here $(a) \Rightarrow (c) \Rightarrow (d) \Rightarrow  (d_\text{o})$ 
is the extended de Finetti theorem \cite{Koes08aPP} which we review in Section 
\ref{section:symmetries} as Theorem \ref{thm:definetti}. We are left 
to prove $(a) \Rightarrow (b) \Rightarrow (c)$. The first implication 
is obtained below as a consequence of Theorem 
\ref{thm:exchange-character}, the second is shown in Theorem  
\ref{thm:sequences-braid}.

A fine point of Theorem \ref{thm:main-1} is that braidability is 
intermediate to distributional symmetries which are of \emph{purely 
probabilistic} nature.  Our main result also makes it clear that 
spreadability is not tied to the symmetric group in the noncommutative
realm since it can be obtained from the much wider context of braid 
groups. In fact we do not exclude the possibility of an even wider 
context where spreadability can be produced. 

An auxiliary result in the braid group context is that we can identify
 the tail algebra as the fixed point algebra of the representation. 
 This presents a braided extension of the Hewitt-Savage Zero-One Law.
 See Theorem \ref{thm:braided-HS}. In the last part of Section 
\ref{section:sequences-braid} we define fixed point algebras for 
certain subgroups in a Galois type manner and apply the noncommutative 
de Finetti theorem to obtain commuting squares. This sets the stage 
for the next step.

Our stochastic processes are stationary and this allows, as an 
additional tool, to introduce the time shift endomorphism. This is 
done in Section \ref{section:endomorphisms-braid}. Within the tower of 
the fixed point algebras this endomorphism can be written as an 
infinite product of automorphisms. It has been investigated by one of 
us (see \cite{Gohm04a}) how such product representations can be used 
to study operator theoretic and probabilistic structures, and the 
present setting is a particularly neat example for that. In Appendix \ref{section:appendix} we develop some refinements of this theory which put the setting of this paper into a wider context. In particular, from Theorem 
\ref{thm:tower-reconstruction} it becomes clear why we think of the 
tower of fixed point algebras as a `Galois type' structure and that 
it is unavoidable to choose exactly this tower to express the probabilistic structures associated to braids. In fact, if we have commuting squares for the braids as in Section \ref{section:sequences-braid}, then the associated tower is necessarily the tower of fixed point algebras. See Theorem \ref{thm:tower-reconstruction} and Corollary \ref{cor:tower-reconstruction}.

The results of Section \ref{section:sequences-braid} can be 
extended to obtain triangular towers of commuting squares from 
which we deduce that the time shift is a noncommutative 
Bernoulli shift (see Definition \ref{def:bernoulli}). 
\begin{Theorem}\label{thm:main-3}
The limit
\[
\alpha := \lim_{n \to \infty} 
          \rho(\sigma_1 \sigma_2 \cdots \sigma_{n-1}\sigma_n)
\]
exists on 
\[
\cA_\infty^\rho := \bigcup_{n \in \Nset} \bigcap_{k\ge n} 
\cA^{\rho(\sigma_{k})}
\] 
in the pointwise strong operator topology 
and defines an endomorphism of $\cA_\infty^\rho$ such that 
$\alpha^n \circ \iota_0 = \iota_n$. Then $\alpha$ restricted to 
$\bigvee_{n\in\Nset_0} \alpha^n(\cC_0 \vee \cA^{\rho(\Bset_\infty)})$ 
is a full Bernoulli shift over $\cA^{\rho(\Bset_\infty)}$ with 
generator $\cB_0 := \cC_0 \vee \cA^{\rho(\Bset_\infty)}$. Further 
we have
\[
\cA^\alpha = \cB^\tail = \cA^{\rho,\tail} = \cA^{\rho(\Bset_\infty)}.
\]
Here $\cA^\alpha$ is the fixed point algebra of $\alpha$ and
$\cB^\tail$ resp. $\cA^{\rho,\tail}$ are the tail algebras for
$\cB_0,\, \alpha(\cB_0),\, \alpha^2(\cB_0) \ldots$ resp.
$\cA^\rho_0,\, \alpha(\cA^\rho_0),\, \alpha^2(\cA^\rho_0) \ldots$
\end{Theorem}
This is proved after Theorem \ref{thm:endo-braid-i}. By modifying the 
product representations we also obtain examples of stationary random
 sequences which are order independent but not spreadable.

In Section \ref{section:presentation} we give a new presentation for 
the braid groups $\Bset_n$ in terms of the  generators 
$\set{\gamma_i}{1 \le i \le n-1}$ subject to the relations
\begin{eqnarray*}
\gamma_l \gamma_{l-1} (\gamma_{l-2}\cdots\gamma_k) \gamma_l 
=  \gamma_{l-1} (\gamma_{l-2}\cdots\gamma_k) \gamma_l \gamma_{l-1}  
\qquad \text{for $0 < k < l < n$}. 
\end{eqnarray*}

\begin{figure}[h]
\setlength{\unitlength}{0.2mm}
\begin{picture}(500,120)
\savebox{\artin}(20,20)[1]{\masterartin} 
\savebox{\artininv}(20,20)[1]{\masterartininv} 
\savebox{\strandr}(20,20)[1]{\masterstrandr} 
\savebox{\strandl}(20,20)[1]{\masterstrandl} 
\savebox{\horizontaldots}(20,20)[1]{\masterhorizontaldots}
\put(00,80){\usebox{\strandl}}  
\put(00,80){\usebox{\strandr}}  
\put(20,80){\usebox{\strandr}}
\put(40,80){\usebox{\strandr}}
\put(60,80){\usebox{\strandr}}
\put(00,60){\usebox{\strandl}}
\put(00,60){\usebox{\strandr}}  
\put(20,60){\usebox{\strandr}}
\put(40,60){\usebox{\strandr}}
\put(60,60){\usebox{\strandr}}
\put(00,40){\usebox{\artin}}  
\put(20,40){\usebox{\strandr}}
\put(40,40){\usebox{\strandr}}
\put(60,40){\usebox{\strandr}}
\put(80,40){\usebox{\horizontaldots}}
\put(00,20){\usebox{\strandl}}
\put(00,20){\usebox{\strandr}}  
\put(20,20){\usebox{\strandr}}
\put(40,20){\usebox{\strandr}}
\put(60,20){\usebox{\strandr}}
\put(00,0){\usebox{\strandl}}
\put(00,0){\usebox{\strandr}}  
\put(20,0){\usebox{\strandr}}
\put(40,0){\usebox{\strandr}}
\put(60,0){\usebox{\strandr}}
\put(200,80){\usebox{\strandl}}  
\put(200,80){\usebox{\strandr}}  
\put(220,80){\usebox{\strandr}}
\put(240,80){\usebox{\strandr}}
\put(260,80){\usebox{\strandr}}
\put(200,60){\usebox{\artin}}  
\put(220,60){\usebox{\strandr}}
\put(240,60){\usebox{\strandr}}
\put(260,60){\usebox{\strandr}}
\put(200,40){\usebox{\strandl}}  
\put(220,40){\usebox{\artin}}
\put(240,40){\usebox{\strandr}}
\put(260,40){\usebox{\strandr}}
\put(280,40){\usebox{\horizontaldots}}
\put(200,20){\usebox{\artininv}}  
\put(220,20){\usebox{\strandr}}
\put(240,20){\usebox{\strandr}}
\put(260,20){\usebox{\strandr}}
\put(200,0){\usebox{\strandl}}
\put(200,0){\usebox{\strandr}}  
\put(220,0){\usebox{\strandr}}
\put(240,0){\usebox{\strandr}}
\put(260,0){\usebox{\strandr}}
\put(400,80){\usebox{\artin}}  
\put(420,80){\usebox{\strandr}}
\put(440,80){\usebox{\strandr}}
\put(460,80){\usebox{\strandr}}
\put(400,60){\usebox{\strandl}}  
\put(420,60){\usebox{\artin}}
\put(440,60){\usebox{\strandr}}
\put(460,60){\usebox{\strandr}}
\put(400,40){\usebox{\strandl}}  
\put(420,40){\usebox{\strandl}}
\put(440,40){\usebox{\artin}}
\put(460,40){\usebox{\strandr}}
\put(480,40){\usebox{\horizontaldots}}
\put(400,20){\usebox{\strandl}}  
\put(420,20){\usebox{\artininv}}
\put(440,20){\usebox{\strandr}}
\put(460,20){\usebox{\strandr}}
\put(400,0){\usebox{\artininv}}  
\put(420,0){\usebox{\strandr}}
\put(440,0){\usebox{\strandr}}
\put(460,0){\usebox{\strandr}}
\end{picture}
\caption{Braid diagrams of $\gamma_1$, 
$\gamma_2$ and  $\gamma_3$ (left to right)}
\label{figure:squareroot}
\end{figure}
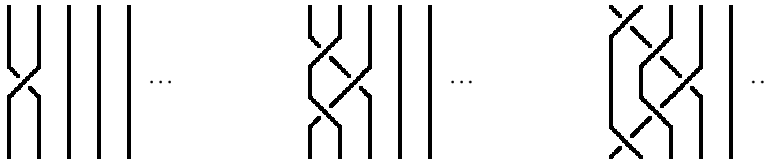

We shall see that the first three $\gamma_i$'s are depicted as the 
geometric braids shown in Figure \ref{figure:squareroot}.
This new presentation may be regarded to be intermediate between the 
Artin presentation \cite{Arti1925a} and the Birman-Ko-Lee presentation 
\cite{BKL98a}. We name it the `square root of free generator' 
presentation, since the squares $\set{\gamma_i^2}{1 \le i \le n-1}$ 
generate the free group $\Fset_{n-1}$. This is followed by the discussion 
of various shifts on $\Bset_\infty$ which are all cocycle perturbations 
of the shift in the Artin generators but which can be better understood 
by considering the new generators. We apply this to characterize relative 
conjugacy classes in $\Bset_\infty$.

These results are a preparation for Section 
\ref{section:left-regular-rep} in which we study in detail the left 
regular representation of $\Bset_\infty$ and some of the corresponding 
stochastic processes. The group von Neumann algebra is a 
non-hyperfinite $II_1$-factor and we can also control the fixed point 
algebras occurring in our theory. Using the Artin generators we get a 
random sequence which is conditionally independent but not spreadable. 
On the other hand the square roots of free generators turn out to be 
spreadable and we speculate about a braided extension of free 
probability theory which is suggested by this picture.

Section \ref{section:examples} discusses a few other examples of braid 
group representations such as the Gaussian representation, Hecke 
algebras and $R$-matrices. While these examples are well known it is 
interesting to reinterpret their properties in the context of 
braidability and our general theory. Of course the reader may take her 
favorite braid group representation and investigate what our theory is 
able to tell about it. For example it was tempting at this point to go 
straight into the braid group representations in Jones' subfactor 
theory. But we only give a short hint in this direction because we 
felt that this is a topic on its own which is better postponed to a 
consecutive paper.

\subsection*{Acknowledgment} The authors are grateful to the anonymous referee for several comments and suggestions helping us to improve the clearness of our presentation.   

\section{Distributional symmetries}
\label{section:symmetries}
A noncommutative version of de Finetti's theorem was obtained by one of 
us \cite{Koes08aPP}. We report here some of the obtained results, as 
far as they are needed for the present paper.

Distributional symmetries of random objects lead to deep structural
results in probability theory and the reader is referred to Kallenberg's 
monograph \cite{Kalle05a} for a recent account on this classical subject. 
Here we are interested to study some of these basic symmetries in the 
context of noncommutative random objects; and we will constrain ourselves 
to infinite sequences of noncommutative random variables (in the sense 
introduced above).     

We start with the introduction of some equivalence relations to prepare 
the definition of stationarity, spreadability and exchangeability in the 
broad sense of distributional symmetries. 
\begin{Notation}\normalfont
The group $\Sset_\infty$ is the inductive limit of the symmetric 
groups $\Sset_n$, $n \ge 2$, where $\Sset_n$ is generated on $\Nset_0$ 
by the transpositions $\pi_{i}\colon (i-1,i) \to (i,i-1)$ with 
$1 \le i <n$. We write $\pi_0$ for the identity of $\Sset_\infty$. 
By $[n]$ we denote the ordered set $\{1,2,\ldots,n\}$. 
\end{Notation}
The symmetric group $\Sset_\infty$ is presented by the 
transpositions $(\pi_i)_{i \in \Nset}$, subject to the relations  
\begin{align}
&&\pi_i \pi_{j} \pi_i &= \pi_{j} \pi_i \pi_{j} 
&\text{if $ \; \mid i-j \mid\, = 1 $;}
&& \tag{B1} \label{eq:B1-sg}\\
&&\pi_i \pi_j &= \pi_j \pi_i  
&\text{if $ \; \mid i-j \mid\, > 1 $;}&& \tag{B2} \label{eq:B2-sg}\\
&&\pi_i^2 &= \pi_0   &\text{for all $i \in \Nset$.}
 &&&\tag{S} \label{eq:S-sg}
\end{align}

\begin{Definition}\normalfont \label{def:relations}
Let $\ii,\jj\colon [n] \to \Nset_0$ be two $n$-tuples. 
\begin{enumerate}
\item
$\ii$ and $\jj$ are \emph{translation equivalent}, in symbols:
$\ii \sim_\theta \jj$, if there exists $k \in \Nset_0$ such that 
\[
\ii = \theta^k \circ \jj 
\qquad \text{or}\qquad\theta^k \circ \ii = \jj .
\]
Here denotes $\theta$ the right translation $m \mapsto m+1$ on 
$\Nset_0$. 
\item
$\ii$ and $\jj$ are \emph{order equivalent}, in symbols: 
$\ii \sim_o \jj$, if there exists a permutation $\pi \in \Sset_\infty$ 
such that
\[
\ii = \pi \circ \jj  
\qquad  \text{and} \qquad \pi|_{\jj([n])} 
\text{ is order preserving.} 
\] 
\item
$\ii$ and $\jj$ are \emph{symmetric equivalent}, in symbols: 
$\ii \sim_\pi \jj$, if there exists a permutation $\pi \in \Sset_\infty$ 
such that 
\[
\ii  =  \pi \circ \jj.  
\]  
\end{enumerate}
\end{Definition}
We have the implications $ (\ii \sim_\theta \jj) \Rightarrow 
(\ii \sim_o \jj) \Rightarrow (\ii \sim_\pi \jj)$.  
\begin{Remark}\normalfont
Order equivalence was introduced in the context of noncommutative 
probability in \cite{KoSp07a}. Our present formulation is equivalent 
to that given in \cite{KoSp07a}. 
\end{Remark}
For the notation of mixed higher moments of random variables,
it is convenient to use Speicher's notation of multilinear maps. 
\begin{Notation} \normalfont \label{notation:mlm}
Let $\iota\equiv (\iota_i)_{i \in \Nset_0} \colon (\cM_0, \psi_0) 
\to (\cM, \psi)$ be given. We put, for $\ii\colon [n] \to \Nset_0$, 
$\mathbf{a} = (a_1, \ldots, a_{n}) \in \cM_0^n$ and $n \in \Nset$, 
\begin{eqnarray}
\iota[\ii; \mathbf{a}] &:=& 
\iota_{\ii(1)}(a_1) \iota_{\ii(2)}(a_2) 
\cdots \iota_{\ii(n)}(a_{n}),  \label{eqn:iota1}\\
\psi_\iota[\ii; \mathbf{a}] &:=& 
\psi\big(\iota[\ii; \mathbf{a}]  \big) \label{eqn:iota2}. 
\end{eqnarray}
\end{Notation}
Next we define the distributional symmetries in terms of the mixed 
moments of a sequence of random variables. 
\begin{Definition}\label{def:ds}
A sequence of random variables $\iota\equiv (\iota_i)_{i \in \Nset_0} 
\colon (\cM_0, \psi_0) \to (\cM, \psi)$, also called 
\emph{(noncommutative) random sequence}, is 
\begin{enumerate}
\item
\emph{exchangeable} if, for any $n \in \Nset$,  
$\psi_\iota[\ii; \cdot\,] = \psi_\iota[\jj; \cdot\,]$ 
whenever $\ii \sim_\pi \jj$; 
\item
\emph{spreadable} if, for any $n \in \Nset$, 
$\psi_\iota[\ii; \cdot\,] = \psi_\iota[\jj; \cdot\,]$ 
whenever $\ii \sim_o \jj$;
\item
\emph{stationary} if, for any $n \in \Nset$,  
$\psi_\iota[\ii; \cdot\,] = \psi_\iota[\jj; \cdot\,]$ 
whenever $\ii \sim_\theta \jj$.       
\end{enumerate}
\end{Definition}
It is obvious from Definition \ref{def:relations} that we have the 
implications (i) $\Rightarrow$ (ii) $\Rightarrow$ (iii).  
\begin{Remark}\normalfont
`Spreadable' is also called `contractable' in the literature, for 
example in \cite{Kalle05a}, and is also in close contact with 
`subsymmetric' in Banach space theory \cite{JuPaXu07a}. 
Here the first notion is more suitable,
since `contractable' in the context of distributional 
symmetries should not be confused with the notion of a contraction in 
the context of operator theory. 
\end{Remark}
Roughly speaking, de Finetti's celebrated theorem states that 
exchangeable infinite commutative random sequences are mixed i.i.d. 
Inspired by \cite{Kalle05a}, a noncommutative dual version of this 
result has been obtained by one of the authors. 
\begin{Theorem}[\cite{Koes08aPP}]\label{thm:definetti}
Let $\scrI$ be a random sequence with (identically distributed) random 
variables 
\[
\iota\equiv (\iota_i)_{i \in \Nset_0} \colon 
(\cA_0, \varphi_0) \to (\cA, \varphi)
\]
and tail algebra 
\[
\cA^{\tail} := \bigcap_{n\ge 0} \bigvee_{k \ge n} \iota_k(\cA_0).
\] 
Consider the following conditions:
\begin{enumerate}
\item[(a)]
$\scrI$ is exchangeable;
\item[(b)]
$\scrI$ is spreadable; 
\item[(c)]
$\scrI$ is stationary and full $\cA^\tail$-independent; 
\item[(c$_\text{o}$)]
$\scrI$ is stationary and order $\cA^\tail$-independent. 
\end{enumerate}
Then we have the implications (a) $\Rightarrow$ (b) $\Rightarrow$ 
                              (c) $\Rightarrow$ (c$_\text{o}$).
\end{Theorem}
See Definitions \ref{def:independence} and \ref{def:order-independence}
for our concept of independence. 
If the von Neumann algebras considered are commutative, then one finds 
a dual version of the extended de Finetti theorem stated in 
\cite[Theorem 1.1]{Kalle05a}. Note that the implication 
(a) $\Rightarrow$ (b) is obvious from Definition \ref{def:ds}, so
is (c) $\Rightarrow$ (c$_{\text{o}}$) from Definition \ref{def:order-independence}. 
The implication (b) $\Rightarrow$ (c) is established by means from noncommutative 
ergodic theory. For the proof and a more-in-depth discussion of this result 
the reader is referred to \cite{Koes08aPP}. 

The following characterization of exchangeability motivates our 
notion of braidability as introduced in Definition 
\ref{def:braidability}. In fact, we have exactly the definition of 
braidability if in Theorem \ref{thm:exchange-character}(b) the 
symmetric group $\Sset_\infty$ is replaced by the braid group 
$\Bset_\infty$.
\begin{Definition}\normalfont
A random sequence $\scrI$ is said to be \emph{minimal} if 
$\cA = \bigvee_{n \ge 0} \iota_n(\cA_0)$.
\end{Definition}
\begin{Theorem}\label{thm:exchange-character}
The following are equivalent for a minimal random sequence $\scrI$:
\begin{enumerate}
\item[(a)] $\scrI$ is exchangeable;
\item[(b)] There exists a representation of the symmetric group,
$
  \rho\colon  \Sset_\infty \to \Aut{\cA,\varphi},
$
such that the properties {\normalfont(PR)} and \normalfont{(L)} are 
satisfied:
\begin{align}
&&&&\iota_n &= \rho(\pi_n \pi_{n-1} \cdots \pi_1 )\iota_0 
&&\text{for all $n \ge 1$};&&&&\tag{PR} \label{eq:PR}\\
&&&& \iota_0    & = \rho(\pi_n) \iota_0
&& \text{if $n \ge 2$.}   & &&&  \tag{L} \label{eq:L}
\end{align} 
\end{enumerate}
\end{Theorem}
\begin{proof}
`(a) $\Rightarrow$ (b)': The minimality of the random sequence ensures 
that the monomials $\iota[\ii; \mathbf{a}]$ with $n$-tuples 
$\ii \colon [n] \to \Nset_0$ and $\mathbf{a} \in (\cA_0)^n$,
$n \in \Nset$, are a weak* total set (see Notation \ref{notation:mlm}).
By exchangeability for every $\pi \in \Sset_\infty$
\[
\varphi_\iota[\ii; \mathbf{a}] 
= \varphi_\iota[\pi \circ \ii; \mathbf{a}],
\]
hence
\[
\rho(\pi)\colon 
\iota[\ii; \mathbf{a}] \mapsto \iota[\pi \circ \ii; \mathbf{a}]
\]
is well defined and extends to an element of $\Aut{\cA, \varphi}$.
Then $\rho \colon \Sset_\infty \to \Aut{\cA, \varphi}$ is a 
representation and properties {\normalfont(PR)} and \normalfont{(L)}
are easily verified.

`(b) $\Rightarrow$ (a)': 
Because $\Sset_\infty$ is generated by the $\pi_j$ it is enough to prove 
from (b) that for all $j\in\Nset,\,\mathbf{a} \in (\cA_0)^n$
\[
\varphi_\iota[\ii; \mathbf{a}] 
= \varphi_\iota[\pi_j \circ \ii; \mathbf{a}]. 
\]
As shown in Lemma \ref{lemma:braid2} (in the more general situation of 
braidability) $\rho(\pi_j) \iota_k = \iota_k$ if $k \notin \{j-1,j\}$.
Further $\rho(\pi_j)\iota_{j-1} = \iota_j$ (by {\normalfont(PR)})
and $\rho(\pi_j) \iota_j = \rho(\pi_j)^{-1} \iota_j= \iota_{j-1}$
(by {\normalfont(S)}: $\pi^2_j = \pi_0$). Summarizing, for all 
$k \in \Nset_0$
\[
\rho(\pi_j) \iota_k = \iota_{\pi_j(k)}.
\]
Using this together with $\rho(\pi_j) \in \Aut{\cA, \varphi}$
we obtain
\[
\varphi_\iota[\ii; \mathbf{a}] 
= \varphi(\iota_{i_1}(a_1) \cdots \iota_{i_n}(a_n))
= \varphi(\rho(\pi_j)\iota_{i_1}(a_1) \cdots \rho(\pi_j)\iota_{i_n}(a_n))
= \varphi_\iota[\pi_j \circ \ii; \mathbf{a}],
\]
which is what we wanted to prove.
\end{proof}
We give the proof of the first implication stated in our main result. 
\begin{proof}[Proof of Theorem \ref{thm:main-1} (a)\,$\Rightarrow$(b)]
We need to show that exchangeability implies braidability. 
Comparing the formulations of Theorem \ref{thm:exchange-character}(b) 
and Definition \ref{def:braidability}, this is accomplished by the
canonical epimorphism $\widehat{\phantom{x}}\colon \Bset_\infty \to 
\Sset_\infty$ satisfying $\widehat{\sigma_i} = \pi_i$ for all $i \in 
\Nset$. 
\end{proof}
Finally, we will need the following noncommutative generalization of the
Kolmogorov Zero-One Law. 
\begin{Theorem}[\cite{Koes08aPP}]\label{thm:kol-0-1}
Suppose the random sequence $\scrI$ is order $\cN$-independent 
with $\cN \subset \cA^\tail$. Then we have $\cN = \cA^\tail$. In
particular, an order $\Cset$-independent random sequence has a 
trivial tail algebra. 
\end{Theorem}
We will make use of this result within the proof of a `braided' 
noncommutative version of the Hewitt-Savage Zero-One Law (see Theorem 
\ref{thm:braid-independence}).
\section{Random variables generated by the braid group $\Bset_\infty$}
\label{section:sequences-braid}
This section is devoted to the construction of spreadable random sequences 
from braid group representations and the study of some of their properties.
Our results give an application for the noncommutative version of the de 
Finetti's theorem, Theorem \ref{thm:definetti}. Moreover this section 
provides the proof of Theorem \ref{thm:main-1}. 
 
Throughout this section, let $\rho\colon \Bset_\infty \to \Aut{\cA,\varphi}$ 
be a given representation on the probability space $(\cA,\varphi)$. For the 
construction of random sequences we are interested in the subgroups 
\[
\Bset_{n,\infty}:= \langle \sigma_k \mid  n \le k <\infty \rangle
\] 
of $\Bset_\infty$ and the corresponding fixed point algebras 
\[
\cA^{\rho(\Bset_{n,\infty})}
:= \set{x \in \cA}{\rho(\sigma)(x) = x 
\text{ for all } \sigma \in \Bset_{n,\infty}}.
\] 
These algebras provide us with a tower of von Neumann algebras:
\[
\cA^{\rho(\Bset_\infty)}  
= \cA^{\rho(\Bset_{1,\infty})} \subset \cA^{\rho(\Bset_{2,\infty})}
\subset \cdots   \subset   \cA^{\rho(\Bset_{n,\infty})} 
\subset \cdots \subset \cA_\infty^\rho 
:= \bigvee_{n\in \Nset} \cA^{\rho(\Bset_{n,\infty})}.  
\]
For short, we write $\cA_n^\rho := \cA^{\rho(\Bset_{n+2,\infty})}$ 
for $n \in \Nset$ so that the above tower can be written as
\[
\cA^{\rho(\Bset_\infty)} 
= \cA_{-1}^\rho \subset \cA_0^\rho \subset \cA_1^\rho 
\subset \cdots \subset \cA_{n-2}^\rho 
\subset \cdots \subset \cA_\infty^\rho.
\] 
These fixed point algebras give us a framework for the following 
construction of spreadable random sequences. We need some preparation. 
\begin{Lemma} \label{lemma:braid2}
Consider the braidable random sequence $\scrI_{\cA_0^\rho}$
and let\, $m,n \in \Nset_0$. If $n \notin \{m,m+1\}$, then we have
\[
\rho(\sigma_n) \iota_{m} =  \iota_{m}.
\]
\end{Lemma}
\begin{proof}
For $n=0$ this is trivial because $\rho(\sigma_0)$ is the identity.
We observe that $\cA_0^\rho = \cA^{\rho(\Bset_{2,\infty})} \subset 
\cA^{\rho(\sigma_n)}$ for $n\ge 2$. Because $\iota_m = \rho(\sigma_m 
\cdots \sigma_1 \sigma_0)_{|\cA_0^\rho}$, it is sufficient to show 
that
\begin{eqnarray*}
\sigma_n (\sigma_m \sigma_{m-1}\cdots \sigma_1 \sigma_0)
= \begin{cases}
(\sigma_m \sigma_{m-1}\cdots \sigma_1 \sigma_0) \sigma_{n} 
                                    & \text{if $n > m+1$}\\
(\sigma_m \sigma_{m-1}\cdots \sigma_1 \sigma_0) \sigma_{n+1} 
                                   & \text{if $0 \not= n < m$}
\end{cases}
\end{eqnarray*}
This is obvious for $n > m+1$ by \eqref{eq:B2}. The remaining case 
$0 \not= n < m$ follows from \eqref{eq:B1} and \eqref{eq:B2}:
\begin{eqnarray*}
\sigma_n (\sigma_{m} \sigma_{m-1}\cdots \sigma_{1}\sigma_{0})
&=& \sigma_{m} \sigma_{m-1} \cdots 
  \sigma_{n+2} \sigma_n \sigma_{n+1} \sigma_n \sigma_{n-1}\cdots
  \sigma_{1}\sigma_{0}\\
&=& \sigma_{m} \sigma_{m-1} \cdots 
  \sigma_{n+2} \sigma_{n+1} \sigma_{n} \sigma_{n+1} \sigma_{n-1}\cdots
  \sigma_{1}\sigma_{0} \\
&=& (\sigma_{m} \sigma_{m-1} \cdots
  \sigma_{1}\sigma_{0}) \sigma_{n+1}. 
\end{eqnarray*}
\end{proof}
As discussed after Definition \ref{def:braidability} it is useful to 
introduce some flexibility here by considering $\varphi$-conditioned
subalgebras $\cC_0 \subset \cA_0^\rho$.
\begin{Theorem}\label{thm:sequences-braid}
A braidable random sequence $\scrI_{\cC_0}$ is spreadable.
\end{Theorem}
\begin{proof}
Clearly it is enough to give the proof for the maximal case $\cC_0 = \cA_0^\rho$.
Using the multi-linear maps 
\[
\varphi_\iota[\ii; \cdot\,] \colon (\cA_0^\rho)^n \to \Cset
\]
(see Notation \ref{notation:mlm}), we need to show that, for any $n \in 
\Nset$, we have $\varphi_\iota[\ii; \cdot\,] = \varphi_\iota[\jj; \cdot\,]$ 
whenever $\ii \sim_o \jj$. Note that $\ii \sim_o \jj$ if and only if there 
exists a finite sequence of order-equivalent $n$-tuples $(\ii_k)_{k=1,\ldots,K} 
\colon [n] \to \Nset_0$ satisfying the following conditions:
\begin{enumerate}
\item
$\ii_1 = \ii$ and $\ii_K = \jj$;
\item
for each $k \in \{1,\ldots,K-1\}$, there exists a nonempty subset $A \subset [n]$ 
such that $\ii_{k}|_{[n]\backslash A} = \ii_{k+1}|_{{[n] \backslash A}}$, and 
such that $\ii_{k}(A) = \{l\}$ and  $\ii_{k+1}(A) = \{l^\prime\}$ with 
$|l-l^\prime|\le 1$ for some $l,l^\prime \in \Nset_0$. 
\end{enumerate}
Thus it is sufficient to prove that $\varphi_\iota[\ii_k; \cdot\,] = 
\varphi_\iota[\ii_{k+1}; \cdot\,]$. Let $\ii_{k}$ and $\ii_{k+1}$ be two 
$n$-tuples meeting the above conditions for some set $A$ and nonnegative
integers $l$ and $l^\prime$. The case $l=l^\prime$ is trivial. It is sufficient 
to consider the case $l^\prime = l+1$ (otherwise reverse the order). We note 
$l+1 \notin \ii_k([n])$ for later purposes. Since $\varphi = \varphi \circ 
\rho(\sigma_{l+1})$, we obtain that, for some $n$-tuple $\mathbf{a} = 
(a_1, \ldots,a_n ) \in (\cA_0^\rho)^n$,
\begin{eqnarray*}
\varphi_\iota[\ii_k; \mathbf{a}] 
&=& \varphi\big(\iota_{\ii_{k}(1)}(a_1) \cdots \iota_{\ii_{k}(n)}(a_n)\big)\\
&=& \varphi\big( \rho(\sigma_{l+1})
\iota_{\ii_{k}(1)}(a_1)\cdots\rho(\sigma_{l+1}) \iota_{\ii_{k}(n)}(a_n)\big) 
\end{eqnarray*}
We consider each factor $\rho(\sigma_{l+1})\iota_{\ii_{k}(j)}(a_j)$ 
separately, for fixed $j\in [n]$. Since $l+1 \notin \ii_k([n])$, one 
of the following two cases occurs:\\
Case $l+1 \notin \{\ii_{k}(j),\ii_{k}(j)+1\}$: we conclude $j \notin A$ 
and $\rho(\sigma_{l+1})\iota_{\ii_k(j)} = \iota_{\ii_k(j)}=
\iota_{\ii_{k+1}(j)}$ with Lemma \ref{lemma:braid2}.\\
Case $l+1 = \ii_{k}(j)+1$:  we infer $j \in A$ and 
$\rho(\sigma_{l+1})\iota_{\ii_k(j)} = \rho(\sigma_{l+1})\iota_{l}
= \iota_{l+1} = \iota_{\ii_k(j)+1} = \iota_{\ii_{k+1}(j)}$ from the 
definition of the random variable $\iota_l$ and the relation between 
$\ii_k$ and $\ii_{k+1}$.\\
Altogether, we conclude that
\begin{eqnarray*}
\rho(\sigma_{l+1})\iota_{\ii_k(1)}(a_1) \cdots  
   \rho(\sigma_{l+1}) \iota_{\ii_k(n)}(a_n) 
= \iota_{\ii_{k+1}(1)}(a_1) \cdots   \iota_{\ii_{k+1}(n)}(a_n)
\end{eqnarray*} 
and thus $\varphi_\iota[\ii_k; \mathbf{a}] = \varphi_\iota[\ii_{k+1}; 
\mathbf{a}]$. Now a finite induction on $k \in \{1,\ldots,K\}$ shows 
that $\varphi_\iota[\ii; \mathbf{a}] = \varphi_\iota[\ii_{1}; 
\mathbf{a}] = \cdots = \varphi_\iota[\ii_K; \mathbf{a}] = 
\varphi_\iota[\jj; \mathbf{a}]$.
\end{proof}
\begin{Remark}\normalfont \label{rem:implementation}
The proof actually shows that for order equivalent tuples $\ii \sim_o \jj$ 
there always exists a braid $\tau \in \Bset_\infty$ such that $\rho(\tau) 
(\iota[\ii; a]) = \iota[\jj; a]$ for all $a \in (\cA_0^\rho)^n$. Note 
however that while it is always possible to construct a representation of 
$\Sset_\infty$ from an exchangeable sequence it is not clear at the present 
state which additional probabilistic conditions would allow us to construct 
braid group representations and braidability from spreadable sequences.
\end{Remark}
We are now in the position to apply the noncommutative extended de Finetti 
theorem \ref{thm:definetti}.
\begin{Theorem}\label{thm:braid-independence}
A braidable random sequence $\scrI_{\cC_0}$ is stationary and full 
$\cC^\tail$-independent.
\end{Theorem}
\begin{proof}
The random sequence $\scrI_{\cC_0}$ is spreadable by Theorem 
\ref{thm:sequences-braid}. Thus its stationarity and full 
$\cC^\tail$-independence follows directly from the implication 
(b) $\Rightarrow$ (c) of Theorem \ref{thm:definetti}.
\end{proof}
Another immediate implication of Theorem \ref{thm:sequences-braid} and the 
noncommutative version of de Finetti's theorem, Theorem \ref{thm:definetti}, 
is a noncommutative generalized version of the famous Hewitt-Savage Zero-One 
Law. More precisely, in the context of exchangeable commutative infinite 
random sequences and representations of the symmetric group $\Sset_\infty$, 
the tail algebra of the random sequence is identified as the fixed point 
algebra of $\Sset_\infty$ (see \cite{Kalle05a}, for example). Now the 
Hewitt-Savage Zero-One Law states that these two algebras are trivial if the 
random sequence is (order) $\Cset$-independent. With Theorem 
\ref{thm:sequences-braid} at our disposal, the tail algebra 
\[
\cC^{\tail} := \bigcap_{n\ge 0}\bigvee_{k \ge n} \iota_k(\cC_0),
\] 
is identified in the much broader context of braid group representations and 
we obtain a `braided' extension of the Hewitt-Savage Zero-One Law.  
\begin{Theorem}\label{thm:braided-HS} 
A braidable random sequence $\scrI_{\cC_0}$ satisfies
\[
\cC^{\tail} \subset \cA^{\rho(\Bset_\infty)}.
\]
Suppose $\cA^{\rho(\Bset_\infty)} \subset \cC_0 \,\,(\subset \cA_0^\rho)$.
Then we have an equality
\[
\cC^{\tail} = \cA^{\rho(\Bset_\infty)}.
\]
In particular, these two algebras are trivial if the random sequence 
$\scrI_{\cC_0}$ is order $\Cset$-independent. 
\end{Theorem}
Note that the assumption $\cA^{\rho(\Bset_\infty)} \subset \cC_0$ is superfluous for the maximal choice $\cC_0 = \cA_0^\rho$.
\begin{proof}
To prove $\cC^{\tail} \subset \cA^{\rho(\Bset_\infty)}$ it 
suffices to show that $\cC^{\tail} \subset \cA^{\rho(\sigma_l)}$ for any 
$l \in \Nset$. But $\cC^{\tail} \subset \bigvee_{k > l} 
\iota_k(\cC_0)$, hence the assertion follows from
$ \rho(\sigma_l)\iota_{k} =\iota_{k}$ for all $k > l$, by Lemma 
\ref{lemma:braid2}. 

Now assume $\cA^{\rho(\Bset_\infty)} \subset \cC_0$.
We verify $ \cA^{\rho(\Bset_\infty)} \subset 
\cC^{\tail}$. Indeed, because 
$\cA^{\rho(\Bset_\infty)} \subset \cC_0 \subset
\cA^{\rho(\Bset_{2,\infty})} = \cA_0^\rho$, we have $\cA^{\rho(\Bset_\infty)} 
= \iota_k(\cA^{\rho(\Bset_\infty)}) \subset \iota_k(\cC_0)$ for all 
$k$. This implies  $\cA^{\rho(\Bset_\infty)}\subset 
\bigcap_{n\ge 0} \bigvee_{k \ge n} \iota_k(\cC_0) = \cC^{\tail}$.   

Finally,  we conclude $\Cset \simeq  \cC^{\tail} = 
\cA^{\rho(\Bset_\infty)}$ from the $\Cset$-independence 
of $\iota$ by applying Theorem \ref{thm:kol-0-1}. 
\end{proof}
\begin{Remark}\normalfont
The assumptions do \emph{not} require the global $\rho(\Bset_2)$-invariance 
of the von Neumann algebra $\iota_0(\cC_0) \vee \iota_1(\cC_0) = 
\cC_0 \vee \rho(\sigma_1)(\cC_0)$. This invariance property is automatic if 
the representation $\rho$ is a representation of the symmetric group 
$\Sset_\infty$, or in other words, if we have $\rho(\sigma_n^2)= \id$ for 
all $n \in \Nset$:
\[
\rho(\sigma_1)\Big(\cC_0 \vee \rho(\sigma_1)(\cC_0)\Big) 
= \rho(\sigma_1)(\cC_0) \vee \rho(\sigma_1^2)(\cC_0)
= \rho(\sigma_1)(\cC_0) \vee \cC_0.
\] 
\end{Remark}
Our next result states that from braid group representations we can produce 
commuting squares (see Appendix \ref{section:appendix}) in the tower of fixed point algebras.
\begin{Theorem}\label{thm:braid-tower}
Assume that the probability space $(\cA,\varphi)$ is equipped with the  
representation $\rho\colon \Bset_\infty \to \Aut{\cA,\varphi}$ and let 
$\cA_{n-1}^\rho:= \cA^{\rho(\Bset_{n+1,\infty})}$, the fixed point
algebra of $\rho(\Bset_{n+1,\infty})$ (with $n \in \Nset_0$).  
Then 
\[
\cA_{-1}^\rho \subset \cA_{0}^\rho 
\subset \cA_{1}^\rho \subset \cdots \subset \cA
\] 
is a tower of von Neumann algebras such that, for all $n \in \Nset_0$,
$\rho(\sigma_{n+1})$ restricts to an automorphism of $\cA_{n+1}^\rho$ and
\begin{eqnarray*}
\begin{matrix}
\rho(\sigma_{n+1})(\cA_{n}^\rho) &\subset & \cA_{n+1}^\rho\\
  \cup    &         &  \cup\\
\cA_{n-1}^\rho & \subset &   \cA_{n}^\rho
\end{matrix}
\end{eqnarray*}
is a commuting square.
\end{Theorem}
\begin{proof}
The global invariance of $\cA_{n+1}^\rho$ under the action of 
$\rho(\sigma_{n+1})$ is concluded from $\cA_{n+1}^\rho = 
\cA^{\rho(\Bset_{n+3,\infty})}$ and relation \eqref{eq:B2}. 
The existence of the conditional expectations needed to define 
a commuting square follows along the lines sketched in Appendix \ref{section:appendix}. 

We claim for $n=0$ the order $\cA_{-1}^\rho$-independence of $\cA_0^\rho$ 
and $\rho(\sigma_1)\cA_0^\rho$. Indeed, the corresponding random sequence 
$\scrI_{\cA_0^\rho}$ enjoys $\iota_0(\cA_0^\rho)=\cA_0^\rho$ and 
$\iota_1(\cA_0^\rho)= \rho(\sigma_1)\cA_0^\rho$. Thus we deduce from
Theorem \ref{thm:braid-independence} combined with Theorem 
\ref{thm:braided-HS} that $\cA_0^\rho$ and $\rho(\sigma_1)\cA_0^\rho$ are 
$\cA_{-1}^\rho$-independent. This establishes our claim. 
 
To treat the general case we consider, for $n \ge 0$ fixed, the 
`$n$-shifted' representation $\rho_n$ of $\Bset_{\infty}$, defined by the 
multiplicative extension of  
\[
\rho_n(\sigma_k):= 
\begin{cases}
\rho(\sigma_{n+k}) & \text{if $k >0$}\\
\rho(\sigma_{0})=\id   & \text{if $k=0$}
\end{cases}.
\]
Thus we have $\cA^{\rho_n(\Bset_\infty)}= \cA_{-1}^{\rho_n}=
\cA_{n-1}^\rho$. Now we obtain the $\cA_{n-1}^\rho$-independence of 
$\cA_{n}^\rho$ and $\rho(\sigma_{n+1}) (\cA_{n}^\rho)$ along the same 
lines of arguments as before, based on the `$n$-shifted' random sequence 
$\iota^{(n)} \equiv (\iota^{(n)}_{k})_{k\ge 0}\colon 
(\cA_{n}^\rho,\varphi|_{\cA_{n}^\rho}) \to (\cA,\varphi)$ with 
$\iota^{(n)}_{k}= \rho_n(\sigma_k \sigma_{k-1}\cdots  \sigma_1\sigma_0)$. 
\end{proof}
\begin{Remark}\normalfont \label{rem:rs-shifted}
The proof of the above theorem uses that, given the braid group 
representation $\rho$, one can easily produce `$n$-shifted' representations 
$\rho_n := \rho \circ \sh^n$, where the shift $\sh\colon \Bset_\infty \to 
\Bset_\infty$ sends the generator $\sigma_i$ to $\sigma_{i+1}$. Note that 
this endomorphism $\sh$ is injective (see \cite[Lemma 3.3]{Deho00a}). The 
prospect of passing to a shifted representation is of interest 
when the fixed point algebra $\cA^{\rho(\Bset_{2,\infty})}$ turns out to 
be trivial or too small for the required task, and a fixed point algebra 
$\cA^{\rho(\Bset_{n+2,\infty})}$, with $n \in \Nset$, fulfills the requirements.
This idea is used in Section \ref{section:left-regular-rep}, Theorem
\ref{thm:lrr-shifted}.
\end{Remark}
\begin{Remark}\normalfont \label{rem:rs-inverse} 
The spreadable random sequence $\iota$ from Theorem 
\ref{thm:sequences-braid} is induced by positive braids of the form 
$\sigma_n \sigma_{n-1}\cdots \sigma_1$ for $n \in \Nset$. Using the group 
automorphism $\inv \colon \Bset_\infty \to \Bset_\infty$ which sends the 
generator $\sigma_i$ to $\sigma_{i}^{-1}$ for all $i \in \Nset$ and given 
the representation $\rho$, we obtain the representation $\rho^{\inv} := 
\rho \circ \inv$. Thus the random variables
\[
\iota_n^{\inv}= \rho(\sigma_n^{-1} \sigma_{n-1}^{-1}\cdots 
\sigma_{1}^{-1} \sigma_{0}^{-1})_{|_{\cA_0^\rho}}, \qquad n \in \Nset_0,
\] 
define another spreadable random sequence, since 
$\cA^{\rho^{\inv}(\Bset_{2,\infty})} = \cA^{\rho(\Bset_{2,\infty})}$ 
and consequently Theorem \ref{thm:sequences-braid} applies. The random 
sequences $\iota$ and $\iota^{\inv}$ have the same tail algebra. This 
is easily concluded from Theorem \ref{thm:braided-HS} and 
$\cA^{\rho^{\inv}(\Bset_\infty)} = \cA^{\rho(\Bset_\infty)}$. More 
generally, we have $\cA^{\rho^{\inv}(\Bset_{n,\infty})} = 
\cA^{\rho(\Bset_{n,\infty})}$ for all $n \in \Nset$. Thus the commuting 
squares constructed in Theorem \ref{thm:braid-tower} from starting with 
the representation $\rho^{\inv}$ are just those coming from the representation 
$\rho$, but now with $\rho(\sigma_n^{-1})$ in the upper left corner, 
instead of $\rho(\sigma_n)$. 
\end{Remark}
\section{Endomorphisms generated by the braid group $\Bset_\infty$}
\label{section:endomorphisms-braid}
Suppose $(\cA,\varphi)$ is equipped with the representation 
$\rho\colon \Bset_\infty \to \Aut{\cA,\varphi}$. Then the fixed point algebras 
$\cA_{n-2}^\rho:= \cA^{\rho(\Bset_{n,\infty})}$, with $n \ge 1$ and 
$\Bset_{n,\infty}=\langle \sigma_{n}, \sigma_{n+1},\ldots \rangle$, provide a tower 
\[
 \cA^{\rho(\Bset_\infty)} 
= \cA_{-1}^\rho \subset \cA_{0}^\rho \subset \cA_{1}^\rho \subset \cdots 
                                      \subset \cA_\infty^\rho \subset \cA,
\]
where $\cA_\infty^\rho$ denotes the weak closure of 
$\cA^{\rho, \alg}_\infty:=\bigcup_{k} \cA_k^\rho$. 
\begin{Definition}\normalfont \label{def:generating}
The representation $\rho\colon \Bset_\infty \to (\cA, \varphi)$ has the 
\emph{generating property} if 
\[
\cA = \cA_\infty^\rho.
\]
\end{Definition}
If a representation of $\Bset_\infty$ has the generating property then $\cA$ is also generated by the fixed point algebras $\cA^{\rho(\sigma_n)}$. This generating property is not always fulfilled from the outset (see Proposition \ref{prop:non-generating}), but we may always restrict a representation to a generating one.
\begin{Proposition}
The representation $\rho: \Bset_\infty \to \Aut{\cA,\varphi}$ restricts to 
the generating representation $\rho^{\operatorname{res}}\colon \Bset_\infty 
\to \Aut{\cA^\rho_\infty, \varphi^\rho_\infty}$ such that 
$\rho(\sigma_i)(\cA_\infty^\rho) \subset \cA_\infty^\rho$ and 
$E_{\cA_\infty^\rho} E_{\cA^{\rho(\sigma_i)}} = E_{\cA^{\rho(\sigma_i)}}E_{\cA_\infty^\rho}$ 
(for all $i \in \Nset$). 
\end{Proposition}
\begin{proof}
We study the action of $\rho(\sigma_i)$ on $\cA_n^\rho$ and 
$\cA_\infty^\rho$. Since $\cA_n^\rho =  \bigcap_{k \ge n}
\cA^{\rho(\sigma_{k+2})}$, it holds $ \rho(\sigma_i)(\cA_n^\rho) \subset 
\cA_n^\rho$ for $i \in \{1, \ldots, n \} \cup \{n+2, n+3, \ldots \}$. If 
$i = n+1$, then $\cA_n^\rho \subset \cA_{n+1}^\rho$ implies 
$\rho(\sigma_{n+1})\cA_n^\rho \subset \rho(\sigma_{n+1})\cA_{n+1}^\rho 
\subset \cA_{n+1}^\rho$. From this we conclude that $\rho(\sigma_i) 
\cA_n^\rho \subset \cA_{n+1}^\rho$ and therefore $\rho(\sigma_i) 
(\cA_\infty^\rho) \subset \cA_{\infty}^\rho$ for all $i \in \Nset$. A 
similar argument ensures the inclusion $\rho(\sigma_i^{-1}) (\cA_\infty^\rho) 
\subset \cA_{\infty}^{\rho}$ for all $i \in \Nset$. Consequently,
$\cA_\infty^\rho$ is globally invariant under the action of 
$\rho(\Bset_\infty)$ and the representation $\rho: \Bset_\infty \to 
\Aut{\cA,\varphi}$ restricts to the representation $\rho^{\operatorname{res}}
\colon \Bset_\infty \to \Aut{\cA^\rho_\infty, \varphi_\infty^\rho}$ which, 
by construction, has the generating property.

That $E_{\cA_\infty^\rho}$ and $E_{\cA^{\rho(\sigma_i)}}$ commute is concluded 
by routine arguments from 
\[
E_{\cA_\infty^\rho} \rho(\sigma_i) E_{\cA_\infty^\rho} 
= \rho(\sigma_i) E_{\cA_\infty^\rho}
\]
and thus $E_{\cA_\infty^\rho} \rho(\sigma_i) =  \rho(\sigma_i)  
E_{\cA_\infty^\rho}$, and by an application of the mean ergodic theorem
(as in \cite[Theorem 8.3]{Koes08aPP}, for example):
\[ 
E_{\cA^{\rho(\sigma_i)}} 
= \lim_{N\to \infty }\frac{1}{N} \sum_{n=0}^{N-1} \rho(\sigma_i^n).
\]
Here the limit is taken in the pointwise strong operator topology. 
\end{proof}

The following proposition gives a method to construct new braid group representations from a given (simpler) one. In this way we can find many interesting examples with and without the generating property.        
\begin{Proposition}\label{prop:non-generating}
Given the representation $\rho\colon \Bset_\infty \to \Aut{\cA,\varphi}$  suppose $\gamma \in \Aut{\cA,\varphi}$ is an automorphism commuting with all $\rho(\sigma_i)$'s. Then the multiplicative extension of 
\[
\rho_\gamma(\sigma_i):=\begin{cases} 
                       \gamma \rho(\sigma_i) & \text{if $i > 0$}\\
                       \id & \text{if $i =0$}
                       \end{cases}
\] 
defines another representation of $\Bset_\infty$ in $\Aut{\cA,\varphi}$ such that:
\begin{enumerate}
\item[(i)]
the restriction $\rho^{\operatorname{res}}_\gamma$ of $\rho_\gamma$ to the fixed point algebra $\cA^\gamma$ has the generating property if $\rho$ has the generating property, and $\rho^{\operatorname{res}}_\gamma$ coincides with the restriction of $\rho$ to $\cA^\gamma$;
\item[(ii)] 
the representation $\rho_{\gamma}$ does not have the generating property  
if all $\rho(\sigma_i)$'s are $N$-periodic but $\gamma^N \neq \id$ for some $N \in \Nset$.
\end{enumerate}  
\end{Proposition}

\begin{proof}
An elementary calculation shows that $\rho_\gamma$ satisfies the
braid relations
\begin{align*}
\rho_\gamma(\sigma_i)\rho_\gamma(\sigma_j)\rho_\gamma(\sigma_i)
&= \rho_\gamma(\sigma_j)\rho_\gamma(\sigma_i)\rho_\gamma(\sigma_j)
& \text{for $|i-j|=1$};\\
\rho_\gamma(\sigma_i)\rho_\gamma(\sigma_j)
&= \rho_\gamma(\sigma_j)\rho_\gamma(\sigma_i)
& \text{for $|i-j|>1$}.
\end{align*}
Clearly, $\varphi \circ \rho_\gamma(\sigma)= \varphi$. So $\rho_\gamma$ is a representation from $\Bset_\infty$ into $\Aut{\cA,\varphi}$. \\
(i) Since $\gamma$ commutes with all $\rho(\sigma_i)$'s and $\rho_\gamma(\sigma_i)$'s, both representations $\rho$ and $\rho_\gamma$ restrict to $\cA^\gamma$. An elementary calculation shows that these two restrictions coincide; and we denote them both by $\rho_\gamma^{\operatorname{res}}$. We show next that $\cA_k^{\rho_\gamma^{\operatorname{res}}} = \cA_k^\rho \cap \cA^\gamma$.
For this purpose let $E_k$ and $E_\gamma$ be the $\varphi$-preserving conditional expectations from $\cA$ onto $\cA_k^{\rho} = \cA^{\rho(\Bset_{k+2,\infty})}$ resp.~$\cA^\gamma$. Since all $\rho(\sigma_i)$'s and $\gamma$ commute, we conclude that $\gamma E_k = E_k \gamma$ for all $k$.  But this entails $E_k E_\gamma = E_\gamma E_k$ by an application of the mean ergodic theorem  (similar as done for Proposition 3.2). 
Consequently, $\cA_k^{\rho_\gamma^{\operatorname{res}}} = \cA_k^\rho \cap \cA^\gamma$. Finally, the generating property of $\rho$ implies that $\lim_{k \to \infty} E_k = \id$ (in the pointwise $\sot$-sense). Thus $\lim_{k \to \infty} E_\gamma E_k = E_\gamma$. So  $\rho_\gamma^{\operatorname{res}}$ has the generating property.\\
(ii)
We infer from the $N$-periodicity of the $\rho(\sigma_i)$'s that
$
\rho_\gamma(\sigma_i)^N = \gamma^N. 
$  
Consequently,
\[
 \cA^{\rho_\gamma(\sigma_i)} \subset \cA^{\rho_\gamma(\sigma_i)^N} = \cA^{\gamma^N}
\]
and, passing to the intersections of the fixed point algebras  $\cA^{\rho_\gamma(\sigma_i)}$, 
\[
\cA^{\rho_\gamma}_k = \cA^{\rho_\gamma(\Bset_{k+2, \infty})} \subset \cA^{\gamma^N}.
\]
Altogether this gives the inclusions $ \cA^{\rho_\gamma}_\infty  \subset \cA^{\gamma^N} \subset \cA$. So $\rho_\gamma$ does not have the generating property if $\gamma$ is not $N$-periodic.   
\end{proof}
Let us consider some concrete examples.
\begin{Example}\normalfont
Period $N=1$ means that $\rho$ is trivial, i.e.~$\rho(\sigma_i) = \id$ for all $i$. Now any non-trivial $\gamma \in \Aut{\cA,\varphi}$ gives a representation $\rho_\gamma$ without the generating property. 
The simplest example is a (classical) probability space with two points,
each with probability $\frac{1}{2}$, on which $\gamma$ acts by interchanging the two points. Here we find $\cA_\infty^{\rho_\gamma} \simeq \Cset \neq \Cset^2 = \cA$.   
\end{Example}
\begin{Example}\normalfont
The case of period $N=2$ covers the representations of $\Sset_\infty$ in $\Aut{\cA,\varphi}$ and to obtain a non-generating representation one just needs to find a non-idempotent $\gamma \in \Aut{\cA,\varphi}$ which commutes with all $\rho(\sigma_i)$'s. 

Interesting examples come from infinite tensor products $\bigotimes_\Nset M_2$ with product states. The canonical tensor product flips on neighboring factors provide us with a state-preserving representation $\rho$ of $\Sset_\infty$, and thus of $\Bset_\infty$ with period $N=2$. Now implement the automorphism $\gamma$ as a Xerox action, in other words: as the infinite tensor product $\gamma = \bigotimes_{\Nset} \gamma_0$, where $\gamma_0$ is a state-preserving automorphism of $M_2$. It is easy to check that $\gamma$ commutes with all $\rho(\sigma_i)$'s.
Since $\gamma^2 \neq \id$ if and only if $\gamma_0^2 \neq \id$, we have plenty of choices for $\gamma_0$ such that $\rho_\gamma$ does not have the generating property. On the other hand $\rho$ itself clearly has the generating property and by Proposition \ref{prop:non-generating}(i) this is inherited by the restriction $\rho_\gamma^{\operatorname{res}}$ to the fixed point algebra of $\gamma$.
\end{Example}
\begin{Example}\normalfont
In a non-tracial situation we always have a non-trivial modular automorphism group commuting with $\rho$ (see Appendix \ref{section:appendix}) which gives further possibilities to apply
Proposition \ref{prop:non-generating}.
\end{Example}

From now on we will assume without loss of generality that the 
representation $\rho$ has the generating property:
\[
\cA = \cA_\infty^\rho.
\]
This allows us to define the following endomorphism. 
Due to the fixed point properties of the tower with respect to the 
$\rho(\sigma_k)$'s and the weak denseness of $\cA^{\rho,\alg}_\infty$, 
it is easily verified that 
\begin{align}
\alpha&:= \sotlim_{n\to \infty}\rho(\sigma_1\sigma_2 \cdots \sigma_n)
\tag{PR-0} \label{item:bpr-0}
\end{align}
exists pointwise in $\cA$ and defines an adapted endomorphism $\alpha$ 
of $\cA$ with a product representation. This is discussed in detail in 
Appendix \ref{section:appendix}, see Definition \ref{def:adapted-end}. In particular 
(for $k,n \in \Nset_0$)
\begin{align}
\rho(\sigma_k) (\cA_n^\rho) &= \cA_n^\rho \quad\quad (k \le n) \tag{PR-1}  \label{item:bpr-i}\\
\rho(\sigma_{k})|_{\cA_n^\rho} & = \id |_{\cA_n^\rho} 
\quad (k \ge n+2) \tag{PR-2} \label{item:bpr-ii}
\end{align} 
Next we address how the endomorphism $\alpha$ relates to the spreadable 
(and thus stationary) random sequence $\iota\equiv (\iota_n)_{n\in\Nset_0}$ 
from Theorem \ref{thm:sequences-braid}. For this purpose we need an 
elementary result on the Artin generators of the braid group.
\begin{Lemma}\label{lemma:br}
If $\sigma_1,\ldots,\sigma_m$ are Artin generators of the braid group $\Bset_{m+1}$ 
then
\[
\sigma_1 \sigma_2 \cdots \sigma_{m-1} \sigma_m \sigma_{m-1} \cdots 
\sigma_2  \sigma_1\;=\;
\sigma_m \sigma_{m-1} \cdots \sigma_2 \sigma_1 \sigma_2  \cdots
\sigma_{m-1} \sigma_m.
\]
In words: Pyramids up and down are the same.
\end{Lemma}
\begin{proof}
For $m=2$ this is \eqref{eq:B1}. The general case follows by induction.
\end{proof}
\begin{Proposition}\label{prop:endo-rs} 
The endomorphism $\alpha$ for $(\cA,\varphi)$, given by
\begin{eqnarray*}
\alpha&= \lim_{n\to \infty}\rho(\sigma_1\sigma_2 \cdots \sigma_n),
\end{eqnarray*}
and the random sequence $\iota\equiv (\iota_n)_{n\in \Nset_0} \colon (\cA_0^\rho,\varphi_0^\rho) \to (\cA,\varphi)$,
given by 
\begin{eqnarray*}
\iota_n = \rho(\sigma_n \sigma_{n-1}\cdots \sigma_1 \sigma_0)|_{\cA_0^\rho},
\end{eqnarray*}
are related by
\begin{eqnarray}
\alpha^n|_{\cA_0^\rho} &=& \iota_n . \label{eq:endo-rs} 
\end{eqnarray}
for all $n\in \Nset_0$. (Above we have put $\varphi_0^\rho = \varphi|_{\cA_0^\rho}$.) 
\end{Proposition}
\begin{proof}
This is trivial for $n=0$, since 
$\iota_0 = \rho(\sigma_0)$ and $\sigma_0$ is the identity in $\Bset_\infty$.
Using induction and Lemma \ref{lemma:br}, we get for $x \in \cA_0^\rho= \cA^{\rho(\Bset_{2,\infty})}$ 
\begin{eqnarray*}
\alpha^{n+1}(x) &=& \alpha \alpha^{n}(x) = \alpha \, \iota_n (x)\\
&=& \rho(\sigma_1 \sigma_2 \cdots \sigma_{n+1}) \rho(\sigma_n \sigma_{n-1} \cdots \sigma_0) (x)\\
&=& \rho(\sigma_{n+1} \sigma_n \cdots \sigma_{2}\sigma_{1} \sigma_{2} \cdots \sigma_n \sigma_{n+1})(x)\\
&=& \rho(\sigma_{n+1} \sigma_n \cdots \sigma_{1}\sigma_{0})(x)\\
&=& \iota_{n+1}(x).
\end{eqnarray*}
\end{proof}
We can interpret Proposition 
\ref{prop:endo-rs} by saying that $\alpha$ implements the time evolution 
of the stationary process associated to the random sequence $\iota$.
Note that, even with the generating property of the representation, 
the minimal part
\[
\bigvee_{n \in \Nset_0} \alpha^n(\cA_0^\rho) 
= \bigvee_{n \in \Nset_0} \iota_n(\cA_0^\rho)
\]
may be strictly contained in $\cA$. 
\begin{Theorem} \label{thm:endo-braid-i}
Assume that the probability space $(\cA,\varphi)$ is equipped with the 
generating representation $\rho\colon \Bset_\infty \to \Aut{\cA,\varphi}$ 
and let $\cA_{n-1}^\rho := \cA^{\rho(\Bset_{n+1,\infty})}$, with $n\in \Nset_0$.
Then one obtains a triangular tower of inclusions such that each cell forms 
a commuting square: 
\begin{eqnarray*}
\setcounter{MaxMatrixCols}{20}
\begin{matrix}
\cA_{-1}^\rho &\subset&   \cA_0^\rho & \subset & \cA_1^\rho & \subset &\cA_2^\rho & \subset &\cA_3^\rho & \subset  & \cdots & \subset & \cA\\
        &&          \cup  &         & \cup  &         & \cup &         & \cup  &       & & & \cup  \\
        &&   \cA_{-1}^\rho&\subset&\alpha(\cA_0^\rho)&\subset&\alpha(\cA_1^\rho)&\subset&\alpha(\cA_2^\rho)&\subset& \cdots & \subset & \alpha(\cA)\\
         &&               &         & \cup  &         & \cup   &         & \cup  &       & & & \cup \\
              &&&&   \cA_{-1}^\rho& \subset & \alpha^2 (\cA_0^\rho)  & \subset & \alpha^2(\cA_1^\rho)
              & \subset  & \cdots & \subset & \alpha^2(\cA) \\
       &&&&&&   \cup           &         & \cup  &       &&  & \cup \\
        &&&&&&   \vdots           &         & \vdots  &       &&  & \vdots 
\end{matrix}
\setcounter{MaxMatrixCols}{10}
\end{eqnarray*}
\end{Theorem}
\begin{proof}[Proof of Theorem \ref{thm:endo-braid-i}]
All inclusions stated in the triangular tower follow from the adaptedness 
property $\alpha(\cA_{n-1}^\rho) \subset \cA_n^\rho$ (for all n) of the 
endomorphism $\alpha$, which is an immediate consequence of (PR-1) and 
(PR-2) (see also Appendix \ref{section:appendix}). We are left to prove that all its cells are 
commuting squares. We know already from Theorem \ref{thm:braid-tower} that, 
for any $n \ge 1$, 
\begin{eqnarray*}
\begin{matrix}
\cA_{n-1}^\rho &\subset &  \cA_n^\rho \\
  \cup    &         &  \cup\\
\cA_{n-2}^\rho & \subset &   \rho(\sigma_n)(\cA_{n-1}^\rho)  
\end{matrix}
\end{eqnarray*}
is a commuting square. Introducing the automorphism 
$\gamma_n := \rho(\sigma_1 \cdots \sigma_{n})$,
\begin{eqnarray*}
\begin{matrix}
\gamma_{n-1}(\cA_{n-1}^\rho) &\subset &  \gamma_{n-1}(\cA_n^\rho) \\
  \cup    &         &  \cup\\
\gamma_{n-1}(\cA_{n-2}^\rho) & \subset &  \gamma_{n-1}\,\rho(\sigma_n)(\cA_{n-1}^\rho)  
\end{matrix}
\end{eqnarray*}
is obviously a commuting square. Let us consider counterclockwise the 
corners of this diagram, starting with the lower left corner. One 
readily verifies: 
\begin{eqnarray*}
\text{(\ref{item:bpr-0}) \& (\ref{item:bpr-ii})}
\quad &\Rightarrow& \quad \gamma_{n-1}(\cA_{n-2}^\rho)=\alpha(\cA_{n-2}^\rho);\\
\text{(\ref{item:bpr-0}) \& (\ref{item:bpr-ii})} 
\quad &\Rightarrow& \quad                  \gamma_{n-1}\rho(\sigma_n)\,(\cA_{n-1}^\rho)=\alpha(\cA_{n-1}^\rho);\\
\text{(\ref{item:bpr-i})}  \quad &\Rightarrow& \quad 
 \gamma_{n-1}(\cA_{n}^\rho)=\gamma_n \, \rho(\sigma_n)^{-1}(\cA_{n}^\rho) 
  = \gamma_n (\cA_{n}^\rho) = \cA_n^\rho;\\
\text{(\ref{item:bpr-i})} \quad &\Rightarrow &\quad   \gamma_{n-1}(\cA_{n-1}^\rho)=\cA_{n-1}^\rho. 
\end{eqnarray*}
Summarizing this corner discussion, we have shown that, for any $n \ge 1$,
\begin{eqnarray*}
\begin{matrix}
\cA_{n-1}^\rho &\subset &  \cA_n^\rho \\
  \cup    &         &  \cup\\
\alpha(\cA_{n-2}^\rho) & \subset &   \alpha(\cA_{n-1}^\rho)  
\end{matrix}
\end{eqnarray*}
is a commuting square. Inductively acting with $\alpha$ on the corners 
of such a commuting square, we conclude further for any $k,n \in \Nset$:
\begin{eqnarray*}
\begin{matrix}
\alpha^{k-1}(\cA_{n-1}^\rho) &\subset &  \alpha^{k-1}(\cA_n^\rho) \\
  \cup    &         &  \cup\\
\alpha^{k}(\cA_{n-2}^\rho) & \subset &   \alpha^{k}(\cA_{n-1}^\rho)  
\end{matrix}
\end{eqnarray*}
is a commuting square. But this is the general form of a cell in the 
triangular tower of inclusions. Cells involving the column on the right are 
also commuting squares by the generating property of the representation.
\end{proof}
\begin{proof}[Proof of Theorem \ref{thm:main-3}]
Here we consider a $\varphi$-conditioned subalgebra $\cC_0$ 
of $\cA^\rho_0$. If $\cA^\rho_{-1} = \cA^{\rho(\Bset_\infty)}$
is contained in $\cC_0$ then, with the results of Theorem
\ref{thm:braid-tower} or \ref{thm:endo-braid-i}, we can apply
Theorem \ref{thm:endo-pr} to obtain a Bernoulli shift over
$\cA^\rho_{-1}$ with generator $\cC_0$. If $\cA^\rho_{-1}$
is not contained in $\cC_0$ then we have to use 
$\cB_0 = \cC_0 \vee \cA^\rho_{-1}$ as generator. Now in
$\cA^\alpha = \cB^\tail = \cA^{\rho,\tail} = \cA^{\rho(\Bset_\infty)}$
the first two equalities follow from Theorem \ref{thm:endo-pr}
while the last equality is Theorem \ref{thm:braided-HS}.
\end{proof}
\begin{Remark} \normalfont \label{rem:galois}
Given a braid group representation Theorem \ref{thm:endo-braid-i}
shows that we find commuting squares in the tower of fixed point algebras 
and, with Theorem \ref{thm:endo-pr} and Theorem \ref{thm:main-3}, a 
corresponding Bernoulli shift. In Theorem 
\ref{thm:tower-reconstruction} we show that we can reconstruct the tower of 
fixed point algebras if we are given the commuting squares. This is what we 
mean when we think of this tower as a `Galois type' tower.
\end{Remark}
Next we produce a more general family of examples for adapted endomorphisms 
with product representations.  
\begin{Corollary}\label{cor:endo-braid-ii}
Under the assumptions of Theorem \ref{thm:endo-braid-i}, let a sequence 
$\mathbf{\varepsilon} = (\varepsilon_k)_{k\in \Nset} \in \{1, -1\}^\Nset$ 
be given. Then 
\begin{align}
\alpha_{\mathbf{\varepsilon}} (x)&:= \sotlim_{n\to \infty} \rho(\sigma_1^{\varepsilon_1} \sigma_2^{\varepsilon_2} \cdots \sigma_n^{\varepsilon_n})(x)\tag{PR-0'} \label{item:bpr-0-e}
\end{align}
defines an endomorphism for $(\cA,\varphi)$ such that, for all $k,n \ge 0$, 
\begin{align}
\rho(\sigma_k^{\varepsilon_k}) (\cA_n^\rho) &= \cA_n^\rho 
\quad\quad (k \le n); 
\tag{PR-1'}  \label{item:bpr-i-e}\\
\rho(\sigma_k^{\varepsilon_k})|_{\cA_n^\rho} & = \id |_{\cA_n^\rho} 
\quad (k \ge n+2).
\tag{PR-2'} \label{item:bpr-ii-e}
\end{align} 
Moreover, one obtains a family of triangular towers of commuting squares, 
indexed by the sequence $\varepsilon$:  
\begin{eqnarray*}
\setcounter{MaxMatrixCols}{20}
\begin{matrix}
\cA_{-1}^\rho &\subset&   \cA_0^\rho & \subset & \cA_1^\rho & \subset &\cA_2^\rho & \subset &\cA_3^\rho & \subset  & \cdots & \subset & \cA\\
        &&          \cup  &         & \cup  &         & \cup &         & \cup  &       & & & \cup  \\
        &&  \alpha_\varepsilon(\cA_{-1}^\rho)&\subset&\alpha_\varepsilon(\cA_0^\rho)&\subset&\alpha_\varepsilon(\cA_1^\rho)&\subset&\alpha_\varepsilon(\cA_2^\rho)&\subset& \cdots & \subset & \alpha_\varepsilon(\cA)\\
         &&               &         & \cup  &         & \cup   &         & \cup  &       & & & \cup \\
              &&&&   \vdots &  & \vdots  &     & \vdots
              &  &   &   & \vdots 
 \end{matrix}
\setcounter{MaxMatrixCols}{10}
\end{eqnarray*}
\end{Corollary}
\begin{proof} 
We can argue in the same way as for Theorem \ref{thm:endo-braid-i}.
In particular, to prove that each cell forms a commuting square, we take 
advantage of the following observation:
\begin{eqnarray*}
\begin{matrix}
\cA_{n-1}^\rho &\subset &  \cA_n^\rho \\
  \cup    &         &  \cup\\
\cA_{n-2}^\rho & \subset &   \rho(\sigma_n)(\cA_{n-1}^\rho)  
\end{matrix}
\end{eqnarray*}
is a commuting square if and only if 
\begin{eqnarray*}
\begin{matrix}
\cA_{n-1}^\rho &\subset &  \cA_n^\rho \\
  \cup    &         &  \cup\\
\cA_{n-2}^\rho & \subset &   \rho(\sigma_n^{-1})(\cA_{n-1}^\rho)  
\end{matrix}
\end{eqnarray*}
is a commuting square. 
\end{proof}
Corollary \ref{cor:endo-braid-ii} provides a rich source of stationary order 
$\cA^{\rho(\Bset_\infty)}$-independent random sequences 
(by Theorem \ref{thm:endo-pr})
which come from braid 
group representations and in general are no longer spreadable. Thus one obtains an interesting class of random sequences for which the implication 
(c$_\text{o}$) $\Rightarrow$ (b) in the noncommutative de Finetti theorem, Theorem \ref{thm:definetti}, fails to be true.  

\begin{Remark}\normalfont\label{rem:garside}
Equation \eqref{eq:endo-rs} can be expressed in terms of fundamental braids.
This suggests that Garside structures may be relevant for possible generalizations of braidability. We close this section with some information on this connection for the interested reader.

For a solution of the word and conjugacy problem for braid groups
Garside introduced in \cite{Gars69a} the fundamental braid 
\[
\Delta_n := (\sigma_{1} \sigma_2 \cdots \sigma_{n-1}) (\sigma_{1} \sigma_2 \cdots \sigma_{n-2})
             \cdots (\sigma_1 \sigma_2)(\sigma_1)
\]
in $\Bset_n$. Since $\Delta_n^2$ generates the center $Z(\Bset_n)$ of $\Bset_n$, 
the fundamental braid $\Delta_n$ is also called the `square root' of the center. 
In a recent new approach to the word and conjugacy problem, Birman et.\ al.\ 
introduced in \cite{BKL98a} the fundamental braid
\[
\delta_n:= \sigma_{n-1}\sigma_{n-2}\cdots \sigma_1,
\]
which satisfies $\delta_n^n = \Delta_n^2$ and thus can be thought of to be the 
`$n$th root' of the center. Since 
\[
\rho(\Delta_n)|_{\cA_0^\rho} = \alpha^{n-1}|_{\cA_0^\rho} \quad\text{and}\quad
\rho(\delta_n)|_{\cA_0^\rho} = \iota_{n-1},
\]  
the relationship between the random sequence $\iota$ and the endomorphism $\alpha$,
\[
\alpha^{n-1}|_{\cA_0^\rho} = \iota_{n-1},
\]
reads in terms of the fundamental braids $\Delta_n$ and $\delta_n$ as     
\[
\rho(\Delta_n)|_{\cA_0^\rho}= \rho(\delta_n)|_{\cA_0^\rho}.
\]
Furthermore, we have 
\begin{eqnarray}\label{eq:delta-Delta}
\rho(\Delta_n^2)\rho^{\inv} (\delta_n)|_{\cA_0^\rho}= \rho(\delta_n)|_{\cA_0^\rho},
\end{eqnarray}
which defines also a spreadable random sequence according to Remark 
\ref{rem:rs-inverse}. Indeed, it follows from the definition of 
$\Delta_n$ and Lemma \ref{lemma:braid2} that
\[
\rho(\Delta_n^2)|_{\cA_0^\rho}=  \rho(\Delta_n) \rho(\delta_n)|_{\cA_0^\rho} = \rho(\sigma_1\sigma_2\cdots \sigma_{n-1})\rho(\delta_n)|_{\cA_0^\rho}.
\]
Next we multiply with $\rho^{\inv} (\delta_n) = \rho(\inv(\delta_n))$ 
from the left. We note that 
\[
\inv(\delta_n) = \sigma_{n-1}^{-1}\sigma_{n-2}^{-1} \cdots \sigma_1^{-1}
\] 
and use that $\Delta_n^2$ generates the
center of $\Bset_n$. This gives \eqref{eq:delta-Delta}. 

We also conclude from the above discussion that the random 
sequences $\iota$ and $\iota^{\inv} \equiv (\iota^{\inv}_n)_{n \in \Nset_0}$ 
(from Remark \ref{rem:rs-inverse}) are connected according to
\[
\rho(\Delta_{n+1}^2)\iota^{\inv}_{n} = \iota_{n} \quad \text{for all $n \in \Nset$.}
\]
Finally, we point out to the reader that the Artin presentation and the Birman-Ko-Lee presentation are presently the only known presentations of $\Bset_n$ which possess a Garside structure \cite{Bir08a,Deho08a}. 
\end{Remark}
\section{Another braid group presentation,\\ $k$-shifts and braid handles}
\label{section:presentation}
We begin by providing a presentation for braid groups in terms of generators 
$\gamma_i$ which will play a distinguished role within our investigations of 
spreadability for random sequences coming from braid group representations. 
This presentation may be regarded as being `intermediate' between the Artin 
presentation \cite{Arti1925a} and the Birman-Ko-Lee presentation \cite{BKL98a}. 
\begin{Theorem}[Square root of free generators presentation]  
\label{thm:sqrt-rep} 
The braid group $\Bset_n$ (for $n \ge 3$) is presented by the generators 
$\set{\gamma_i}{1 \le i \le n-1}$ subject to the defining relations
\begin{align}\tag{EB} \label{eq:EB} 
       \gamma_{l} \gamma_{l-1} (\gamma_{l-2} \gamma_{l-3} \cdots \gamma_{k+1}\gamma_{k}) \gamma_l 
     = \gamma_{l-1}  (\gamma_{l-2} \gamma_{l-3}\cdots \gamma_{k+1}\gamma_{k}) \gamma_{l}\gamma_{l-1} 
\end{align}
for $0 < k < l <n$.
\end{Theorem}
It will become apparent in the proof that the two sets of generators 
$\{\gamma_i\}$ and $\{\sigma_i\}$ are related by the formulas:
\begin{eqnarray} 
\gamma_k^{} &=& (\sigma_{1}^{}\cdots \sigma_{k-1}^{})\sigma_{k}^{}(\sigma_{k-1}^{-1}\cdots \sigma_{1}^{-1}); \label{eq:gamma2sigma}\\
\sigma_k^{} &=& (\gamma_{1}^{-1}\cdots \gamma_{k-1}^{-1}) \gamma_k^{}(\gamma_{k-1}^{} \cdots \gamma_1^{}). \label{eq:sigma2gamma} 
\end{eqnarray}
\begin{figure}[h]
\setlength{\unitlength}{0.2mm}
\begin{picture}(410,130)
\savebox{\artin}(20,20)[1]{\masterartin} 
\savebox{\artininv}(20,20)[1]{\masterartininv} 
\savebox{\strandr}(20,20)[1]{\masterstrandr} 
\savebox{\strandl}(20,20)[1]{\masterstrandl} 
\savebox{\horizontaldots}(20,20)[1]{\masterhorizontaldots}
\put(0,80){\usebox{\artin}}  
\put(20,80){\usebox{\strandr}}
\put(40,80){\usebox{\strandr}}
\put(60,80){\usebox{\strandr}}
\put(80,80){\usebox{\strandr}}
\put(100,80){\usebox{\strandr}}
\put(0,60){\usebox{\strandl}}
\put(20,60){\usebox{\artin}}
\put(40,60){\usebox{\strandr}}
\put(60,60){\usebox{\strandr}}
\put(80,60){\usebox{\strandr}}
\put(100,60){\usebox{\strandr}}
\put(0,40){\usebox{\strandl}}
\put(20,40){\usebox{\strandl}}
\put(40,40){\usebox{\artin}}
\put(60,40){\usebox{\strandr}}
\put(80,40){\usebox{\strandr}}
\put(100,40){\usebox{\strandr}}
\put(120,40){\usebox{\horizontaldots}}
\put(0,20){\usebox{\strandl}}
\put(20,20){\usebox{\artininv}}
\put(40,20){\usebox{\strandr}}
\put(60,20){\usebox{\strandr}}
\put(80,20){\usebox{\strandr}}
\put(100,20){\usebox{\strandr}}
\put(0,0){\usebox{\artininv}}
\put(20,0){\usebox{\strandr}}
\put(40,0){\usebox{\strandr}}
\put(60,0){\usebox{\strandr}}
\put(80,0){\usebox{\strandr}}
\put(100,0){\usebox{\strandr}}
\put(300,90){\usebox{\artin}}  
\put(320,90){\usebox{\strandr}}
\put(340,90){\usebox{\strandr}}
\put(360,90){\usebox{\strandr}}
\put(380,90){\usebox{\strandr}}
\put(400,90){\usebox{\strandr}}
\put(300,70){\usebox{\strandl}}
\put(320,70){\usebox{\artin}}
\put(340,70){\usebox{\strandr}}
\put(360,70){\usebox{\strandr}}
\put(380,70){\usebox{\strandr}}
\put(400,70){\usebox{\strandr}}
\put(300,50){\usebox{\strandl}}
\put(320,50){\usebox{\strandl}}
\put(340,50){\usebox{\artin}}
\put(360,50){\usebox{\strandr}}
\put(380,50){\usebox{\strandr}}
\put(400,50){\usebox{\strandr}}
\put(420,50){\usebox{\horizontaldots}}
\put(300,30){\usebox{\strandl}}
\put(320,30){\usebox{\strandl}}
\put(340,30){\usebox{\artin}}
\put(360,30){\usebox{\strandr}}
\put(380,30){\usebox{\strandr}}
\put(400,30){\usebox{\strandr}}
\put(300,10){\usebox{\strandl}}
\put(320,10){\usebox{\artininv}}
\put(340,10){\usebox{\strandr}}
\put(360,10){\usebox{\strandr}}
\put(380,10){\usebox{\strandr}}
\put(400,10){\usebox{\strandr}}
\put(300,-10){\usebox{\artininv}}
\put(320,-10){\usebox{\strandr}}
\put(340,-10){\usebox{\strandr}}
\put(360,-10){\usebox{\strandr}}
\put(380,-10){\usebox{\strandr}}
\put(400,-10){\usebox{\strandr}}
\end{picture}
\caption{ 
Braid diagrams of $\gamma_3 = \sigma_1 \sigma_2 \sigma_3 \sigma^{-1}_2 \sigma^{-1}_1$ (left) and
$\gamma_3^2 = \sigma_1 \sigma_2 \sigma_3^2 \sigma^{-1}_2 \sigma^{-1}_1$ (right)}
\end{figure}
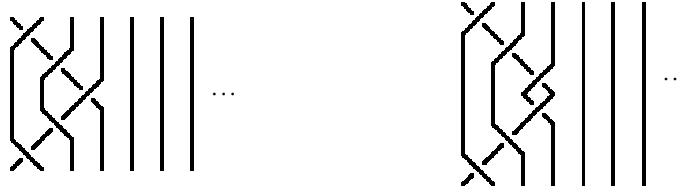\label{figure:generators}

The geometrical picture will be further discussed later. 
We note that $\gamma_1^2, \ldots, \gamma_{n-1}^2$ generate 
the free group $\Fset_{n-1}$ (cf.~\cite[Section 1.4]{Bir75a} 
and \cite[Lecture 5]{Jone91a}). This connection to the free 
group motivated us to call this presentation as titled in 
Theorem \ref{thm:sqrt-rep}. Similarly, we will say that the $\gamma_i$'s 
are \emph{square roots of free generators}. 
\begin{proof}
We start with the relations $\eqref{eq:B1}$ and $\eqref{eq:B2}$ for the 
Artin generators $\sigma_i$ and suppose that the $\gamma_j$'s are defined 
according to \eqref{eq:gamma2sigma}. Then a straightforward computation
yields 
\begin{eqnarray*}
\gamma_{l}^{} \gamma_{l-1}^{} \gamma_{l-2}^{} \gamma_{l-3}^{} \cdots \gamma_{k+1}^{}\gamma_{k}^{} 
= \sigma_1^{} \sigma_2^{}  \cdots \sigma_{l-1}^{}\sigma_l^{} \sigma_{k-1}^{-1} \sigma_{k-2}^{-1} \cdots \sigma_2^{-1}\sigma_{1}^{-1}. 
\end{eqnarray*}
We multiply the equation above with $\gamma_l$ from the right and see that 
the left hand side of \eqref{eq:EB} equals
\begin{align}\label{eq:srbraid-i}
(\sigma_1^{} \sigma_2^{}  \cdots \sigma_{l-1}^{}\sigma_l^{}) (\sigma_{k}^{} \sigma_{k+1}^{}\sigma_{k+2}^{}\cdots  \sigma_{l-2}^{}\sigma_{l-1}^{}\sigma_l^{}) (\sigma_{l-1}^{-1} \sigma_{l-2}^{-1}\cdots \sigma_{1}^{-1}). 
\end{align}
Similarly, we obtain for the right hand side of \eqref{eq:EB} the expression
\begin{align}\label{eq:srbraid-ii}
(\sigma_1^{} \sigma_2^{}  \cdots \sigma_{l-1}^{}) (\sigma_{k}^{} \sigma_{k+1}^{}\sigma_{k+2}^{}\cdots  \sigma_{l-2}^{}\sigma_{l-1}^{}\sigma_l^{}) (\sigma_{l-2}^{-1} \sigma_{l-3}^{-1}\cdots \sigma_{1}^{-1}). 
\end{align} 
Thus the formulas \eqref{eq:srbraid-i} and \eqref{eq:srbraid-ii} are equal 
if and only if 
\[
\sigma_l^{} (\sigma_{k}^{} \sigma_{k+1}^{}\sigma_{k+2}^{}\cdots  \sigma_{l-2}^{}\sigma_{l-1}^{}\sigma_l^{}) \sigma_{l-1}^{-1}
= (\sigma_{k}^{} \sigma_{k+1}^{}\sigma_{k+2}^{}\cdots  \sigma_{l-2}^{}\sigma_{l-1}^{}\sigma_l^{}).
\]
Indeed, an application of \eqref{eq:B1} and \eqref{eq:B2} to the left hand side 
of this equation shows 
\begin{eqnarray*}
\sigma_l^{} (\sigma_{k}^{} \sigma_{k+1}^{}\sigma_{k+2}^{}\cdots  \sigma_{l-2}^{}\sigma_{l-1}^{}\sigma_l^{}) \sigma_{l-1}^{-1}
&=&   (\sigma_{k}^{} \sigma_{k+1}^{}\sigma_{k+2}^{}\cdots  \sigma_{l-2}^{}\sigma_l^{} \sigma_{l-1}^{}\sigma_l^{}) \sigma_{l-1}^{-1}\\
&=&  (\sigma_{k}^{} \sigma_{k+1}^{}\sigma_{k+2}^{}\cdots  \sigma_{l-2}^{}\sigma_{l-1}^{} \sigma_{l}^{}\sigma_{l-1}^{}) \sigma_{l-1}^{-1}.
\end{eqnarray*}
Thus the relations  \eqref{eq:B1} and \eqref{eq:B2} for the $\sigma_i$'s imply 
the relations \eqref{eq:EB} for the $\gamma_i's$. 
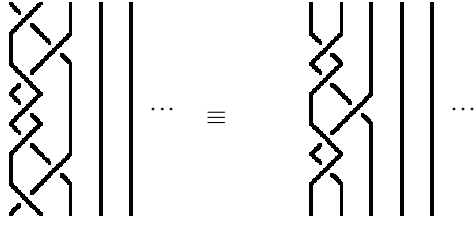
\begin{figure}[h]
\setlength{\unitlength}{0.2mm}
\begin{picture}(280,140)
\savebox{\artin}(20,20)[1]{\masterartin} 
\savebox{\artininv}(20,20)[1]{\masterartininv} 
\savebox{\strandr}(20,20)[1]{\masterstrandr} 
\savebox{\strandl}(20,20)[1]{\masterstrandl} 
\savebox{\horizontaldots}(20,20)[1]{\masterhorizontaldots}
\put(0,120){\usebox{\artin}}  
\put(20,120){\usebox{\strandr}}
\put(40,120){\usebox{\strandr}}
\put(60,120){\usebox{\strandr}}
\put(0,100){\usebox{\strandl}}
\put(20,100){\usebox{\artin}}
\put(40,100){\usebox{\strandr}}
\put(60,100){\usebox{\strandr}}
\put(0,80){\usebox{\artininv}}
\put(20,80){\usebox{\strandr}}
\put(40,80){\usebox{\strandr}}
\put(60,80){\usebox{\strandr}}
\put(0,60){\usebox{\artin}}
\put(20,60){\usebox{\strandr}}
\put(40,60){\usebox{\strandr}}
\put(60,60){\usebox{\strandr}}
\put(80,60){\usebox{\horizontaldots}}
\put(0,40){\usebox{\artin}}  
\put(20,40){\usebox{\strandr}}
\put(40,40){\usebox{\strandr}}
\put(60,40){\usebox{\strandr}}
\put(0,20){\usebox{\strandl}}
\put(20,20){\usebox{\artin}}
\put(40,20){\usebox{\strandr}}
\put(60,20){\usebox{\strandr}}
\put(0,0){\usebox{\artininv}}
\put(20,0){\usebox{\strandr}}
\put(40,0){\usebox{\strandr}}
\put(60,0){\usebox{\strandr}}
\put(130,60){$\equiv$}
\put(200,120){\usebox{\strandl}}
\put(200,120){\usebox{\strandr}}
\put(220,120){\usebox{\strandr}}
\put(240,120){\usebox{\strandr}}
\put(260,120){\usebox{\strandr}}
\put(200,100){\usebox{\artin}}
\put(220,100){\usebox{\strandr}}
\put(240,100){\usebox{\strandr}}
\put(260,100){\usebox{\strandr}}
\put(200,80){\usebox{\artin}}
\put(220,80){\usebox{\strandr}}
\put(240,80){\usebox{\strandr}}
\put(260,80){\usebox{\strandr}}
\put(200,60){\usebox{\strandl}}
\put(220,60){\usebox{\artin}}
\put(240,60){\usebox{\strandr}}
\put(260,60){\usebox{\strandr}}
\put(280,60){\usebox{\horizontaldots}}
\put(200,40){\usebox{\artininv}}
\put(220,40){\usebox{\strandr}}
\put(240,40){\usebox{\strandr}}
\put(260,40){\usebox{\strandr}}
\put(200,20){\usebox{\artin}}
\put(220,20){\usebox{\strandr}}
\put(240,20){\usebox{\strandr}}
\put(260,20){\usebox{\strandr}}
\put(200,00){\usebox{\strandl}}
\put(200,00){\usebox{\strandr}}
\put(220,00){\usebox{\strandr}}
\put(240,00){\usebox{\strandr}}
\put(260,00){\usebox{\strandr}}
\end{picture}
\caption{Braid relation  $\gamma_2 \gamma_1 \gamma_2 = \gamma_1 \gamma_2 \gamma_1$}
\label{figure:braidrelation-i}
\end{figure}

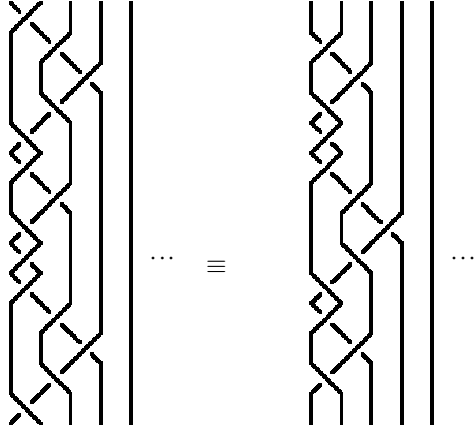
\begin{figure}[h]
\setlength{\unitlength}{0.2mm}
\begin{picture}(280,290)
\savebox{\artin}(20,20)[1]{\masterartin} 
\savebox{\artininv}(20,20)[1]{\masterartininv} 
\savebox{\strandr}(20,20)[1]{\masterstrandr} 
\savebox{\strandl}(20,20)[1]{\masterstrandl} 
\savebox{\horizontaldots}(20,20)[1]{\masterhorizontaldots}
\put(0,260){\usebox{\artin}}  
\put(20,260){\usebox{\strandr}}
\put(40,260){\usebox{\strandr}}
\put(60,260){\usebox{\strandr}}
\put(0,240){\usebox{\strandl}}
\put(20,240){\usebox{\artin}}
\put(40,240){\usebox{\strandr}}
\put(60,240){\usebox{\strandr}}
\put(0,220){\usebox{\strandl}}
\put(20,220){\usebox{\strandl}}
\put(40,220){\usebox{\artin}}
\put(60,220){\usebox{\strandr}}
\put(0,200){\usebox{\strandl}}
\put(20,200){\usebox{\artininv}}
\put(40,200){\usebox{\strandr}}
\put(60,200){\usebox{\strandr}}
\put(0,180){\usebox{\artininv}}
\put(20,180){\usebox{\strandr}}
\put(40,180){\usebox{\strandr}}
\put(60,180){\usebox{\strandr}}
\put(0,160){\usebox{\artin}}  
\put(20,160){\usebox{\strandr}}
\put(40,160){\usebox{\strandr}}
\put(60,160){\usebox{\strandr}}
\put(0,140){\usebox{\strandl}}
\put(20,140){\usebox{\artin}}
\put(40,140){\usebox{\strandr}}
\put(60,140){\usebox{\strandr}}
\put(0,120){\usebox{\artininv}}
\put(20,120){\usebox{\strandr}}
\put(40,120){\usebox{\strandr}}
\put(60,120){\usebox{\strandr}}
\put(0,100){\usebox{\artin}}
\put(20,100){\usebox{\strandr}}
\put(40,100){\usebox{\strandr}}
\put(60,100){\usebox{\strandr}}
\put(80,100){\usebox{\horizontaldots}}
\put(0,80){\usebox{\artin}}  
\put(20,80){\usebox{\strandr}}
\put(40,80){\usebox{\strandr}}
\put(60,80){\usebox{\strandr}}
\put(0,60){\usebox{\strandl}}
\put(20,60){\usebox{\artin}}
\put(40,60){\usebox{\strandr}}
\put(60,60){\usebox{\strandr}}
\put(0,40){\usebox{\strandl}}
\put(20,40){\usebox{\strandl}}
\put(40,40){\usebox{\artin}}
\put(60,40){\usebox{\strandr}}
\put(0,20){\usebox{\strandl}}
\put(20,20){\usebox{\artininv}}
\put(40,20){\usebox{\strandr}}
\put(60,20){\usebox{\strandr}}
\put(0,0){\usebox{\artininv}}
\put(20,0){\usebox{\strandr}}
\put(40,0){\usebox{\strandr}}
\put(60,0){\usebox{\strandr}}
\put(130,100){$\equiv$}
\put(200,260){\usebox{\strandl}}
\put(200,260){\usebox{\strandr}}
\put(220,260){\usebox{\strandr}}
\put(240,260){\usebox{\strandr}}
\put(260,260){\usebox{\strandr}}
\put(200,240){\usebox{\artin}}
\put(220,240){\usebox{\strandr}}
\put(240,240){\usebox{\strandr}}
\put(260,240){\usebox{\strandr}}
\put(200,220){\usebox{\strandl}}
\put(220,220){\usebox{\artin}}
\put(240,220){\usebox{\strandr}}
\put(260,220){\usebox{\strandr}}
\put(200,200){\usebox{\artininv}}
\put(220,200){\usebox{\strandr}}
\put(240,200){\usebox{\strandr}}
\put(260,200){\usebox{\strandr}}
\put(200,180){\usebox{\artin}}
\put(220,180){\usebox{\strandr}}
\put(240,180){\usebox{\strandr}}
\put(260,180){\usebox{\strandr}}
\put(200,160){\usebox{\artin}}
\put(220,160){\usebox{\strandr}}
\put(240,160){\usebox{\strandr}}
\put(260,160){\usebox{\strandr}}
\put(200,140){\usebox{\strandl}}
\put(220,140){\usebox{\artin}}
\put(240,140){\usebox{\strandr}}
\put(260,140){\usebox{\strandr}}
\put(200,120){\usebox{\strandl}}
\put(220,120){\usebox{\strandl}}
\put(240,120){\usebox{\artin}}
\put(260,120){\usebox{\strandr}}
\put(200,100){\usebox{\strandl}}
\put(220,100){\usebox{\artininv}}
\put(240,100){\usebox{\strandr}}
\put(260,100){\usebox{\strandr}}
\put(280,100){\usebox{\horizontaldots}}
\put(200,80){\usebox{\artininv}}
\put(220,80){\usebox{\strandr}}
\put(240,80){\usebox{\strandr}}
\put(260,80){\usebox{\strandr}}
\put(200,60){\usebox{\artin}}
\put(220,60){\usebox{\strandr}}
\put(240,60){\usebox{\strandr}}
\put(260,60){\usebox{\strandr}}
\put(200,40){\usebox{\strandl}}
\put(220,40){\usebox{\artin}}
\put(240,40){\usebox{\strandr}}
\put(260,40){\usebox{\strandr}}
\put(200,20){\usebox{\artininv}}
\put(220,20){\usebox{\strandr}}
\put(240,20){\usebox{\strandr}}
\put(260,20){\usebox{\strandr}}
\put(200,00){\usebox{\strandl}}
\put(200,00){\usebox{\strandr}}
\put(220,00){\usebox{\strandr}}
\put(240,00){\usebox{\strandr}}
\put(260,00){\usebox{\strandr}}
\end{picture}
\caption{Braid relation  $\gamma_3 \gamma_2 \gamma_1 \gamma_3 
= \gamma_2 \gamma_1 \gamma_3 \gamma_2$}
\label{figure:braidrel-ii}
\end{figure}
Conversely, suppose the group $G$ is generated by $\gamma_1, \gamma_2, 
\cdots \gamma_{n-1}$ subject to the relations \eqref{eq:EB}. We denote by 
$\gamma_0$ the identity of $G$ and show that, as an intermediate step, the 
`new' generators
\[
\tsigma_i := \gamma_{i-1}^{-1}\gamma_i^{} \gamma_{i-1}^{}, \qquad  0 < i <n,  
\]
satisfy the relations \eqref{eq:B1} and \eqref{eq:B2}. We begin with proving 
$\tsigma_i \tsigma_j = \tsigma_{j} \tsigma_{i}$ for $|i-j|>1$ and assume 
without loss of generality $i+1<j$. For this purpose, multiply \eqref{eq:EB} 
by $\gamma_{l-1}^{-1}$ from the left side and by $\gamma_{l}^{-1}$ from the 
right side to obtain
\[ 
(\gamma_{l-1}^{-1} \gamma_{l} \gamma_{l-1}) 
(\gamma_{l-2} \gamma_{l-3} \cdots \gamma_{k+1}\gamma_{k}) 
= 
(\gamma_{l-2} \gamma_{l-3}\cdots \gamma_{k+1}\gamma_{k}) 
(\gamma_{l}\gamma_{l-1} \gamma_{l}^{-1}).
\]    
The relations \eqref{eq:EB} include the braid relations $\gamma_{l} 
\gamma_{l-1} \gamma_l = \gamma_{l-1} \gamma_{l} \gamma_{l-1}$, which entails
\begin{align}\label{eq:srbraid-iii} 
\tsigma_l = \gamma_{l-1}^{-1} \gamma_{l}^{} \gamma_{l-1}^{} 
=  \gamma_{l}^{} \gamma_{l-1}^{} \gamma_{l}^{-1}. 
\end{align}
This establishes the following commutation relations:  
\begin{eqnarray}\label{eq:gamma2sigma-cr}
\tsigma_l (\gamma_{l-2} \gamma_{l-3} \cdots \gamma_{k+1}\gamma_{k}) 
= (\gamma_{l-2} \gamma_{l-3} \cdots \gamma_{k+1}\gamma_{k}) \tsigma_l
\qquad \text{for $0 < k+1 < l$}.
\end{eqnarray}
It follows that $\tsigma_l$ commutes with all $\gamma_i$ with
$i \le l-2$ and hence also with all $\tsigma_i$ with $i \le l-2$. 
This establishes the relation \eqref{eq:B2} for the $\tsigma_i$'s. 

We are left to prove $\tsigma_l \tsigma_{k} \tsigma_l = \tsigma_k \tsigma_l \tsigma_k$ for $|k-l|=1$ and assume
without loss of generality that $l=k+1$. We have already shown that $\tsigma_{k+1}$ and $\gamma_{k-1}$ commute.
Thus, also using \eqref{eq:srbraid-iii},
\begin{eqnarray*}
\tsigma_k \tsigma_{k+1}\tsigma_k 
&=& (\gamma_{k-1}^{-1} \gamma_{k}^{}\gamma_{k-1}^{})\tsigma_{k+1}^{} (\gamma_{k-1}^{-1} \gamma_{k}^{}\gamma_{k-1}^{})\\
&=& (\gamma_{k-1}^{-1} \gamma_{k}^{}) \tsigma_{k+1}^{} (\gamma_{k}^{}\gamma_{k-1}^{})
= (\gamma_{k-1}^{-1} \gamma_{k}^{})  \gamma_{k}^{-1} \gamma_{k+1}^{}\gamma_{k}^{}  (\gamma_{k}^{}\gamma_{k-1}^{})\\
&=& \gamma_{k-1}^{-1} \gamma_{k+1}^{}\gamma_{k}^{2} \gamma_{k-1}^{}.
\end{eqnarray*}  
On the other hand,
\begin{eqnarray*}
\tsigma_{k+1}^{} \tsigma_{k}^{}\tsigma_{k+1}^{}
&=& \tsigma_{k+1}^{} (\gamma_{k-1}^{-1} \gamma_k \gamma_{k-1}^{})\tsigma_{k+1}^{}\\
&=& \gamma_{k-1}^{-1} \tsigma_{k+1}^{}  \gamma_k \tsigma_{k+1}^{} \gamma_{k-1}^{}\\
&=& \gamma_{k-1}^{-1} (\gamma_{k}^{-1}\gamma_{k+1}^{} \gamma_{k}^{})   
\gamma_k (\gamma_{k}^{-1}\gamma_{k+1}^{} \gamma_{k}^{})  \gamma_{k-1}^{}\\
&=& \gamma_{k-1}^{-1} \gamma_{k}^{-1}(\gamma_{k+1}^{} \gamma_{k}^{} 
\gamma_{k+1}^{}) \gamma_{k}^{}  \gamma_{k-1}^{}\\
&=& \gamma_{k-1}^{-1} \gamma_{k}^{-1}(\gamma_{k}^{} \gamma_{k+1}^{} 
\gamma_{k}^{}) \gamma_{k}^{}  \gamma_{k-1}^{}\\
&=& \gamma_{k-1}^{-1} \gamma_{k+1}^{}\gamma_{k}^{2} \gamma_{k-1}^{}.
\end{eqnarray*}  
Altogether we have shown that the relations \eqref{eq:EB} for the $\gamma_i$'s imply
the relations \eqref{eq:B1} and  \eqref{eq:B2} for the $\tsigma_j$'s. 
Finally, we note that the commutation relations \eqref{eq:gamma2sigma-cr} imply
\[
\tsigma_l =  (\gamma_{1}^{-1} \gamma_{2}^{-1} \cdots \gamma_{l-3}^{-1}\gamma_{l-2}^{-1})  \tsigma_l^{} (\gamma_{l-2}^{} \gamma_{l-3}^{} \cdots \gamma_{2}^{}\gamma_{1}^{}).
\] 
This establishes that the $\sigma_k$'s defined by \eqref{eq:sigma2gamma} 
satisfy the braid relations $\eqref{eq:B1}$ and $\eqref{eq:B2}$.
\end{proof}
\begin{Remark}\normalfont
The relations \eqref{eq:EB} are `extended' versions of the Artin braid relation 
\[
\gamma_{l} \gamma_{l-1} \gamma_{l} =  \gamma_{l-1} \gamma_{l} \gamma_{l-1}.
\] 
On the other hand, the generators $\gamma_l$ are a subset of 
$\set{\sigma_{s,t}}{1 \le s < t \le n}$, the generators of the Birman-Ko-Lee 
presentation (cf.\ \cite{BKL98a}). Following the notation of \cite{BiBr05a}, 
these generators are of the form  
\[
\sigma_{s,t} = (\sigma_{t-1}\sigma_{t-2}\cdots \sigma_{s+1})\sigma_s (\sigma_{s+1}^{-1}\cdots \sigma_{t-2}^{-1}\sigma_{t-1}^{-1})
\]
and $\gamma_l = \sigma_{1,l+1}$ is easily verified in the geometric picture. 
\end{Remark}
\begin{Remark}\normalfont
We are indebted to Patrick Dehornoy \cite{Deho08a} for the observation that 
the `square root of free generators presentation' is closely related to the 
Sergiescu presentations associated to planar graphs \cite{Se93a}. In fact, 
the relations \eqref{eq:EB} for $\gamma_1,\ldots,\gamma_{n-1}$ are also 
satisfied for the $n-1$ generators associated to a star-shaped graph with 
$n-1$ edges. But then we conclude from Theorem \ref{thm:sqrt-rep} together 
with the Hopfian property of $\Bset_n$ (see \cite{MaKaSo76a}, Section 3.7 
and 6.5) that $\gamma_1,\ldots,\gamma_{n-1}$ satisfy all relations associated 
to a star-shaped graph with $n-1$ edges. This shows that there is much more 
symmetry involved than originally stated, for example we have 
\[ 
\gamma_j \gamma_k \gamma_j = \gamma_k \gamma_j \gamma_k \quad\quad\quad 1 \le j,k \le n-1 
\] 
\[ \gamma_l \gamma_k \gamma_j \gamma_l = \gamma_k \gamma_j \gamma_l \gamma_k = \gamma_j \gamma_l \gamma_k \gamma_j \quad\quad 1 \le j < k < l \le n-1 
\] 
etc. It is a nice exercise to give direct algebraic proofs of such 
relations from \eqref{eq:EB}.

We would also like to thank Joan Birman \cite{Bir08a} who pointed out to us 
the more recent work of Han and Ko on positive braid group presentations from 
linearly spanned graphs which contains another alternative (minimal) collection 
of relations which is equivalent to \eqref{eq:EB} (see Lemma 3.3 and its proof 
in \cite{HaKo02a} for further details).
\end{Remark}
\begin{Remark}\normalfont
The group $\Bset_n$ is also presented by the set of generators $\set{\tgamma_i}{1 \le i \le n-1}$
subject to the defining relations
\begin{align}\tag{$\widetilde{\mathrm{EB}}$} \label{eq:EBtilde} 
      \tgamma_{l} (\tgamma_{k} \tgamma_{k+1}\cdots \tgamma_{l-2}) \tgamma_{l-1} \tgamma_{l} 
=   \tgamma_{l-1} \tgamma_l(\tgamma_{k} \tgamma_{k+1}\cdots \tgamma_{l-2}) \tgamma_{l-1}     
\end{align} 
for $0 < k < l <n$.
This is immediate from Proposition \ref{thm:sqrt-rep} if one lets $\tgamma_i:= \inv(\gamma_i^{-1})$ or if one 
notes that the transformation $\tau \mapsto \inv(\tau^{-1})$ reverses the ordering of the letters of the word $\tau$
(given in terms of the $\sigma_i^{\pm 1}$'s).

The Garside fundamental braid $\Delta_n$ (see Remark \ref{rem:garside})
has a simple form in terms of the 
$\gamma_k$'s: 
\[
\Delta_n= (\gamma_{n-1} \cdots \gamma_1)(\gamma_{n-2}\cdots \gamma_1) \cdots (\gamma_2 \gamma_1) \gamma_1.  
\]
Also the fundamental braid $\delta_n$ takes a simple form, but now in the 
terms of the generators $\tgamma_i$:
\[
\delta_n= \tgamma_{n-1} \tgamma_{n-2}\cdots \tgamma_2 \tgamma_1.
\] 
A more detailed investigation of these 
presentations with relations \eqref{eq:EB} or \eqref{eq:EBtilde} may be 
of interest on its own for the word and conjugacy problem in braid groups.
\end{Remark}
We continue with the introduction of several closely related shifts on 
$\Bset_\infty$ and relate them to the geometric operation of inserting 
`handles' into braid diagrams. This geometric approach turns out to be 
fruitful for results on orbits generated by the action of a shift on braids, 
as they are needed in Section \ref{section:left-regular-rep}. For further 
material on braid handles, the reader is referred to \cite{Deho97b} and 
\cite[Chapter III]{Deho00a}. 
\begin{Definition}\normalfont \label{def:braidshifts} 
The \emph{shift} $\sh$ is given by the endomorphism on $\Bset_\infty$ 
defined by
\[
\sh(\sigma_n) = \sigma_{n+1}
\] 
for all $n \in \Nset$. The \emph{$m$-shift} $\sh_m$ on $\Bset_\infty$, with 
fixed $m\in \Nset$, is given by the endomorphism 
\[
\sh_m(\tau) := \sigma_m \sigma_{m-1} \cdots 
\sigma_1 \sh(\tau) \sigma_1^{-1} \cdots \sigma_{m-1}^{-1}\sigma_m^{-1}.
\]
\end{Definition}
\begin{Lemma}\label{lem:shift_1}
The endomorphisms $\sh$ and $\sh_m$ on $\Bset_\infty$ are injective for all 
$m \in \Nset$. Moreover there exists, for every $\tau \in \Bset_\infty$, 
some $n \in \Nset$ such that 
\begin{eqnarray}\label{eq:shift-1}
\sh_m(\tau) = (\sigma_{m+1}^{-1}\sigma_{m+2}^{-1} \cdots 
\sigma_{n-1}^{-1}\sigma_n^{-1})\, \tau \,(\sigma_n^{} \sigma_{n-1}^{} 
\cdots \sigma_{m+2}^{} \sigma_{m+1}^{}).
\end{eqnarray}
\end{Lemma}
\begin{proof}
A geometric proof for the injectivity of $\sh$ is given in 
\cite[Chapter I, Lemma 3.3]{Deho00a}. The endomorphism $\sh_m$ is the 
composition of an injective endomorphism and automorphisms, and thus 
injective. We observe that, using the braid relations \eqref{eq:B1} and  
\eqref{eq:B2},
\begin{eqnarray*}
\sigma_{i+1} = \big(\sigma_1^{-1} \sigma_2^{-1} 
\cdots \sigma_{i-1}^{-1}\sigma_i^{-1} \sigma_{i+1}^{-1}\big) 
\sigma_i \big(\sigma_{i+1} \sigma_i \sigma_{i-1} \cdots \sigma_{2} 
\sigma_1\big) 
= \sh(\sigma_{i}).   
\end{eqnarray*}
Now let $\tau \in \Bset_\infty$ be given. Then there exists some 
$n \in \Nset$ such that $\tau \in \Bset_n$ and 
\[
\sh(\tau) = (\sigma_1^{-1} \sigma_2^{-1} 
 \cdots \sigma_{n-1}^{-1}\sigma_n^{-1}) \tau (\sigma_n^{} \sigma_{n-1}^{} 
 \cdots \sigma_2^{} \sigma_1^{}).
\] 
This proves \eqref{eq:shift-1}, since $\sh_m(\tau) = 
\sigma_m \cdots \sigma_1 \sh(\tau) \sigma_1^{-1} \cdots \sigma_m^{-1}$. 
\end{proof}
It is advantageous at this point to use the picture of geometric braids 
(cf. \cite[Section I.1]{Deho00a}) and to introduce the notion of braid 
handles to visualize the action of the shifts $\sh$ and $\sh_m$. Here the 
action of $\sh$ on some initial braid $\tau \in \Bset_\infty$ corresponds 
to inserting a new strand on the left of the initial braid diagram. Since 
$\sh_1(\tau) = \sigma_1 \sh(\tau) \sigma_1^{-1}$, we see that the action 
of $\sh_1$ corresponds to inserting a new strand between the first and 
second strand of the braid diagram for $\tau$ and above the other strands. 
This looks like an upper handle, as illustrated in Figure 
\ref{figure:one-shift}. 

\begin{figure}[h]
\setlength{\unitlength}{0.2mm}
\begin{picture}(400,100)
\savebox{\artin}(20,20)[1]{\masterartin} 
\savebox{\artininv}(20,20)[1]{\masterartininv} 
\savebox{\strandr}(20,20)[1]{\masterstrandr} 
\savebox{\strandl}(20,20)[1]{\masterstrandl} 
\savebox{\horizontaldots}(20,20)[1]{\masterhorizontaldots}
\put(0,60){\usebox{\strandl}}
\put(20,60){\usebox{\strandl}}
\put(40,60){\usebox{\strandl}}
\put(60,60){\usebox{\strandl}}
\put(80,60){\usebox{\strandl}}
\put(100,60){\usebox{\strandl}}
\put(120,60){\usebox{\strandl}}
\put(100,40){\usebox{\strandl}}
\put(120,40){\usebox{\strandl}}
\linethickness{1pt}
\put(0,60){\line(1,0){85}}
\put(0,60){\line(-1,0){5}}
\put(0,20){\line(1,0){85}}
\put(0,20){\line(-1,0){5}}
\put(85,20){\line(0,1){40}}
\put(-5,20){\line(0,1){40}}
\put(40,35){$\tau$}
\put(100,20){\usebox{\strandl}}
\put(120,20){\usebox{\strandl}}
\put(120,30){\usebox{\horizontaldots}}
\put(0,00){\usebox{\strandl}}
\put(20,00){\usebox{\strandl}}
\put(40,00){\usebox{\strandl}}
\put(40,00){\usebox{\strandr}}
\put(60,00){\usebox{\strandr}}
\put(80,00){\usebox{\strandr}}
\put(100,00){\usebox{\strandr}}
\put(185,35){$\sh_1$}
\put(165,30){\vector(1,0){60}}
\put(233,68){$\sigma_1$}
\put(260,60){\usebox{\artin}}
\put(280,60){\usebox{\strandr}}
\put(300,60){\usebox{\strandr}}
\put(320,60){\usebox{\strandr}}
\put(340,60){\usebox{\strandr}}
\put(360,60){\usebox{\strandr}}
\linethickness{1pt}
\put(280,60){\line(1,0){85}}
\put(280,60){\line(-1,0){5}}
\put(280,20){\line(1,0){85}}
\put(280,20){\line(-1,0){5}}
\put(275,20){\line(0,1){40}}
\put(365,20){\line(0,1){40}}
\put(310,35){$\sh(\tau)$}
\put(260,40){\usebox{\strandl}}
\put(360,40){\usebox{\strandr}}
\put(260,20){\usebox{\strandl}}
\put(360,20){\usebox{\strandr}}
\put(380,30){\usebox{\horizontaldots}}
\put(233,0){$\sigma_1^{-1}$}
\put(260,0){\usebox{\artininv}}
\put(280,00){\usebox{\strandr}}
\put(300,00){\usebox{\strandr}}
\put(320,00){\usebox{\strandr}}
\put(340,00){\usebox{\strandr}}
\put(360,00){\usebox{\strandr}}
\end{picture}
\caption{Braid diagram of the action $\tau \mapsto \sh_1(\tau)
= \sigma_1 \sh(\tau)\sigma_1^{-1}$}
\label{figure:one-shift}
\end{figure}
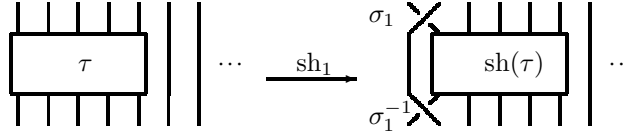

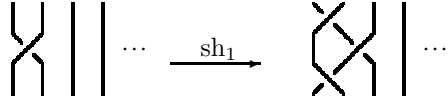
\begin{figure}[h]
\setlength{\unitlength}{0.2mm}
\begin{picture}(260,80)
\savebox{\artin}(20,20)[1]{\masterartin} 
\savebox{\artininv}(20,20)[1]{\masterartininv} 
\savebox{\strandr}(20,20)[1]{\masterstrandr} 
\savebox{\strandl}(20,20)[1]{\masterstrandl} 
\savebox{\horizontaldots}(20,20)[1]{\masterhorizontaldots}
\put(0,40){\usebox{\strandl}}
\put(20,40){\usebox{\strandl}}
\put(40,40){\usebox{\strandl}}
\put(60,40){\usebox{\strandl}}
\put(0,20){\usebox{\artin}}
\put(20,20){\usebox{\strandr}}
\put(40,20){\usebox{\strandr}}
\put(60,20){\usebox{\horizontaldots}}
\put(0,00){\usebox{\strandl}}
\put(20,00){\usebox{\strandl}}
\put(40,00){\usebox{\strandl}}
\put(40,00){\usebox{\strandr}}
\put(125,25){$\sh_1$}
\put(105,20){\vector(1,0){60}}
\put(200,40){\usebox{\artin}}
\put(220,40){\usebox{\strandr}}
\put(240,40){\usebox{\strandr}}
\put(200,20){\usebox{\strandl}}
\put(220,20){\usebox{\artin}}
\put(240,20){\usebox{\strandr}}
\put(260,20){\usebox{\horizontaldots}}
\put(200,00){\usebox{\artininv}}
\put(220,00){\usebox{\strandr}}
\put(240,00){\usebox{\strandr}}
\end{picture}
\caption{Braid diagram of the action 
$\sigma_1 \mapsto \sh_1(\sigma_1)= \sigma_1 \sigma_2 \sigma_1^{-1}$}
\label{figure:handle-i}
\end{figure}

The new strand is not entangled with the strands of the initial braid $\tau$. 
In particular, $\sigma_1 \mapsto \sh_1(\sigma_1) = \sigma_1^{} \sigma_2 
\sigma_1^{-1}$ (see Figure \ref{figure:handle-i}) and $\sh_1(\sigma_i) = 
\sigma_{i+1}$ for $i\ge 2$. Thus the action  of $\sh_1$ coincides with the 
action of $\sh$ on $\sigma_i$ for $i \ge 2$.     

Generalizing the above discussion and motivated by 
\cite[Definition 3.3]{Deho00a}, we introduce:
\begin{Definition}\normalfont \label{def:handle}
An \emph{upper $m$-handle} is a braid of the form $\sh_m(\tau)$ with 
$m\ge 1$ and $\tau \in \Bset_\infty$.
\end{Definition}
\begin{figure}[h]
\setlength{\unitlength}{0.2mm}
\begin{picture}(400,170)
\savebox{\artin}(20,20)[1]{\masterartin} 
\savebox{\artininv}(20,20)[1]{\masterartininv} 
\savebox{\strandr}(20,20)[1]{\masterstrandr} 
\savebox{\strandl}(20,20)[1]{\masterstrandl} 
\savebox{\horizontaldots}(20,20)[1]{\masterhorizontaldots}
\put(0,140){\usebox{\strandl}}
\put(20,140){\usebox{\strandl}}
\put(40,140){\usebox{\strandl}}
\put(60,140){\usebox{\strandl}}
\put(80,140){\usebox{\strandl}}
\put(100,140){\usebox{\strandl}}
\put(120,140){\usebox{\strandl}}
\put(0,120){\usebox{\strandl}}
\put(20,120){\usebox{\strandl}}
\put(40,120){\usebox{\strandl}}
\put(60,120){\usebox{\strandl}}
\put(80,120){\usebox{\strandl}}
\put(100,120){\usebox{\strandl}}
\put(120,120){\usebox{\strandl}}
\put(0,100){\usebox{\strandl}}
\put(20,100){\usebox{\strandl}}
\put(40,100){\usebox{\strandl}}
\put(60,100){\usebox{\strandl}}
\put(80,100){\usebox{\strandl}}
\put(100,100){\usebox{\strandl}}
\put(120,100){\usebox{\strandl}}
\put(100,80){\usebox{\strandl}}
\put(120,80){\usebox{\strandl}}
\linethickness{1pt}
\put(0,100){\line(1,0){85}}
\put(0,100){\line(-1,0){5}}
\put(0,60){\line(1,0){85}}
\put(0,60){\line(-1,0){5}}
\put(85,60){\line(0,1){40}}
\put(-5,60){\line(0,1){40}}
\put(40,75){$\tau$}
\put(100,60){\usebox{\strandl}}
\put(120,60){\usebox{\strandl}}
\put(120,70){\usebox{\horizontaldots}}
\put(0,40){\usebox{\strandl}}
\put(20,40){\usebox{\strandl}}
\put(40,40){\usebox{\strandl}}
\put(40,40){\usebox{\strandr}}
\put(60,40){\usebox{\strandr}}
\put(80,40){\usebox{\strandr}}
\put(100,40){\usebox{\strandr}}
\put(0,20){\usebox{\strandl}}
\put(20,20){\usebox{\strandl}}
\put(40,20){\usebox{\strandl}}
\put(40,20){\usebox{\strandr}}
\put(60,20){\usebox{\strandr}}
\put(80,20){\usebox{\strandr}}
\put(100,20){\usebox{\strandr}}
\put(0,0){\usebox{\strandl}}
\put(20,0){\usebox{\strandl}}
\put(40,0){\usebox{\strandl}}
\put(40,0){\usebox{\strandr}}
\put(60,0){\usebox{\strandr}}
\put(80,0){\usebox{\strandr}}
\put(100,0){\usebox{\strandr}}
\put(185,75){$\sh_3$}
\put(165,70){\vector(1,0){60}}
\put(233,148){$\sigma_3$}
\put(260,140){\usebox{\strandl}}
\put(280,140){\usebox{\strandl}}
\put(300,140){\usebox{\artin}}
\put(320,140){\usebox{\strandr}}
\put(340,140){\usebox{\strandr}}
\put(360,140){\usebox{\strandr}}
\put(233,128){$\sigma_2$}
\put(260,120){\usebox{\strandl}}
\put(280,120){\usebox{\artin}}
\put(300,120){\usebox{\strandr}}
\put(320,120){\usebox{\strandr}}
\put(340,120){\usebox{\strandr}}
\put(360,120){\usebox{\strandr}}
\put(233,108){$\sigma_1$}
\put(260,100){\usebox{\artin}}
\put(280,100){\usebox{\strandr}}
\put(300,100){\usebox{\strandr}}
\put(320,100){\usebox{\strandr}}
\put(340,100){\usebox{\strandr}}
\put(360,100){\usebox{\strandr}}
\linethickness{1pt}
\put(280,100){\line(1,0){85}}
\put(280,100){\line(-1,0){5}}
\put(280,60){\line(1,0){85}}
\put(280,60){\line(-1,0){5}}
\put(275,60){\line(0,1){40}}
\put(365,60){\line(0,1){40}}
\put(310,75){$\sh(\tau)$}
\put(260,80){\usebox{\strandl}}
\put(360,80){\usebox{\strandr}}
\put(260,60){\usebox{\strandl}}
\put(360,60){\usebox{\strandr}}
\put(380,70){\usebox{\horizontaldots}}
\put(230,40){$\sigma_1^{-1}$}
\put(260,40){\usebox{\artininv}}
\put(280,40){\usebox{\strandr}}
\put(300,40){\usebox{\strandr}}
\put(320,40){\usebox{\strandr}}
\put(340,40){\usebox{\strandr}}
\put(360,40){\usebox{\strandr}}
\put(230,20){$\sigma_2^{-1}$}
\put(260,20){\usebox{\strandl}}
\put(280,20){\usebox{\artininv}}
\put(300,20){\usebox{\strandr}}
\put(320,20){\usebox{\strandr}}
\put(340,20){\usebox{\strandr}}
\put(360,20){\usebox{\strandr}}
\put(230,0){$\sigma_3^{-1}$}
\put(260,0){\usebox{\strandl}}
\put(280,0){\usebox{\strandl}}
\put(300,0){\usebox{\artininv}}
\put(320,0){\usebox{\strandr}}
\put(340,0){\usebox{\strandr}}
\put(360,0){\usebox{\strandr}}
\end{picture}
\caption{Braid diagram of the action $\tau \mapsto \sh_3(\tau)
= \sigma_3 \sigma_2 \sigma_1 \sh(\tau)\sigma_1^{-1} \sigma_2^{-1}\sigma_3^{-1}$}
\label{figure:m-shift}
\end{figure}
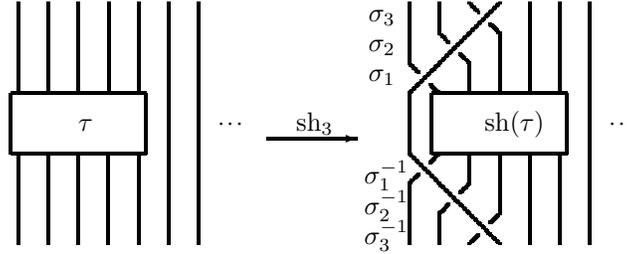
We discuss the geometric interpretation of the shift $\sh_m$ for $m \ge 1$. 
The action of $\sh_m$ on the $\sigma_i$'s can easily be identified to be
\begin{eqnarray*}
\sh_m(\sigma_i) =
\begin{cases}
\sigma_i & \text{ if $i < m$;}\\
\sigma_i \sigma_{i+1} \sigma_i^{-1}  & \text{ if $i = m$;}\\
\sigma_{i+1} &  \text{ if $i > m$.}
\end{cases}
\end{eqnarray*} 
The action of $\sh_m$ turns the initial braid $\tau$ into an upper 
$m$-handle and geometrically corresponds to inserting an upper strand 
between the $m$-th and $(m+1)$-th strand of the initial braid diagram
(see Figure \ref{figure:m-shift} for $m=3$). Further simplification occurs 
if we replace the generators $\sigma_1, \sigma_2, \ldots $ by the square 
roots of free generators $\gamma_1, \gamma_2, \ldots $ and express the 
initial braid $\tau$ in terms of the $\gamma_i^{\pm 1}$. 
\begin{Lemma}\label{lem:shift_2} 
\begin{eqnarray*}
\sh_m(\gamma_i) =
\begin{cases}
\gamma_i & \text{ if $i < m$;}\\
\gamma_{i+1}  & \text{ if $i \ge m$.}
\end{cases}
\end{eqnarray*}
In particular, for all $n \in \Nset$, 
\[
\sh_1(\gamma_n)= \gamma_{n+1}.
\]
\end{Lemma}
We see that as soon as we work in the square root of free generators 
presentation the $1$-shift is nothing but the shift in the generators while 
the $m$-shifts for $m>1$ are partial shifts. The upper $m$-handle $\sh_m(\tau)$ 
can be expressed in terms of the generators $\gamma_i^{\pm 1}$ without
the appearance of $\gamma_m^{\pm 1}$. 

Now we provide results on relative conjugacy classes which will be needed in 
Section \ref{section:left-regular-rep}.
\begin{Definition} \normalfont \label{def:totalwidth}
The \emph{total width} $\tw(\tau)$ of a braid $\tau \in \Bset_\infty$ is 
defined as the function $\tw\colon \Bset_\infty \to \Nset_0$, where 
$\tw(\tau)$ is the minimal number $n \in \Nset_0$ such that $\tau \in \Bset_{n+1}$.
\end{Definition}
In other words: if $\tw(\tau) = n$, then we can express the braid $\tau$ as a
word in the symbols $\sigma_0, \sigma_1^{\pm 1}, \ldots, \sigma_{n-1}^{\pm 1}, 
\sigma_{n}^{\pm 1}$, but not as a word in $\sigma_0,  \sigma_1^{\pm 1}, \ldots, 
\sigma_{n-1}^{\pm 1}$. Note that $\tw(\tau)=0$ if and only if $\tau = \sigma_0$.
\begin{Remark}\normalfont
The \emph{total width} $\tw(\tau)$ should not be confused with the \emph{width} 
of a non-trivial braid $\tau$ introduced in \cite[Definition 3.16]{Deho00a}. 
The latter is the number of strands `really involved in $\tau$' which is less 
or equal $\tw(\tau)$. For example, $\sigma_6\sigma_7^{-1} \sigma_9^2$ has the 
width $(9-6+2)=5$, but the total width 9.
\end{Remark}
\begin{Proposition}\label{prop:infinite_cc}
Suppose $\sigma_0  \not= \tau \in \Bset_\infty$ and let $k,l \in \Nset_0$. 
Then we have:
\[
\tw(\sh_1(\tau)) = \tw(\tau) + 1.
\]
It particular, the set $\set{(\sh_1)^k(\tau)}{k \in \Nset_0}$ is infinite. 
If $\tw(\tau) \ge m$ then these assertions remain valid for
$\sh_m$ instead of $\sh_1$. 
\end{Proposition}
\begin{proof}
Suppose $\sigma_0 \not=\tau \in \Bset_\infty$ with total width $\tw(\tau) = n$. 
We have $\tw(\sh_1(\tau)) \le n+1$. In the following we will use geometric 
arguments to prove that $\tw(\sh_1(\tau)) = n+1$. We have already seen that 
the action of $\sh_1$ turns the initial braid $\tau$ into an upper 1-handle 
(see Figure \ref{figure:one-shift}). Suppose now that, using the braid 
relations \eqref{eq:B1} and \eqref{eq:B2}, it is possible to write the upper 
1-handle $\sh_1(\tau)$ as a word in $\sigma_1^{\pm 1}, \ldots, \sigma_{n}^{\pm 1}$ 
only. Then it is geometrically clear, because the new upper strand 
is not entangled, that we can use the same operations to write the initial braid 
$\tau$ as a word in $\sigma_1^{\pm 1}, \ldots, \sigma_{n-1}^{\pm 1}$, 
contrary to our assumption about $\tw(\tau) $. Hence $\tw(\sh_1(\tau))= n+1$. This 
proves the proposition for $\sh_1$. The general case can be done similarly. 
The condition $n = \tw(\tau) \ge m$ is needed to ensure that the new upper 
strand is inserted between the $n+1$ strands used to model $\Bset_{n+1}$, 
so that the geometric argument still works. 
\end{proof}
\begin{Proposition}\label{prop:infinite-ccc}
Let $C_m(\tau) := \set{w \tau w^{-1}}{w \in \Bset_{m+1,\infty}}$ denote the 
relative conjugacy class of $\tau \in \Bset_\infty$ for some $m \ge 1$. Then 
we have:
\begin{eqnarray*}
\tau \in \Bset_{m}  
\quad \Longleftrightarrow \quad  C_m(\tau) = \{\tau\}
\quad \Longleftrightarrow \quad  C_m(\tau) \text{ is finite}.  
\end{eqnarray*}
\end{Proposition}
\begin{proof}
The implication `$\Longrightarrow$' is immediate from the braid relation \eqref{eq:B2}. 
Conversely, let $m \ge 1$ be fixed and suppose $\tau \in \Bset_\infty$ with 
total width $\tw(\tau)\ge m$. In other words, we assume that the braid $\tau$ 
is not contained in $\Bset_m$. Now consider the shift $\sh_m$ and note that 
by Lemma \ref{lem:shift_1} always
\[
\sh_m(\bullet) = \sigma_{m+1}^{-1}\sigma_{m+2}^{-1}\cdots \sigma_n^{-1}\sigma_{n+1}^{-1}\bullet \sigma_{n+1}^{} \sigma_n^{}\cdots 
\sigma_{m+2}^{} \sigma_{m+1}^{}
\]
with some $n$ and hence $\set{(\sh_m)^k(\tau)}{k \in \Nset_0} \subset C_m(\tau)$. 
We conclude with Proposition \ref{prop:infinite_cc} that $C_m(\tau)$ is infinite.
\end{proof} 
\begin{Lemma}\label{lem:ICC-index}
$\Bset_\infty$ is an ICC group. Moreover, the inclusion $\Bset_{2,\infty} 
\subset \Bset_\infty$ has infinite group index $[\Bset_\infty \colon \Bset_{2,\infty}]$. 
\end{Lemma}
\begin{proof}
The first assertion is immediate from Proposition \ref{prop:infinite-ccc} 
for $m=1$. For the second assertion note that the cosets 
$(\sigma_1)^n \, \Bset_{2,\infty}$ are all different from each other for all 
$n \in \Nset$. This follows geometrically or from Theorem \ref{thm:lrr}(vi) 
below which implies that they are pairwise orthogonal with respect to the 
trace of the group algebra. 
\end{proof}
\section{An application to the group von Neumann algebra $L(\Bset_\infty)$}
\label{section:left-regular-rep}
This section is devoted to the construction of noncommutative random 
sequences from the left regular representation of the braid group 
$\Bset_\infty$. This will bring us in contact with free probability which 
has been introduced by Voiculescu for the study of free group von Neumann 
algebras. 

The group von Neumann algebra $L(\Bset_\infty)$ is generated by the 
left-regular representation $\set{L_\sigma}{\sigma \in \Bset_\infty}$ 
of $\Bset_\infty$ on the Hilbert space $\ell^2(\Bset_\infty)$, where  
\[
L_\sigma f(\sigma^\prime) := f(\sigma^{-1}\sigma^\prime). 
\]
Let $\delta_\sigma \in \ell^2(\Bset_\infty)$ be the function
\[
\delta_\sigma(\sigma^\prime) = 
\begin{cases}
1 & \text{if $\sigma = \sigma^\prime$}\\
0 & \text{otherwise}.
\end{cases}
\] 
Then the complex linear extension of 
\[
\trace_\infty(L_\sigma) := \langle \delta_{\sigma_0}, 
               L_\sigma \delta_{\sigma_0}\rangle 
\]
defines the normal faithful tracial state $\trace_\infty$ on 
$L(\Bset_\infty)$. We denote its restriction to 
$L(\Bset_n)$ by $\trace_n$.  
\begin{Definition} \label{def:irreducible}\normalfont
The inclusion $\cN \subset \cM$ of two von Neumann algebras is said to be 
\emph{irreducible} if the relative commutant is trivial, i.e.
\[
\cN^\prime \cap \cM \simeq \Cset.
\]
\end{Definition}
\begin{Theorem}\label{thm:irreducibility}
The inclusion $L(\Bset_{2,\infty}) \subset L(\Bset_\infty)$ is irreducible.
\end{Theorem}
\begin{proof}
In analogy with the standard arguments for  ICC groups 
(cf.\ \cite[V.7]{Take03a}), to prove irreducibility of the inclusion 
$L(\Bset_{2,\infty}) \subset L(\Bset_\infty)$ it is sufficient to show 
that the relative conjugacy class 
$C_1(\tau):= \set{w\tau w^{-1}}{w \in \Bset_{2,\infty}}$
is infinite for every $\tau \in \Bset_\infty$ with $\tau \not= \sigma_0$. 
But this follows from Proposition \ref{prop:infinite-ccc}.
\end{proof}
We collect some auxiliary results on $L(\Bset_\infty)$ which are 
of interest on its own. Since we could not find 
proofs in the literature, we provide them here for the convenience 
of the reader. The second author is indebted to Benoit Collins and Thierry 
Giordano for stimulating discussions on Property ($\Gamma$) resp.~Property (T) 
for braid groups. 
\begin{Corollary}\label{cor:factoriality-i} 
\begin{enumerate}
\item
$L(\Bset_\infty)$ is a non-hyperfinite $II_1$-factor; 
\item
$L(\Bset_\infty)$ has Property {\normalfont{($\Gamma$)}}; 
\item
$L(\Bset_\infty)$ does not have Kazhdan Property {\normalfont{(T)}};
\item 
$L(\Bset_{2,\infty}) \subset L(\Bset_\infty)$ is a subfactor inclusion 
with infinite Jones index.
\end{enumerate}
In particular, $L(\Bset_\infty)$ and $L(\Fset_n)$ are non-isomorphic 
($2 \le n \le \infty$). 
\end{Corollary}
\begin{proof}
(i)
The factoriality of $L(\Bset_\infty)$ is immediate from Theorem \ref{thm:irreducibility} 
and the definition of irreducibility. We are left to prove the non-hyperfiniteness. 
$\Bset_3$ and hence $\Bset_\infty$ contains $\Fset_2$ as a subgroup 
(cf. \cite[Lecture 5]{Jone91a}). Since any subgroup of an amenable discrete 
group is amenable (\cite[XIII, Example 4.4(iii)]{Take03c}) and $\Fset_2$ is 
non-amenable (\cite[XIII,Example 4.4(v)]{Take03c}), we conclude that the group 
$\Bset_\infty$ is non-amenable. But this implies the non-hyperfiniteness of 
$L(\Bset_\infty)$ (see \cite[XIII, Theorem 4.10]{Take03c}).  

(ii) Let $[x,y]:= xy -yx$. We conclude from the braid relation \eqref{eq:B2} that
\[
\lim_{k \to \infty} \tau\big([L_{\sigma_k}, x)]^* [L_{\sigma_k}, x)]\big) = 0
\] 
for any $x \in \alg\set{L_\sigma}{\sigma \in \Bset_{\infty}}$. This extends to 
$x \in L(\Bset_\infty)$ by a density argument. Thus $(L_{\sigma_k})_k$ is a 
non-trivial central sequence. 

(iii) $\Bset_\infty$ does not have Property (T). This can be deduced along the 
hints of \cite[Exercise 1.8.14]{BHK07aTA}, compare also \cite{GiHa91a}. But an 
ICC group has Property (T) if and only if its group von Neumann algebra does so 
\cite{Conn82a,CoJo85a}. Thus $L(\Bset_\infty)$ does not have Property (T).

(iv) Since $\sh(\Bset_\infty) = \Bset_{2,\infty}$, we conclude the factoriality 
of $L(\Bset_{2,\infty})$ from the injectivity of $\sh$ (see Lemma \ref{lem:shift_1}). 
The infinity of the Jones index is clear from Lemma \ref{lem:ICC-index} and 
$[\Bset_\infty \colon \Bset_{2,\infty}] = [L(\Bset_\infty) \colon L(\Bset_{2,\infty})]$ 
(see also \cite[Example 2.4]{Take03c}).

Finally, $L(\Bset_\infty)$ and $L(\Fset_n)$ are not isomorphic since $L(\Bset_\infty)$ 
has a non-trivial central sequence and thus cannot be full (see \cite[Theorem 3.8]{Take03c}). 
But the free group factors $L(\Fset_n)$ are full for $2 \le n \le \infty$ 
(see \cite[Theorem 3.9]{Take03c}).
\end{proof}
\begin{Remark}\normalfont
$L(\Bset_n)$ is \emph{not} a factor for $2 \le n <\infty$, since $Z(\Bset_n)$, 
the center of the group $\Bset_n$, is non-trivial.    
\end{Remark}
We are now going to identify the relative commutants 
$\big(L(\Bset_{n,\infty})\big)^\prime \cap L(\Bset_\infty)$ for $n \ge 2$. 
Note that the case $n=2$ is already covered by Proposition \ref{prop:infinite_cc} 
and Theorem \ref{thm:irreducibility}. 
\begin{Theorem} \label{thm:relative-commutants} 
We have $\big(L(\Bset_{n+1,\infty})\big)^\prime \cap L(\Bset_\infty) = 
L (\Bset_{n})$ for all $n \ge 1$. 
\end{Theorem} 
\begin{proof}
The inclusion $L (\Bset_{n}) \subset \big(L(\Bset_{n+1,\infty})\big)^\prime$ 
is clear from the braid relations \eqref{eq:B2}. For the converse inclusion, 
we conclude in analogy to arguments for ICC groups (cf.\ \cite[V.7]{Take03a}). Let 
$x \in L(\Bset_\infty)$ with $x= \sum_{\tau \in \Bset_\infty} x(\tau) L_{\tau}$,
where $x(\tau)$ are scalars such that $\sum_{\tau \in \Bset_\infty} |x(\tau)|^2 < \infty$. 
We have $x \in \big(L(\Bset_{n+1,\infty})\big)^\prime$ if and only if the coefficients 
$x(\tau)$ are constant on the relative conjugacy class $C_{n}(\tau)$. Thus the 
square summability of the non-zero coefficients $x(\sigma)$ implies that 
$C_{n}(\tau)$ is finite for every non-zero coefficient $x(\tau)$. We conclude 
from this that $\tau \in \Bset_n$ by Proposition \ref{prop:infinite-ccc}. This 
shows that $x \in \big(L(\Bset_{n+1,\infty})\big)^\prime$ implies $x\in L(\Bset_n)$.
\end{proof}
We turn our attention to braid group representations on the probability
space $(L(\Bset_\infty), \trace_\infty)$. Here we are interested in considering 
the representation
\[
\rho \colon \Bset_\infty \to \Aut{L(\Bset_\infty)},
\]
defined by $\rho(\sigma):= \Ad(L_\sigma) := L_\sigma \bullet L^*_\sigma$ 
with $\sigma \in \Bset_\infty$. We note that $\trace_\infty$ is automatically
$\rho(\sigma)$-invariant and thus $\rho(\sigma) \in \Aut{L(\Bset_\infty), \trace_\infty}$. 
Furthermore, the representation $\rho$ has the generating property 
(see Definition \ref{def:generating}):
\begin{eqnarray*}
L(\Bset_\infty) &\supset& \bigvee_{n \ge 0} \big(L(\Bset_\infty)\big)^{\rho(\Bset_{n+2,\infty})}                       
                         \supset \bigvee_{n \ge 0} L(\Bset_{n+1}) = L(\Bset_\infty).   
\end{eqnarray*}
At this point we have verified all assumptions of Theorem \ref{thm:endo-braid-i}. 
We next identify the fixed point algebras as they appear in Theorem \ref{thm:endo-braid-i}. 
We have from Theorem \ref{thm:irreducibility} that 
\[
\big(L(\Bset_\infty)\big)^{\rho(\Bset_{\infty})} = \cZ(L(\Bset_\infty)) \simeq \Cset,
\]
and from Theorem \ref{thm:relative-commutants} that, for all $n\ge 0$,
\[
\big(L(\Bset_\infty)\big)^{\rho(\Bset_{n+2,\infty})} 
= \big(L(\Bset_{n+2,\infty})\big)^\prime \cap L(\Bset_\infty)  = L(\Bset_{n+1}).  
\]
Furthermore, a straightforward computation shows that the action of the 
endomorphism $\alpha$ (from  Theorem \ref{thm:endo-braid-i}) 
comes from the shift $\sh$ on $\Bset_\infty$; more precisely:
\[
\alpha(L_\tau) = L_{\sh(\tau)}.
\]
Let us summarize the discussion above.
\begin{Theorem} \label{thm:lrr}
Consider the probability space $\big(L(\Bset_\infty), \trace_\infty\big)$, 
equipped with the representation  $\rho\colon \Bset_\infty \to 
\Aut{L(\Bset_\infty), \trace_\infty}$ given by $\tau \mapsto \rho(\tau):= \Ad L_{\tau}$. 
We arrive at the following conclusions:
\begin{enumerate}
\item \label{item:lrr-i}
$\rho$ has the generating property; 
\item \label{item:lrr-ii}
$\big(L(\Bset_\infty)\big)^{\rho(\Bset_\infty)} \simeq \Cset;$
\item \label{item:lrr-iii}
$\big(L(\Bset_\infty)\big)^{\rho(\Bset_{n+2,\infty})} = L(\Bset_{n+1})$ for all $n \ge 0$;
\item \label{item:lrr-iv}
the map 
$\alpha(x) := \sotlim_{n \to \infty} \rho(\sigma_1 \sigma_2 \cdots \sigma_n)(x)$ is an
      endomorphism for $(L(\Bset_\infty), \trace_\infty)$ such that $\alpha(L_\tau) = L_{\sh(\tau)}$, with 
      $\tau \in \Bset_\infty$;
\item \label{item:lrr-v}
Each cell of the triangular tower of inclusions is a commuting square: 
\begin{eqnarray*}
\setcounter{MaxMatrixCols}{20}
\begin{matrix}
\Cset &\subset&   \Cset  & \subset & L(\Bset_2) & \subset & L(\Bset_3) & \subset & L(\Bset_4) & \subset   \cdots \subset &  L(\Bset_\infty)\\
        &&          \cup  &         & \cup  &         & \cup &         & \cup  &        & \cup  \\
        &&   \Cset &\subset& \Cset&\subset&\alpha( L(\Bset_2))&\subset&\alpha( L(\Bset_3))&\subset \cdots  \subset & \alpha( L(\Bset_\infty))\\
         &&               &         & \cup  &         & \cup   &         & \cup  &        & \cup \\
              &&&&   \Cset& \subset &  \Cset  & \subset & \alpha^2( L(\Bset_2))
              & \subset   \cdots  \subset & \alpha^2( L(\Bset_\infty)) \\
       &&&& &&  \cup           &         & \cup  &         & \cup \\
        &&&&&&   \vdots           &         & \vdots  &        & \vdots 
\end{matrix}
\setcounter{MaxMatrixCols}{10}
\end{eqnarray*}
\item \label{item:lrr-vi}
The maps $\iota_n^{(\alpha)}:= \alpha^n|_{L(\Bset_2)}$ define a 
stationary and full $\Cset$-independent random sequence 
\[
\iota^{(\alpha)} \equiv (\iota_n^{(\alpha)})_{n \in \Nset_0} \colon 
\big(L(\Bset_2), \trace_2\big) \to \big(L(\Bset_\infty), \trace_\infty\big).
\] 
In other words, $\alpha$ is a full Bernoulli shift over $\Cset$ with
generator $L(\Bset_2)$ (see Definition \ref{def:bernoulli}). 
\item \label{item:lrr-vii}
The random sequence $\iota^{(\alpha)}$ is not spreadable.
\end{enumerate}
\end{Theorem}
\begin{Remark} \label{rem:lrr} \normalfont
For all $n \in \Nset_0$, 
\[
\iota_n(L_{\sigma_1}) = L_{\sh^n(\sigma_1)} =L_{\sigma_{n+1}}.
\]
Hence the left regular representation of the Artin generators gives us
a full $\Cset$-independent sequence which is not spreadable, and thus
also not braidable by Theorem \ref{thm:main-1}.
The random sequence $\iota^{(\alpha)}$ shows that the class of stationary and 
conditionally full independent random sequences is strictly larger than the 
class of spreadable random sequences in our setting of the noncommutative de 
Finetti theorem, Theorem \ref{thm:definetti}.  
\end{Remark}
\begin{proof} 
(\ref{item:lrr-i}) to (\ref{item:lrr-iv}) are shown above. (\ref{item:lrr-v}) is
immediate from (\ref{item:lrr-ii}), (\ref{item:lrr-iii})
and Theorem \ref{thm:endo-braid-i}. We are left to prove  (\ref{item:lrr-vi}) 
and (\ref{item:lrr-vii}).

The random sequence $\iota^{(\alpha)}$ is stationary, since it is induced by 
the endomorphism $\alpha$. We show first that  $\iota^{(\alpha)}$ is order 
$\Cset$-independent. For this purpose let $I = \{n_1, \ldots, m_1\}$ and 
$J =\{n_2, \ldots, m_2\}$ be `intervals' in $\Nset_0$ with $I < J$, or more 
explicitly: $n_1 \le  m_1 < n_2 \le m_2$. Indeed, $I <J$ implies the 
order $\Cset$-independence of $\bigvee_{k \in I} \alpha^k(L(\Bset_2))$ and 
$\bigvee_{l \in J}\alpha^l(L(\Bset_2))$ by the following arguments. We have 
$\alpha^k(L_{\sigma_1}) = L_{\sh^k(\sigma_1)} =L_{\sigma_{k+1}}$, from which 
we conclude
\[
\bigvee_{k \in I} \alpha^k(L(\Bset_2)) \subset L(\Bset_{m_1 +2})
\] 
and 
\[
\bigvee_{l \in J}\alpha^l(L(\Bset_2)) \subset \alpha^{m_1+1}\big(L(\Bset_{\infty})\big).
\] 
Looking at the triangular tower, this clearly implies the order $\Cset$-independence.
Moreover, this entails the order $\Cset$-independence of the random sequence $\iota^{(\alpha)}$.
  
We still need to show that order $\Cset$-independence upgrades to full $\Cset$-independence, 
valid whenever $I \cap J = \emptyset$ (see Definition \ref{def:order-independence}). 
We will prove this by induction.
Let 
\begin{eqnarray*}
\cA_I &:=& \bigvee_{k=1}^N \cA_{I_k}\qquad \text{with } \cA_{I_k} =  \bigvee_{i \in I_k} \alpha^i(L(\Bset_2))\\
\cA_J &:=& \bigvee_{l=1}^N \cA_{J_l}\qquad \text{with }  \cA_{J_l} =  \bigvee_{j \in J_l} \alpha^j(L(\Bset_2)).
\end{eqnarray*}
where the non-empty finite `intervals' $\{I_k\}_{k=1,...,N}$ and $\{J_l\}_{l=1,...,N}$ satisfy
\[
 I_1 < J_1 < I_2 < J_2 < \cdots < I_N < J_N.
\] 
Since $\alpha$ comes from the symbolic shift $\sh$ on the Artin generators $\sigma_i$ (see (\ref{item:lrr-iv})), we know from the assumptions on the `intervals' and the braid relation \eqref{eq:B2} that $\cA_{I_k}$ and $\cA_{I_{k^\prime}}$ commute for $k \neq k^\prime$. So do  $\cA_{J_l}$ and $\cA_{J_{l^\prime}}$ 
for $l \neq l^\prime$. We conclude from this by a simple induction on $k$ and $l$ that
\begin{eqnarray*}
\cA_{I^{(k)}} &=& \overline{\cA_{I_k} \cdot \cA_{I^{(k-1)}}}^{\,\sot}  
\qquad \text{for } I^{(k)}:=\bigcup_{k^\prime=1}^k I_{k^\prime},\\
\cA_{J^{(l)}} &=& \overline{\cA_{J^{(l-1)}} \cdot \cA_{J_{l}}}^{\,\sot}
 \qquad \text{for } J^{(l)}=\bigcup_{l^\prime=1}^l J_{l^\prime}. 
\end{eqnarray*}
By linearity and $\sot$-density arguments, it is sufficient to consider elements $x:=x^{(N)} \in \cA_{I}$ and
$y:= y^{(N)} \in \cA_{J}$ of a product form which is inductively defined by 
\begin{align*}
&&&&x^{(k)} &:= x_k x^{(k-1)} &\qquad &\text{with $x_k \in \cA_{I_k}$ and $x^{(k-1)} \in \cA_{I^{(k-1)}}$ },&&&&\\
&&&&y^{(l)} &:= y^{(l-1)}y_l &\qquad & \text{with $y_l \in \cA_{J_l}$ and $y^{(l-1)} \in \cA_{J^{(l-1)}}$ }.&&&&
\end{align*}
This puts us into the position to use order $\Cset$-independence in the next calculation:
\begin{eqnarray*}
\trace_{\infty}(x^{(k)}y^{(k)}) 
&=& \trace_{\infty}\big(x_k x^{(k-1)} y^{(k-1)}y_k\big)\\
&=& \trace_{\infty}\big(x_k x^{(k-1)} y^{(k-1)}\big) \trace_{\infty}( y_k) \\
&=& \trace_{\infty}(x_k) \trace_{\infty}\big( x^{(k-1)} y^{(k-1)}\big) \trace_{\infty}( y_k) 
 \end{eqnarray*}
We iterate this factorization and then, after everything is factorized, we undo it
for $x$ and $y$ separately. This gives 
\begin{eqnarray*}
\trace_{\infty}(x y ) 
= \prod_{k=1}^N\trace_{\infty}(x_k)  \prod_{l=1}^N\trace_{\infty}( y_l)
= \trace_{\infty}(x) \trace_{\infty}(y). 
 \end{eqnarray*}
Finally, the factorization properties on $\sot$-total sets of $\cA_I$ and $\cA_J$ extend linearly
to $\cA_I$ and $\cA_J$ by approximation. Doing so the proof of full $\Cset$-independence
of the random sequence $\iota^{(\alpha)}$ is completed. 

We next verify that $\iota_n(L_{\sigma_1}) = L_{\sh^n(\sigma_1)} =L_{\sigma_{n+1}}$. Indeed, this is 
obvious from (\ref{item:lrr-iv}) and the definition of the shift $\sh$. 

We are left to prove the non-spreadability of $\iota^{(\alpha)}$. For this purpose consider
the two words 
\[
w_1^{}:= \sigma_1^{} \sigma_2^{} \sigma_1^{} \sigma_2^{-1} \sigma_1^{-1}\sigma_2^{-1}
\quad \text{and}\quad 
w_2^{}:=\sigma_1^{} \sigma_3^{} \sigma_1^{} \sigma_3^{-1} \sigma_1^{-1} \sigma_3^{-1}.
\]
We note that these two words have order-equivalent index tuples 
\[
(1,2,1,2,1,2)\quad \text{and} \quad (1,3,1,3,1,3).
\]
Since the braid relations \eqref{eq:B1} and \eqref{eq:B2} imply 
\[
\trace_\infty(L_{w_1^{}})= \trace_\infty(L_{\sigma_0^{}}) = 1 \quad \text{and} \quad 
\trace_\infty(L_{w_2^{}})= \trace_\infty(L_{\sigma_1^{}\sigma_3^{-1}}) = 0,
\] 
we conclude that $\iota^{(\alpha)}$ is not spreadable. 
\end{proof}
Actually there is an abundance of non-spreadable random sequences which 
are still order $\Cset$-independent. 
\begin{Corollary}\label{cor:lrr}
Under the assumptions of Theorem \ref{thm:lrr},
let the sequence $\mathbf{\varepsilon} = (\varepsilon_k)_{k\in \Nset} \in \{1, -1\}^\Nset$ be given.
Then 
\begin{align*}
\alpha_{\mathbf{\varepsilon}} (x)&:= \sotlim_{n\to \infty} \rho(\sigma_1^{\varepsilon_1} \sigma_2^{\varepsilon_2} \cdots \sigma_n^{\varepsilon_n})(x) 
\end{align*}
defines an endomorphism for $(L(\Bset_\infty), \trace_\infty)$ such that one obtains a family of triangular towers of commuting squares, indexed by the sequence $\varepsilon$:  
\begin{eqnarray*}
\setcounter{MaxMatrixCols}{15}
\begin{matrix}
\Cset &\subset&   \Cset & \subset & L(\Bset_2) & \subset &L(\Bset_3) & \subset &L(\Bset_4) & \subset   \cdots  \subset & 
L(\Bset_\infty)\\
        &&          \cup  &         & \cup  &         & \cup &         & \cup  &       & \cup  \\
        &&  \Cset&\subset&\Cset &\subset&\alpha_\varepsilon(L(\Bset_2))&\subset&\alpha_\varepsilon(L(\Bset_3))&\subset \cdots \subset & \alpha_\varepsilon(L(\Bset_\infty))\\
         &&               &         & \cup  &         & \cup   &         & \cup  &       & \cup \\
              &&&&   \vdots &  & \vdots  &     & \vdots
              &     & \vdots 
 \end{matrix}
\setcounter{MaxMatrixCols}{10}
\end{eqnarray*}
In particular, $\iota_n^{(\alpha_\varepsilon)}:= \alpha^n_\varepsilon|_{L(\Bset_2)}$ defines a 
family, indexed by $\varepsilon$,  of stationary and order $\Cset$-independent random sequences 
\[
\iota^{(\alpha_\varepsilon)} \equiv (\iota_n^{(\alpha_\varepsilon)})_{n \in \Nset_0} \colon 
\big(L(\Bset_2), \trace_2\big) \to \big(L(\Bset_\infty), \trace_\infty\big).
\] 
In other words, $\alpha_\epsilon$ is an ordered Bernoulli shift over $\Cset$ with
generator $L(\Bset_2)$ 
(see Definition \ref{def:bernoulli}). 
\end{Corollary}
\begin{proof}
Combine Theorem \ref{thm:lrr} and Corollary \ref{cor:endo-braid-ii}.
\end{proof}
Let us now replace $\sh$ by $\sh_1$ and the Artin generators by the square roots of free generators. We can do this systematically in the following way.
Recall that the isomorphism $\inv: \Bset_\infty \to \Bset_\infty$ sends the generator $\sigma_i$ to $\sigma_i^{-1}$, 
for example: $\inv(\sigma_1^{} \sigma_2^3\sigma_3^{-1}\sigma_4^{})= \sigma_1^{-1}\sigma_2^{-3}\sigma_3^{}\sigma_4^{-1}$
and $\inv(\sigma_1^{}\sigma_2^{}\cdots \sigma_{n-1}^{}\sigma_{n}^{}) = \sigma_1^{-1}\sigma_2^{-1}\cdots \sigma_{n-1}^{-1}\sigma_{n}^{-1}$.
\begin{Theorem} \label{thm:lrr-shifted}
Consider the probability space $\big(L(\Bset_\infty), \trace_\infty\big)$, equipped with the 
1-shifted inverse representation  $\rho_1^{\inv}\colon \Bset_\infty \to \Aut{L(\Bset_\infty), \trace_\infty}$ given by
$\tau \mapsto \rho_1^{\inv}(\tau):= \Ad L_{\sh(\inv(\tau))}$. 
Then we conclude:
\begin{enumerate}
\item \label{item:lrr-shifted-i}
$\rho_1^{\inv}$ has the generating property; 
\item \label{item:lrr-shifted-ii}
$\big(L(\Bset_\infty)\big)^{\rho_1^{\inv}(\Bset_\infty)} \simeq \Cset;$
\item \label{item:lrr-shifted-iii}
$\big(L(\Bset_\infty)\big)^{\rho_1^{\inv}(\Bset_{n+2,\infty})} = L(\Bset_{n+2})$ for all $n \ge 0$;
\item \label{item:lrr-shifted-iv}
The map
\begin{eqnarray*}
\beta(x) &:=& \sotlim_{n \to \infty} \rho_1^{\inv}(\sigma_1 \sigma_2 \cdots \sigma_n)(x)\\
&=& \sotlim_{n \to \infty} \rho(\sigma^{-1}_2 \sigma^{-1}_3 \cdots \sigma^{-1}_n)(x)
\end{eqnarray*}
is an endomorphism for $(L(\Bset_\infty), \trace_\infty)$ such that $\beta(L_\tau) = L_{\sh_1(\tau)}$, with $\tau \in \Bset_\infty$;
\item \label{item:lrr-shifted-v}
Each cell of the triangular tower of inclusions is a commuting square: 
\begin{eqnarray*}
\setcounter{MaxMatrixCols}{20}
\begin{matrix}
\Cset &\subset&   L(\Bset_2)  & \subset & L(\Bset_3) & \subset & L(\Bset_4)  & \subset   \cdots  \subset &  L(\Bset_\infty)\\
        &&          \cup  &         & \cup  &              & \cup  &        & \cup  \\
        &&   \Cset &\subset& \beta(L(\Bset_2))&\subset&\beta( L(\Bset_3))&\subset \cdots \subset & \beta( L(\Bset_\infty))\\
         &&               &         & \cup  &                 & \cup  &        & \cup \\
              &&&&   \Cset& \subset &  \beta^2(L(\Bset_2))
              & \subset   \cdots \subset & \beta^2( L(\Bset_\infty)) \\
       &&&&&&   \cup           &         & \cup          \\
        &&&&&&   \vdots           &         & \vdots        
\end{matrix}
\setcounter{MaxMatrixCols}{10}
\end{eqnarray*}
\item \label{item:lrr-shifted-vi}
The maps $\iota_n^{(\beta)}:= \beta^n|_{L(\Bset_2)}$ define a 
braidable, spreadable and full $\Cset$-independent random sequence 
\[
\iota^{(\beta)} \equiv (\iota_n^{(\beta)})_{n \in \Nset_0} \colon \big(L(\Bset_2), \trace_2\big) \to \big(L(\Bset_\infty), \trace_\infty\big).
\] 
In particular, $\beta$ is a full Bernoulli shift over $\Cset$ with
generator $L(\Bset_2)$ 
(see Definition \ref{def:bernoulli}). 
\item \label{item:lrr-shifted-vii}
The random sequence $\iota^{(\beta)}$ is not exchangeable.   
\end{enumerate}
\end{Theorem}
\begin{Remark} \label{rem:lrr-shifted} \normalfont
For all $n \in \Nset_0$, 
\[
\iota_n^{(\beta)}(L_{\gamma_1}) = L_{\sh_1^n(\gamma_1)} =L_{\gamma_{n+1}}.
\]
Hence we have the very remarkable fact that the roots of free generators form 
a spreadable random sequence and behave better in this respect than the usual 
Artin generators.
\end{Remark} 
\begin{proof}
It is easy to see that the state $\trace_\infty$ on $L(\Bset)$ is $\rho_1^{\inv}(\tau)$-invariant for every $\tau \in \Bset_\infty$.
Thus we have a representation $\rho_1^{\inv}\colon \Bset_\infty \to \Aut{L(\Bset_\infty, \trace_\infty)}$. 

(\ref{item:lrr-shifted-i}) The generating property of $\rho^{\inv}_1$ follows from 
\begin{eqnarray*}
L(\Bset_\infty) &\supset& \bigvee_{n \ge 0} \big(L(\Bset_\infty)\big)^{\rho_1^{\inv}(\Bset_{n+2,\infty})}
                         = \bigvee_{n \ge 0} \big(L(\Bset_\infty)\big)^{\rho^{\inv}(\Bset_{n+3,\infty})}\\  
                         &= &\bigvee_{n \ge 0} \big(L(\Bset_\infty)\big)^{\rho(\Bset_{n+3,\infty})}
                         \supset \bigvee_{n \ge 0} L(\Bset_{n+2}) = L(\Bset_\infty).
\end{eqnarray*}  
At this point we have verified all assumptions of Theorem \ref{thm:endo-braid-i}. Since
\[
\big(L(\Bset_\infty)\big)^{\rho_1^{\inv}(\Bset_{n+2,\infty})} = L(\Bset_\infty)^{\rho(\Bset_{n+3,\infty})},
\]
(\ref{item:lrr-shifted-ii}) and (\ref{item:lrr-shifted-iii}) are directly concluded from 
Theorem \ref{thm:irreducibility} and \ref{thm:relative-commutants}.

(\ref{item:lrr-shifted-iv})
The action of the endomorphism $\beta$ on $L_{\tau}$ is identified by   Lemma \ref{lem:shift_1} as
\begin{eqnarray*}
\beta(L_{\tau}) &=& L_{\sigma_2^{-1} \sigma_3^{-1} \cdots \sigma_{n}^{-1}}L_{\tau} L_{\sigma_n^{} \cdots \sigma_{3}^{}\sigma_{2}^{}}
= L_{\sh_1(\tau)}
\end{eqnarray*} 
(with $n$ sufficiently large). 
(\ref{item:lrr-shifted-v}) is immediate from  (\ref{item:lrr-shifted-ii}) and Theorem \ref{thm:endo-braid-i}. 

(\ref{item:lrr-shifted-vi}) 
Braidability of $\iota^{(\beta)}$ follows by definition (with the representation $\rho_1^{\inv}$) and
spreadability is then a consequence of Theorem \ref{thm:sequences-braid}.
Using (\ref{item:lrr-shifted-iv}) and $\sh_1(\gamma_i) = \gamma_{i+1}$
(Lemma \ref{lem:shift_2}) the image of $L_{\gamma_1}$ under 
$\iota^{(\beta)}_n$ is identified to be $\gamma_{n+1}$.  
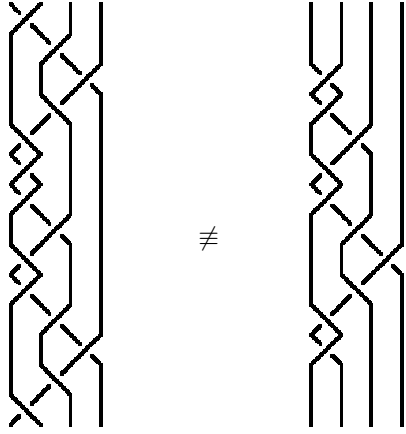
\begin{figure}[h]
\setlength{\unitlength}{0.2mm}
\begin{picture}(280,290)
\savebox{\artin}(20,20)[1]{\masterartin} 
\savebox{\artininv}(20,20)[1]{\masterartininv} 
\savebox{\strandr}(20,20)[1]{\masterstrandr} 
\savebox{\strandl}(20,20)[1]{\masterstrandl} 
\savebox{\horizontaldots}(20,20)[1]{\masterhorizontaldots}
\put(0,260){\usebox{\artin}}  
\put(20,260){\usebox{\strandr}}
\put(40,260){\usebox{\strandr}}
\put(0,240){\usebox{\strandl}}
\put(20,240){\usebox{\artin}}
\put(40,240){\usebox{\strandr}}
\put(0,220){\usebox{\strandl}}
\put(20,220){\usebox{\strandl}}
\put(40,220){\usebox{\artin}}
\put(0,200){\usebox{\strandl}}
\put(20,200){\usebox{\artininv}}
\put(40,200){\usebox{\strandr}}
\put(0,180){\usebox{\artininv}}
\put(20,180){\usebox{\strandr}}
\put(40,180){\usebox{\strandr}}
\put(0,160){\usebox{\artin}}
\put(20,160){\usebox{\strandr}}
\put(40,160){\usebox{\strandr}}
\put(0,140){\usebox{\artin}}  
\put(20,140){\usebox{\strandr}}
\put(40,140){\usebox{\strandr}}
\put(0,120){\usebox{\strandl}}
\put(20,120){\usebox{\artin}}
\put(40,120){\usebox{\strandr}}
\put(0,100){\usebox{\artininv}}
\put(20,100){\usebox{\strandr}}
\put(40,100){\usebox{\strandr}}
\put(0,80){\usebox{\artin}}  
\put(20,80){\usebox{\strandr}}
\put(40,80){\usebox{\strandr}}
\put(0,60){\usebox{\strandl}}
\put(20,60){\usebox{\artin}}
\put(40,60){\usebox{\strandr}}
\put(0,40){\usebox{\strandl}}
\put(20,40){\usebox{\strandl}}
\put(40,40){\usebox{\artin}}
\put(0,20){\usebox{\strandl}}
\put(20,20){\usebox{\artininv}}
\put(40,20){\usebox{\strandr}}
\put(0,0){\usebox{\artininv}}
\put(20,0){\usebox{\strandr}}
\put(40,0){\usebox{\strandr}}
\put(125,120){$\not\equiv$}
\put(200,260){\usebox{\strandl}}
\put(200,260){\usebox{\strandr}}
\put(220,260){\usebox{\strandr}}
\put(240,260){\usebox{\strandr}}
\put(200,240){\usebox{\strandl}}
\put(200,240){\usebox{\strandr}}
\put(220,240){\usebox{\strandr}}
\put(240,240){\usebox{\strandr}}
\put(200,220){\usebox{\artin}}
\put(220,220){\usebox{\strandr}}
\put(240,220){\usebox{\strandr}}
\put(200,200){\usebox{\artin}}
\put(220,200){\usebox{\strandr}}
\put(240,200){\usebox{\strandr}}
\put(200,180){\usebox{\strandl}}
\put(220,180){\usebox{\artin}}
\put(240,180){\usebox{\strandr}}
\put(200,160){\usebox{\artininv}}
\put(220,160){\usebox{\strandr}}
\put(240,160){\usebox{\strandr}}
\put(200,140){\usebox{\artin}}
\put(220,140){\usebox{\strandr}}
\put(240,140){\usebox{\strandr}}
\put(200,120){\usebox{\strandl}}
\put(220,120){\usebox{\artin}}
\put(240,120){\usebox{\strandr}}
\put(200,100){\usebox{\strandl}}
\put(220,100){\usebox{\strandl}}
\put(240,100){\usebox{\artin}}
\put(200,80){\usebox{\strandl}}
\put(220,80){\usebox{\artininv}}
\put(240,80){\usebox{\strandr}}
\put(200,60){\usebox{\artininv}}
\put(220,60){\usebox{\strandr}}
\put(240,60){\usebox{\strandr}}
\put(200,40){\usebox{\artin}}
\put(220,40){\usebox{\strandr}}
\put(240,40){\usebox{\strandr}}
\put(200,20){\usebox{\strandl}}
\put(200,20){\usebox{\strandr}}
\put(220,20){\usebox{\strandr}}
\put(240,20){\usebox{\strandr}}
\put(200,00){\usebox{\strandl}}
\put(200,00){\usebox{\strandr}}
\put(220,00){\usebox{\strandr}}
\put(240,00){\usebox{\strandr}}
\end{picture}
\caption{Four-strand braids $\gamma_3 \gamma_1 \gamma_2 \gamma_3$ (left) and 
$\gamma_1 \gamma_2 \gamma_3 \gamma_1$ (right)}
\label{figure:remove-strands-i}
\end{figure}

(\ref{item:lrr-shifted-vii}) We need to show the non-exchangeability of $\iota^{(\beta)}$. 
From (\ref{eq:EB}) (see Proposition \ref{thm:sqrt-rep}) we have
\[
\gamma_3 \gamma_2 \gamma_1 \gamma_3 
= \gamma_2 \gamma_1 \gamma_3 \gamma_2,
\]
but if we interchange the subscripts $1$ and $2$ then
\[
\gamma_3 \gamma_1 \gamma_2 \gamma_3 
\not= \gamma_1 \gamma_2 \gamma_3 \gamma_1.
\]
In fact, if in the geometric picture we represent the two words
by four-strands braids (see Figure \ref{figure:remove-strands-i}) and 
then remove the third and fourth strand to obtain braids in $\Bset_2$ 
then the left hand side yields the identity $\sigma_0$
but the right hand side yields $\sigma^2_1$. It follows that
\begin{gather*}
\trace_4(L_{\gamma_3} L_{\gamma_2} L_{\gamma_1} L_{\gamma_3} 
L_{\gamma^{-1}_2} L_{\gamma^{-1}_3} L_{\gamma^{-1}_1} L_{\gamma^{-1}_2}) 
= 1 \\
\not= 0 = 
\trace_4(L_{\gamma_3} L_{\gamma_1} L_{\gamma_2} L_{\gamma_3} 
L_{\gamma^{-1}_1} L_{\gamma^{-1}_3} L_{\gamma^{-1}_2} L_{\gamma^{-1}_1}) 
\end{gather*}
and hence $\iota^{(\beta)}$ is not exchangeable.
\end{proof}
\begin{Definition}\normalfont
Let $0 < r < t \le \infty$. We denote by $\Fset_{r,\,t}^{1/2}$ the subgroup 
of $\Bset_\infty$ generated by the square root of free generators 
$\set{\gamma_s}{r \le s \le t}$ and by $\Fset_{r,\,t}$ the subgroup of 
$\Bset_\infty$ with the free generators $\set{\gamma_s^2}{r \le s \le t}$.
\end{Definition}
Clearly $\Fset_{r,\,t} \subset \Fset_{r,\,t}^{1/2}$.
We have $\Fset_{1,\,n}^{1/2} = \Bset_{n+1}$ and $\Fset_{1,\,n} \simeq \Fset_n$, where $\Fset_n$ denotes
the group in $n$ free generators, and hence $\sh_1^k(\Fset_{s,t}^{1/2})= \Fset_{s+k,t+k}^{1/2}$, as well as
$\sh_1^k(\Fset_{s,t})= \Fset_{s+k,t+k}$. 
In the rest of this section we assume that the reader is familiar with some notions in free probability as they have been introduced by Voiculescu. As a reminder see for example \cite{VDN92a}.
\begin{Corollary}\label{cor:srfg}
The square root of free generator presentation $\langle \gamma_i | i \in \Nset \rangle$ of $\Bset_\infty$ gives 
rise to the system of Haar unitaries $\set{L_{\gamma_i}}{i \in \Nset}$ such that 
\[
L(\Fset^{1/2}_{r,\,t}) = \bigvee\set{L_{\gamma_s}}{r \le s \le t}
\] 
and such that one has  the following triangular tower of commuting squares: 
\begin{eqnarray*}
\setcounter{MaxMatrixCols}{20}
\begin{matrix}
\Cset &\subset&   L(\Fset_{1,\,1}^{1/2})  & \subset & L(\Fset_{1,\,2}^{1/2})  & \subset  & 
L(\Fset_{1,\,3}^{1/2})  & \subset &L(\Fset_{1,\,4}^{1/2})  & \subset \cdots   \subset &L(\Fset_{1,\,\infty}^{1/2}) \\
               &&          \cup &&    \cup &&  \cup && \cup && \cup\\
  &&   \Cset & \subset &L(\Fset_{2,\,2}^{1/2})  & \subset &L(\Fset_{2,\,3}^{1/2}) &\subset  & L(\Fset_{2,\,4}^{1/2}) &\subset\cdots  \subset & L(\Fset_{2,\,\infty}^{1/2})\\
        &&            &               & \cup  &  & \cup &       & \cup && \cup \\
        &&   && \Cset &\subset&  L(\Fset_{3,\,3}^{1/2})  &\subset &   L(\Fset_{3,\,4}^{1/2})  &\subset \cdots  \subset &  L(\Fset_{3,\,\infty}^{1/2}) \\
       &&&&               &             & \cup  && \cup && \cup\\
        &&&&              &             & \vdots && \vdots  && \vdots
\end{matrix}
\setcounter{MaxMatrixCols}{10}
\end{eqnarray*}
The squared family $\set{L_{\gamma_i}^2}{i \in \Nset}$ is a free system of Haar unitaries in the sense of Voiculescu 
whose generated triangular tower is a restriction of the above triangular tower:
\begin{eqnarray*}
\setcounter{MaxMatrixCols}{20}
\begin{matrix}
\Cset &\subset&   L(\Fset_{1,\,1})  & \subset & L(\Fset_{1,\,2})  & \subset  & 
L(\Fset_{1,\,3})  & \subset &L(\Fset_{1,\,4})  & \subset \cdots   \subset &L(\Fset_{1,\,\infty}) \\
               &&          \cup &&    \cup &&  \cup && \cup && \cup\\
  &&   \Cset & \subset &L(\Fset_{2,\,2})  & \subset &L(\Fset_{2,\,3}) &\subset  & L(\Fset_{2,\,4}) &\subset\cdots  \subset & L(\Fset_{2,\,\infty})\\
        &&            &               & \cup  &  & \cup &       & \cup && \cup \\
        &&   && \Cset &\subset&  L(\Fset_{3,\,3})  &\subset &   L(\Fset_{3,\,4})  &\subset \cdots  \subset &  L(\Fset_{3,\,\infty}) \\
       &&&&               &             & \cup  && \cup && \cup\\
        &&&&              &             & \vdots && \vdots  && \vdots
\end{matrix}
\setcounter{MaxMatrixCols}{10}
\end{eqnarray*} 
Moreover, each cell in this tower forms a commuting square 
\[
\begin{matrix}C &\subset& D \\\cup && \cup  \\ A &\subset& B\end{matrix},
\] 
such that  $B$ and $C$ are freely independent with amalgamation over $A$. 
\end{Corollary}
\begin{proof}
It is immediate from the definition of the left regular representation that
$L_{\gamma_i}$ is a Haar unitary, i.e., $\trace_\infty(L_{\gamma_i}^n)=0$
for all $n \in \Zset \setminus \{0\}$. 
From Theorem \ref{thm:lrr-shifted} we get the commuting squares for the
$L(\Fset_{s,\,t}^{1/2})$'s and, by restriction, also for the 
$L(\Fset_{s,\,t})$'s. This is independence in the sense of Definition
\ref{def:independence}, but for the squared family more is true: It is a basic result of free probability theory that they are (amalgamated) free. See \cite{VDN92a}, Examples 2.5.8 and 3.8.3.
\end{proof}
It is easily seen that freeness with amalgamation implies conditional independence in the sense of Definition \ref{def:independence}.  The converse fails to be true since our notion of independence is more general. But it would be of interest to determine combinatorial formulas for the mixed moments of square root random sequences such that these formulas extend the combinatorics of noncrossing partitions from free probability.   

The close connection between the square root presentation $\Fset_\infty^{1/2}$ 
and the free group $\Fset_\infty$ invites to ask whether there exists a deeper
parallel between objects considered in free probability theory and the appropriately chosen `square root objects' in a `braided probability theory'. We are going to illustrate this in an example.

Since $L_{\gamma_1}$ and $L_{\gamma_1^2}$ are Haar unitaries, the selfadjoint operators
$L_{\gamma_i}+ L_{\gamma_i}^*$ and $L_{\gamma_i^2}+ L_{\gamma_i^2}^*$ have both the arcsine law on the interval $[-2,2]$
as spectral distribution (see \cite[Lecture 1]{NiSp06a}). But this can be taken one step further:

According to Theorem \ref{thm:main-1} we may consider any von Neumann subalgebra 
$\cC_0$ of $L(\Fset_{1,\,2}^{1/2})$ and the restriction of the random sequence to $\cC_0$ 
will again be a spreadable random sequence. Among the interesting choices is of course 
\[
\cC_0^{1/2} := \bigvee\{L_{\gamma_1} + L_{\gamma_1}^* \}  
\]
This gives the `random sequence' $(\beta^n(\cC_0^{1/2}))_n  \subset L(\Bset_\infty)$ 
which can be understood as a braided counterpart of a free sequence. 
If we put $\cC_0:=\bigvee\{(L_{\gamma_1} + L_{\gamma_1}^*)^2\}$ and observe 
that $(L_{\gamma_1} + L_{\gamma_1}^*)^2 = 2 + L_{\gamma_1^2} + L_{\gamma_1^2}^*$, so the 
`squared random sequence' $(\beta^n(\cC_0))_n\subset L(\Fset_\infty)$ is free.  
\begin{Conjecture}
There exists a braided extension of free probability.
\end{Conjecture}
Such an extension should, of course, be related to some
`braided independence' as a \emph{specific} form of $\Cset$-independence (in the sense of Definition \ref{def:independence}), but it must necessarily lie \emph{beyond} independence with universality rules (in the sense of Speicher) which is completely classified and leaves no room for such an idea, see \cite{Spei97a,BeSc02a, NiSp06a}). 
   
At the moment we have no definite formulation of such a theory but we suggest that there are interesting concrete problems on the way. For example,   
in view of Corollary \ref{cor:srfg},
it is intriguing to ask whether the combinatorics of free probability theory
can be appropriately extended to a combinatorics of free square root presentations. A promising starting point for such investigations are random walks on free square root
presentations, in parallel to Kesten's work on symmetric random walks on groups \cite{Kest59a}. For example
Kesten determined for the symmetric random walk on $\Fset_2$ the spectral distribution of the 
free Laplacian 
\[
\Delta_{\operatorname{free}} = \frac{1}{2}\big(L_{\gamma_1^2} + L_{\gamma_1^2}^* + L_{\gamma_2^2} + L_{\gamma_2^2}^*)
\] 
with respect to the trace $\trace_3$. It is determined by the moment generating function 
\[
\sum_{n=0}^\infty \trace_3(\Delta_{\operatorname{free}}^n)z^n = \frac{2\sqrt{1-12z^2}-1}{1-16z^2}
\]
which can effectively be determined with freeness (see \cite[Lecture 4]{NiSp06a}).

The combinatorics involved for determining the $n$-th moment or the moment generating function of the corresponding braided Laplacian 
\[
\Delta_{\operatorname{braid}}= \frac{1}{2}\big(L_{\gamma_1} + L_{\gamma_1}^* + L_{\gamma_2} + L_{\gamma_2}^*)
\] 
amounts to answer the following question on 3-strand braids: Consider all words $w$ of length $n$ written in the alphabet of 4 letters $\gamma_1^{\pm 1}, \gamma_2^{\pm 1}$. How many words $w$ among the $4^n$-words describe the 
trivial 3-strand braid? An answer to this question immediately gives  the spectral measure of $\Delta_{\operatorname{braid}}$ in terms of moments.
 
Related explicit calculations for randomly growing braids on three strands are contained in \cite{MaMa07a} and involve random Garside normal forms for
the Artin presentation of $\Bset_3$. Unfortunately this approach does not generalize for $\Bset_n$ with $n \ge 4$. It would be of interest to investigate if the square root of free generator presentation, which has no Garside structure for $n \ge 4$ \cite{Bir08a,Deho08a} but coincides for $n=3$ with the Artin presentation, gives an alternative approach to such problems (see also Section \ref{section:presentation}). 

Further background information and additional structures may come from the closeness of our approach to subfactor theory and the recent progress on the connection between subfactors, large random matrices and free probability theory (see \cite{GJS07aPP} and references therein).

\section{Some concrete examples}
\label{section:examples}
In the following we discuss a few concrete examples which are well known 
but which can be looked upon from a new point of view by integrating them 
into our theory of spreadable random sequences from braid group 
representations. It appears that this strategy allows to simplify some 
arguments. Of course there are many other examples. 
\begin{Example}[\emph{Gaussian Representations}] \normalfont 
\label{example:2} 
Choose $2 \leq p \in \Nset$ and
\[
\omega :=
\left\{\begin{array}{cl} exp(2 \pi \i/p) & \mbox{ if p is odd}; \\  
exp(\pi \i/p) & \mbox{ if p is even}. 
\end{array} \right.
\]
Then consider unitaries $(e_i)_{i\in\Nset_0}$ satisfying 
\begin{align*}
&&&& e^p_i &= \1 &&\text{for all $i$};&&\\ 
&&&& e_i e_j &= \omega^2 e_j e_i &&\text{whenever $i<j$}.&& 
\end{align*}
A pair $e_i, e_j$ with $i<j$ can be realized in (and generates) the 
$(p\times p)-$matrices and, taking the weak closure with respect to the 
trace $\trace$, the sequence $(e_i)_{i\in\Nset_0}$ generates the hyperfinite
$II_1$-factor. This is the noncommutative probability space in this class 
of examples.
\\

If we now define $(v_i)_{i\in\Nset}$ by $v_i := \omega e^*_{i-1} e_i$ 
(for all $i \in \Nset$) then
\begin{eqnarray*}
v^p_i &=& \1; \\
v_i v_{i+1} &=& \omega^2 v_{i+1} v_i;\\
v_i v_j &=& v_j v_i \quad \text{if $|i-j|>1$}.
\end{eqnarray*}
We remark that $k \mapsto \omega^{k^2} v^k_i$ is well defined for 
$k$ mod $p$ and so the following sums can always be interpreted as sums 
over the cyclic group $\Zset_p$. Then a direct computation shows that 
$(u_i)_{i\in\Nset}$ defined by
\[
u_i := \frac{1}{\sqrt{p}} \sum^{p-1}_{k=0} \omega^{k^2} v^k_i
\]
are unitary and satisfy the braid relations \eqref{eq:B1} and 
\eqref{eq:B2}. Hence we obtain a unitary representation of
$\Bset_\infty$, called the Gaussian representation in 
\cite[Subsection 5.8]{Jone91a} because of the Gaussian sums in related
computations. As usual, let us consider the endomorphism
\[
\alpha = \lim_{n \to \infty} \Ad(u_1u_2 \cdots u_n).
\]
This endomorphism has already been studied in \cite{Rupp95a} as an
example of a noncommutative Bernoulli shift. Let us summarize what we can say about it from the point of view of our theory. 
\\

From the relations 
\begin{eqnarray*}
e_i v_j =
\begin{cases}
\omega^2 v_j e_i & \text{if $j=i$ or $j=i+1$},\\ 
v_j e_i &  \text{otherwise}, 
\end{cases}
\end{eqnarray*}
it is readily checked that this is 
a product representation for $\alpha$ (satisfying 
(PR-1), (PR-2) of 
Definition \ref{def:adapted-end}) with respect to the tower 
$(\cM_n)^{\infty}_{n=0}$ with $\cM_n$ generated by $e_0,\ldots,e_n$. 
More precisely we have $\alpha(e_i) = e_{i+1}$ for all $i \in \Nset_0$. 
In fact,
\begin{eqnarray*}
u_1 e_0 u^*_1 &=& \frac{1}{p} \sum_{k,k'} \omega^{(k^2-k'^2)} v^k_1 e_0 v^{-k'}_1\\
&=& \frac{1}{p} \sum_{k,k'} \omega^{(k^2-k'^2)} \omega^{-2k} e_0 v^{(k-k')}_1 \\
&=& \sum_\ell \omega^{-\ell^2} e_0 v^\ell_1 \big( \frac{1}{p} \sum_k \omega^{2k(\ell-1)} \big)\\
&=& \omega^{-1} e_0 v_1 = e_1, \quad \mbox{etc.}
\end{eqnarray*}
Iterated application of $\alpha$ to $\cM_0$ produces a braidable sequence 
which by Theorem \ref{thm:main-1} is spreadable. For the sequence 
$e_0, e_1, e_2, \ldots$ this follows also more directly from the fact
that order preserving transformations preserve the commutation relations for 
the $e_n$. 

\begin{Proposition}\label{gauss}
The endomorphism $\alpha$ is a 
full Bernoulli shift (in the sense of Definition \ref{def:bernoulli}) over
\[
\cM_{-1} = \cM^\alpha = \cM^{\tail} = \cM^{\Bset_\infty} = \Cset
\]
with generator
\[
\cM_0 = \linh\{e^k_0\}^{p-1}_{k=0} \simeq \Cset^p. 
\]
Further we have for all $n \ge -1$
\[
\cM_n = \cM \cap \{u_k \colon k \ge n+2\}^\prime.
\]
\end{Proposition}
\begin{proof}
We conclude by Theorem \ref{thm:main-3} that $\alpha$ is a 
full Bernoulli shift over the fixed point algebra $\cM^{\Bset_\infty}$
with generator $\cM_0$. 
We can also check, by direct computation, the commuting squares assumptions of Theorem \ref{thm:endo-pr} with $\cM_{-1}=\Cset$.
Now everything follows from Theorem \ref{thm:endo-pr}
together with 
Theorem \ref{thm:tower-reconstruction} for $\cM_n$. 
\end{proof}

Note that for $p=2$ we have a Clifford algebra with anticommuting $e_n$'s
and we can check that the sequence $e_0, e_1, e_2, \ldots$ is exchangeable. 
This is no longer the case for $p>2$. For example
\[
\trace(e_1 e_2 e^*_1 e^*_2) 
= \omega^2 \not= \omega^{-2} 
= \trace(e_2 e_1 e^*_2 e^*_1).
\]
Further results about more general product representations with respect to 
this tower can be found in \cite{Gohm01a}. 
\end{Example}
\begin{Example}[\emph{Hecke algebras}]\normalfont \label{example:3}
Recall from \cite[Example 3.1]{Jone94a}: The Hecke algebra over $\Cset$ with
parameter $q\in \Cset$ is the unital algebra with generators
$g_0,g_1,\ldots$ and relations
\begin{eqnarray*}
g^2_n &=& (q-1) g_n + q\\ 
g_m g_n &=& g_n g_m  \qquad \text{if $\mid n-m \mid \geq 2$}\\
g_n g_{n+1} g_n &=& g_{n+1} g_n g_{n+1}
\end{eqnarray*}
Then if $q$ is a root of unity it is possible to define an involution
and a trace such that the $g_n$ are unitary, and the tower
$(\cM_n)^{\infty}_{n=0}$ with $\cM_n$ generated by
$g_0,\ldots,g_n$ is embedded into the hyperfinite $II_1$-factor $\cM$
with its trace. In \cite{Jone94a} the commuting square assumptions
required for our Theorem \ref{thm:endo-pr} are checked for the $\Ad(g_n)$'s 
with $\cM_{-1} = \Cset$. As in the previous example we conclude that
\[
\alpha := \lim_{n\to\infty} \Ad(g_1 \cdots g_n)
\]
defines a full Bernoulli shift over $\Cset$ with generator $\cM_0$, that  
\[
\cM^\alpha = \cM^{\tail} = \cM^{\Bset_\infty} = \Cset
\]
by Theorem \ref{thm:main-3}, and that 
\[
\cM_n = \cM \cap \{g_k \colon k \ge n+2\}^\prime
\]
by Theorem \ref{thm:tower-reconstruction}. Such a situation occurs in 
particular in Jones' subfactor theory \cite{Jone83a} if the index 
belongs to the discrete range. 
\end{Example}
\begin{Example}[\emph{$R$-matrices}]\normalfont \label{example:4} 
Now take the tower 
\[
\cM_0 :=M_p \subset \cM_1 := M_p \otimes M_p \subset 
\cM_2 := M_p \otimes M_p \otimes M_p \subset 
\cdots \subset \cM
\]
where $M_p$ denotes the $(p\times p)$-matrices and the embeddings are given by
\mbox{$X \mapsto X \otimes \1$}, 
all sitting inside the hyperfinite $II_1$-factor $\cM$ with trace $\trace$. 
By an $R$-matrix we mean an element $\breve{R}$ of
$M_p \otimes M_p$ satisfying the (constant quantum) Yang-Baxter equation (YBE)
\[
\breve{R}_{12} \breve{R}_{23} \breve{R}_{12} 
= \breve{R}_{23} \breve{R}_{12} \breve{R}_{23}
\]  
where we use a leg notation, i.e., the subscripts describe 
the embedding of $\breve{R}$ into a triple tensor
product of $M_p$'s. See \cite{Jone91a} for the role of $R$-matrices in providing interesting braid group representations. Note that with
$R = P \breve{R}$, where $P$ is the flip operator, we get another familiar form of the YBE
\[
R_{12} R_{13} R_{23} = R_{23} R_{13} R_{12}
\]
which plays an important role in the theory of quantum groups. 
For an overview see \cite{Kas95,Maji95a}.

If $\breve{R}$ is a unitary $R$-matrix then evidently $(u_n)^\infty_{n=1}$ with
$u_n := \breve{R}_{n-1,n}$ satisfy the braid relations and provide us with a
product representation for an adapted endomorphism
$\alpha := \lim_{n\to\infty} \Ad(u_1 \cdots u_n)$. By Theorem 
\ref{thm:main-3} it is (or restricts to) a full Bernoulli shift over 
$\cM^\alpha = \cM^{\tail} = \cM^{\Bset_\infty}$ with generator 
$\cM_0 \vee \cM^{\Bset_\infty}$ 
which produces a braidable random sequence by iterated applications to
$\cM_0$. 
\\
For example we can take $p=2$ and (with respect to a basis
$\delta_0 \otimes \delta_0, \delta_0 \otimes \delta_1,
\delta_1 \otimes \delta_0, \delta_1 \otimes \delta_1$)
\[
\breve{R} =     \begin{pmatrix}
        1 & 0  & 0  & 0  \\ 
        0  & 0  & 1 &  0 \\ 
        0 & 1 & 0  &  0 \\
        0  & 0  &  0 & \omega  
        \end{pmatrix} 
\]
with $|\omega | = 1$. YBE is easily checked directly. 
By computing the first commuting square in Theorem \ref{thm:endo-pr} 
we find that $\cM_{-1}=\Cset$.
This full Bernoulli shift over $\Cset$ with generator $\cM_0$ is called the $\omega$-shift in \cite{Rupp95a}.
The corresponding subfactors have been investigated in \cite{KrSuVa96a}.
Note that these examples include the
usual tensor shift ($\omega = 1$) and the CAR-shift ($\omega = -1$)
in its Jordan-Wigner form \cite{BrRo2}.
\end{Example}
\begin{Example}[\emph{$R$-matrices, non-homogeneous case}]\normalfont 
\label{example:5} 
Slightly varying the construction in Example \ref{example:4} provides 
us with very elementary examples of non-spreadable ordered Bernoulli shifts 
over $\Cset$. Define 
\[
u_n =
        \begin{pmatrix}
        1 &  0 & 0  & 0  \\ 
        0  &  0 & 1 &  0 \\ 
        0  & 1 & 0  &  0 \\
        0 &  0 &  0 & \omega_n  
        \end{pmatrix},
\]
embedded at positions $n-1$ and $n$ as before, and again
$\alpha := \lim_{n\to\infty} \Ad(u_1 \cdots u_n)$. If $\omega_n$
depends on $n$ then $\alpha$ continues to be an ordered Bernoulli shift over 
$\Cset$ (by Theorem \ref{thm:endo-pr}) but may fail to produce spreadable 
random sequences. Take for example $\omega_1 = 1$ and $\omega_2 = -1$, i.e.,
we mix the tensor shift and the CAR-shift. Then for 
\[
\cA_0 \ni x =
        \begin{pmatrix}
        0 & 1 \\ 
        1 & 0 \\  
        \end{pmatrix} 
\] 
we find:
\begin{eqnarray*}
\alpha(x) &=& \Ad(u_1)(x) = \1 \otimes x,\\
\alpha^2(x) &=& \Ad(u_1 u_2)(x) = 
\begin{pmatrix}
        1 &  0 \\ 
        0 & -1   
\end{pmatrix}
\otimes \1 \otimes x,\\
{[x \, \alpha(x)]}^2 &=& [ x \otimes x ]^2 = \1,\\
{[ x\, \alpha^2(x) ]}^2 &=& {\left[ 
   \begin{pmatrix}
        0 & -1 \\ 
        1 &  0 \\  
   \end{pmatrix} 
\otimes \1 \otimes x \right]}^2 = -\1.
\end{eqnarray*}
Therefore $\trace(x\, \alpha(x)\, x\, \alpha(x)) \not=
\trace(x\, \alpha^2(x)\, x\, \alpha^2(x))$, which shows that
already the mixed moments of 4th order are not invariant
under order equivalence.
\end{Example}
\begin{appendix}
\section{Operator algebraic noncommutative probability}
\label{section:appendix}
Here we present in a condensed form some well known and some less well known 
constructions of noncommutative probability theory in the framework of von 
Neumann algebras. This provides a general context for the results of this paper. Similar settings are considered in \cite{Kuem85a,Kuem88a,HKK04a,An06a,JuPa08a,Koes08aPP}.
The last part about connections between product representations 
and Bernoulli shifts (\ref{thm:endo-pr} - \ref{thm:tower-reconstruction}) is new in this form.  
\begin{Definition}\normalfont \label{def:nc-prob-space}
A \emph{(noncommutative) probability space} $(\cM, \psi)$ consists of a von Neumann algebra $\cM$ with separable predual and a faithful normal state $\psi$ on $\cM$. A \emph{$\psi$-conditioned} probability space $(\cM, \psi ,\cM_0)$ consists of a probability space $(\cM, \psi)$ and a von Neumann subalgebra $\cM_0$ of $\cM$ such that the $\psi$-preserving conditional expectation $E_{\cM_0}$ from $\cM$ onto $\cM_0$ exists. As an abbreviation we also say in this case that $\cM_0$ is a \emph{$\psi$-conditioned subalgebra}. 
\end{Definition}
Note that probability spaces always come with a standard representation on a separable Hilbert space via GNS construction for $\psi$. By Takesaki's theorem (see \cite{Take03b}), the  $\psi$-preserving conditional expectation $E_{\cM_0}$ exists if and only if $\sigma_t^\psi(\cM_0) = \cM_0$ for all $t \in \Rset$, where $\sigma^\psi$ is the modular automorphism group associated to $(\cM, \psi)$. Thus the existence of such a conditional expectation is automatic if $\psi$ is a trace.
\begin{Lemma} \label{lemma:automorphism}
A $\psi$-preserving automorphism $\alpha$ always commutes with the modular automorphism group $\sigma^\psi$.
\end{Lemma}
\begin{proof}
The modular automorphism group $\sigma^\psi$ is uniquely characterized by 
the KMS-condition with respect to $\psi$, see \cite{Take03b}. Hence it is 
enough to show that $\alpha^{-1} \sigma_t^\psi \alpha$ also satisfies
the KMS-condition with respect to $\psi$. Using $\psi \circ \alpha = \psi$ 
this is easily done.
\end{proof} 
We denote the set of $\psi$-preserving automorphisms $\alpha$ of $\cM$ by 
$\Aut{\cM,\psi}$. Combined with Takesaki's theorem we conclude that the 
following subalgebras are always $\psi$-conditioned:
\begin{itemize}
\item fixed point algebras of a $\psi$-preserving automorphism \\
(or of a set of such automorphisms)
\item the image of a $\psi$-conditioned algebra under a $\psi$-preserving automorphism
\item algebras generated by $\psi$-conditioned algebras
\end{itemize}
By combining these items all the subalgebras in this paper turn out to be $\psi$-conditioned. 

We also encounter a special type of (non-surjective) $\psi$-preserving endomorphisms. In the probability space $(\cM, \psi)$ consider a tower
\[
\cM_0 \subset \cM_1 \subset \cM_2 \subset \cdots
\]
of $\psi$-conditioned subalgebras such that $\cM = 
\bigvee_{n\in\Nset_0} \cM_n$, and a family of automorphisms
$(\alpha_k)_{k\in\Nset} \subset \Aut{\cM,\psi}$ satisfying (for all $n\in\Nset_0$)
\begin{align}
 \alpha_k\, (\cM_n) = \cM_n & \quad \text{if}\; k \le n
\tag{PR-1}  \label{item:pr-i}\\
 \alpha_k |_{\cM_n} = \id |_{\cM_n} & 
 \quad \text{if}\; k \ge n+2
\tag{PR-2} \label{item:pr-ii}
\end{align}
\begin{Definition}\normalfont \label{def:adapted-end}
Given (\ref{item:pr-i}) and (\ref{item:pr-ii}) then
\begin{align}
\alpha&= \lim_{n\to \infty}\alpha_1\alpha_2 \alpha_3\cdots \alpha_n, \tag{PR-0} \label{item:pr-0}
\end{align}
(in the pointwise strong operator topology)
defines a $\psi$-preserving endomorphism of $\cM$ which we call an \emph{adapted endomorphism with product representation}.
\end{Definition}
So PR stands for product representation. The existence of the limit is easily deduced from (\ref{item:pr-i}) and (\ref{item:pr-ii}), in fact for $n \in \Nset$ and $x \in \cM_{n-1}$
\[
\alpha(x) = \alpha_1 \cdots \alpha_n(x),
\]
a finite product. From the limit formula it also follows that $\alpha$ commutes with the modular automorphism group $\sigma^\psi$,
so that the corresponding remarks above about automorphisms apply here as well. Another immediate consequence from the axioms is that 
for all $n \in \Nset$
\[
\alpha(\cM_{n-1}) \subset \cM_n
\]

Recall that a (noncommutative) random variable is an injective *-homomorphism $\iota$ into a probability space, i.e.
$\iota\colon \cA \to \cM$, see \cite{AFL82a}. 
We also write $\iota \colon (\cA,\psi_0) \rightarrow (\cM,\psi)$,
where $\psi_0 = \psi \circ \iota$.
For us $\cA$ is a von Neumann algebra and we include the property
that $\iota(\cA)$ is $\psi$-conditioned (see Definition \ref{def:nc-prob-space}) into our concept of random variables.

Generalizing terminology from classical probability we may say that
for an adapted endomorphism $\alpha$ the random variables
\begin{align*}
\iota_1 := \alpha\; |_{\cM_0} \colon & \quad \cM_0 \rightarrow \cM_1 \subset \cM \\
\iota_2 := \alpha^2 |_{\cM_0} \colon & \quad \cM_0 \rightarrow \cM_2 \subset \cM \\
\cdots & \\
\iota_n := \alpha^n |_{\cM_0} \colon & \quad \cM_0 \rightarrow \cM_n \subset \cM
\end{align*}
are adapted to the filtration $\cM_0 \subset \cM_1 \subset \cM_2 \subset \cdots$ and $\alpha$ is the time evolution of a stationary process. This explains our terminology. Product representations with (\ref{item:pr-i}) and (\ref{item:pr-ii}) give us a better grasp on adaptedness from an operator theoretic point of view. This idea has been introduced and examined in \cite[Chapter 3]{Gohm04a}, where more comments on the general philosophy and additional details can be found. In fact, the definition in \cite[Chapter 3]{Gohm04a} is a bit more general but Definition \ref{def:adapted-end} above seems to be more handy and, as the results in this paper show, it contains a lot of interesting examples. 

An important class of product representations arises in the following way:
Assume that there exists a sequence of unitaries $(u_n)_{n \in \Nset} \subset \cM^\psi$,
where $\cM^\psi$ is the centralizer, such that
\begin{align}
\alpha &= \lim_{n \to \infty} \Ad(u_1u_2 \cdots u_n)  
\tag{PR-0u} \label{item:pr-0-u}\\
u_n &\in \cM_n  
\tag{PR-1u} \label{item:pr-1-u} \\
u_{n+2} &\in \cM_n^\prime.  
\tag{PR-2u} \label{item:pr-2-u}
\end{align}
for all $n$.
Here $\cM_n^\prime$ denotes the commutant and the limit is taken in the pointwise $\sot$-sense.
It is easy to see that this provides us with an endomorphism $\alpha$ which is adapted with product representation. 

One of the most important procedures in classical probability is conditioning, and the framework of $\psi$-conditioned probability spaces allows us to do that in a natural way also in a noncommutative setting. We use a concept of conditional independence as introduced by K\"ostler in \cite{Koes08aPP}, where a more detailed discussion is given. 
\begin{Definition}\normalfont \label{def:independence}
Let $(\cM,\psi)$ be a probability space with three $\psi$-conditioned
von Neumann subalgebras $\cN,\, \cN_1$ and $\cN_2$. Then $\cN_1$ and $\cN_2$ 
are said to be
\emph{$\cN$-independent} or \emph{conditionally independent} if 
\[
E_{\cN}(xy) = E_{\cN}(x) E_{\cN}(y)
\]
for all $x \in \cN_1 \vee \cN$ and $y\in \cN_2 \vee \cN$. 
\end{Definition}
Note that $\cN_1$ and $\cN_2$ are $(\cN=\Cset)$-independent if and only if  
$
\psi(xy)= \psi(x) \psi(y)
$
for all $x \in \cN_1$ and $y \in \cN_2$ and in this case we recover the definition of K\"ummerer \cite{Kuem88a}. 

One of the attractive features of Definition \ref{def:independence}
from the point of view of operator algebras is the connection to commuting squares.
\begin{Lemma}\normalfont \label{lemma:comm-square}
Suppose $\cN \subset \cN_1 \cap \cN_2 \subset \cM$. Then the following are equivalent:
\begin{enumerate}
\item 
$\cN_1$ and $\cN_2$ are $\cN$-independent.
\item 
The following is a commuting square:
\begin{eqnarray*}
\begin{matrix}
\cN_2  & \subset & \cM  \\
\cup   &         & \cup \\
\cN  & \subset & \cN_1
\end{matrix}
\end{eqnarray*}
\end{enumerate}
\end{Lemma}
The notion of a commuting square has been introduced by Popa
\cite{Popa83b,Popa83a,PiPo86a}. There are many equivalent formulations:
\begin{itemize}
\item[(iii)] $E_{\cN_1}(\cN_2) = \cN$
\item[(iv)] $E_{\cN_1} E_{\cN_2} = E_{\cN}$
\item[(v)] $E_{\cN_1} E_{\cN_2} = E_{\cN_2} E_{\cN_1}$
and $\cN = \cN_1 \cap \cN_2$
\end{itemize}
etc. \\
Note in particular that the equality $\cN = \cN_1 \cap \cN_2$
is automatic in this situation. For further discussion see 
\cite{GHJ89a}, Proposition 4.2.1. The proof of the equivalences above 
is given in the tracial case there but it generalizes immediately to 
$\psi$-conditioned probability spaces. 
\begin{Definition}\normalfont \label{def:order-independence}
A family of $\psi$-conditioned subalgebras $(\cM_n)_{n \in \Nset_0}$, 
is said to be
\begin{enumerate}
\item[(CI)]
\emph{full $\cN$-independent} or \emph{conditionally full independent}, if 
$\bigvee_{i \in I}\cM_i$ and $\bigvee_{j \in J}\cM_j$ are $\cN$-independent 
for all $I,J \subset \Nset_0$ with $I \cap J = \emptyset$.
\item[(CI$_{\text{o}}$)]
\emph{order $\cN$-independent} or \emph{conditionally order independent}, if 
$\bigvee_{i \in I}\cM_i$ and $\bigvee_{j \in J}\cM_j$ are $\cN$-independent 
for all $I,J \subset \Nset_0$ with $I < J$ or $I > J$. (Here the order 
relation $I < J$ means $i < j $ for all $i \in I$ and $j \in J$.) 
\end{enumerate}
\end{Definition}
Clearly, (CI) implies (CI$_{\text{o}}$); but the converse is open in the 
generality of our setting. We will deliberately drop the attributes `full'
or `order' if we want to address conditional independence on an informal
level or if it is clear from the context.   
All these notions of independence translate to random variables by 
saying that random variables are independent if this is true for 
their ranges.   

This notion of conditional independence includes classical, tensor and free 
independence, including their amalgamated variants. Moreover, it goes beyond noncommutative independences with universality rules \cite{Spei97a,BeSc02a,NiSp06a}. It applies to all examples of generalized or noncommutative Gaussian random variables, as long as they respect the properties of a white noise functor \cite{Kuem96a, GuMa02a} and generate von Neumann algebras equipped with a faithful normal state (given by the vacuum vector of the underlying deformed Fock space), as they appear in \cite{Bosp91a,BoSp94a,BKS97a,GuMa02a,Krol02a}. We refer to \cite{HKK04a,Koes08aPP} for a more detailed treatment of conditional 
independence and for further examples coming from quantum probability. 

Now we are in a position to generalize many classical concepts related to 
independence to a noncommutative setting.
\begin{Definition}\normalfont \label{def:bernoulli}
Let $\cB_0$ be a $\psi$-conditioned subalgebra of $(\cM,\psi)$ and $\alpha$ 
a $\psi$-preserving endomorphism. Further let $\cN \subset \cB_0$ be a 
$\psi$-conditioned subalgebra which is pointwise fixed by $\alpha$
(i.e. $\cN \subset \cM^\alpha$). Then $\beta$ defined as the restriction of 
$\alpha$ to $\cB := \bigvee_{n\in\Nset_0} \alpha^n(\cB_0)$ is called a 
\emph{(full/ordered) Bernoulli shift over $\cN$ with generator $\cB_0$} if 
$\big( \beta^n(\cB_0) \big)_{n\in\Nset_0}$ is (full/order) $\cN$-independent.
\end{Definition} 
A trivial example is $\alpha = \id$ with $\cN = \cB_0 = \cM$, but in more
interesting examples the generator is usually small, often finite dimensional.
\begin{Remark}\normalfont
The above definition of a Bernoulli shift is equivalent to that given in
\cite{Koes08aPP}. This relies on the fact that a Bernoulli shift $\beta$ automatically commutes with the modular automorphism group: As required by the definition of independence, the ranges $\beta^n(\cB_0)$ are $\psi$-conditioned and so is $\cB$. Consequently, the modular automorphism group $\sigma_t^\psi$ restricts from $\cM$ to $\cB$.  But the ranges $\beta^n(\cB_0)$ are $\psi$-conditioned if and only if $\beta$ commutes with this restriction of $\sigma_t^\psi$.  
\end{Remark}
\begin{Lemma}\label{lem:bernoulli}
If $\beta$ is a (full/ordered) Bernoulli shift over $\cN$ then 
\[
\cN = \cB^\beta = \cB^\tail
\]
where $\cB^\beta$ is the fixed point algebra of $\beta$ and
$\cB^\tail := \bigcap_{n\in \Nset} \bigvee_{k \ge n } \beta^k(\cB_0)$.
In particular $\cN$ is uniquely determined by the endomorphism $\beta$.
\end{Lemma}
\begin{proof}
See \cite[Corollary 6.9]{Koes08aPP}.
\end{proof}
Noncommutative (ordered) Bernoulli shifts with $\cN = \Cset$ have been 
introduced by K\"ummerer \cite{Kuem88a} and further developed by Rupp 
in \cite{Rupp95a} where the connection to commuting squares is recognized
and schemes similar to the following Theorem \ref{thm:endo-pr} are considered. See also \cite{HKK04a} and \cite{Koes08aPP}.

How can we construct noncommutative Bernoulli shifts? It turns out 
that product representations as introduced in \ref{def:adapted-end} are a 
powerful tool for this task in the framework of conditional order independence.
\begin{Theorem} \label{thm:endo-pr}
In $(\cM,\psi)$ let $(\cM_n)_{n\in \Nset_0}$ be a tower of $\psi$-conditioned subalgebras such that
$\cM = \bigvee_{n\in\Nset_0} \cM_n$ and let $\alpha$ be an adapted endomorphism with product representation 
by factors $(\alpha_n)_{n \in \Nset}$, as in Definition \ref{def:adapted-end}. 
Further let $\cM_{-1}$ be a $\psi$-conditioned subalgebra of $\cM^\alpha \cap \cM_0$, where $\cM^\alpha$ is the fixed point algebra of $\alpha$ and suppose 
that for all $n \in \Nset$,
\begin{eqnarray*}
\begin{matrix}
\cM_{n-1} &\subset &  \cM_n \\
  \cup    &         &  \cup\\
\cM_{n-2} & \subset &   \alpha_n(\cM_{n-1})  
\end{matrix}
\end{eqnarray*}
is a commuting square. Then one obtains a triangular tower of inclusions such 
that all cells form commuting squares: 
\begin{eqnarray*}
\setcounter{MaxMatrixCols}{20}
\begin{matrix}
\cM_{-1} &\subset&   \cM_0 & \subset & \cM_1 & \subset &\cM_2 & \subset &\cM_3 & \subset  & \cdots & \subset & \cM\\
        &&          \cup  &         & \cup  &         & \cup &         & \cup  &       & & & \cup  \\
        &&   \cM_{-1}&\subset&\alpha(\cM_0)&\subset&\alpha(\cM_1)&\subset&\alpha(\cM_2)&\subset& \cdots & \subset & \alpha(\cM)\\
         &&               &         & \cup  &         & \cup   &         & \cup  &       & & & \cup \\
              &&&&   \cM_{-1}& \subset & \alpha^2 (\cM_0)  & \subset & \alpha^2(\cM_1)
              & \subset  & \cdots & \subset & \alpha^2(\cM) \\
       &&&&&&   \cup           &         & \cup  &       &&  & \cup \\
        &&&&&&   \vdots           &         & \vdots  &       &&  & \vdots 
\end{matrix}.
\setcounter{MaxMatrixCols}{10}
\end{eqnarray*}
If $\cB_0$ is a $\psi$-conditioned subalgebra such that $\cM_{-1} \subset \cB_0 \subset \cM_0$ then $\beta$ defined as the restriction 
of $\alpha$ to $\cB := \bigvee_{n\in\Nset_0} \alpha^n(\cB_0)$ is an ordered Bernoulli shift over $\cN$ with generator $\cB_0$, where $\cN$ is given by 
\[
\cN = \cM_{-1} = \cM^\alpha = \cM^\tail = \cB^\beta = \cB^\tail.
\]
\end{Theorem}
\begin{proof}
The proof of Theorem \ref{thm:endo-braid-i} is written up in a way
that easily transfers to the present setting and yields the triangular tower 
of commuting squares given above. Now suppose $I < J$, say $i \le n < j$ for all 
$i \in I$ and $j \in J$. By adaptedness $\bigvee_{i\in I} \beta^i(\cB_0) 
\subset \cM_n$ and $\bigvee_{j \in J} \beta^j(\cB_0) \subset 
\alpha^{n+1}(\cM)$. Inspection of the triangular tower of commuting squares 
shows that these algebras are order $\cM_{-1}$-independent. This proves that $\beta$ is an ordered Bernoulli shift over $\cM_{-1}$, and hence $\cN= \cM_{-1}$. The equalities $\cN = \cB^\beta = \cB^\tail$ follow from Lemma \ref{lem:bernoulli}. The fixed point algebra $\cM^\alpha$ and the tail algebra $\cM^\tail$ cannot be strictly bigger than $\cM_{-1}$ because it is readily checked from (\ref{item:pr-i}), (\ref{item:pr-ii}) and the commuting squares that $\alpha$ maps the $\psi$-orthogonal complement of $\cM_{k-1}$ in $\cM_k$ into the $\psi$-orthogonal complement of $\cM_k$ in $\cM_{k+1}$, for all $k \in \Nset_0$. But this means that the isometry on the GNS Hilbert space induced by $\alpha$ is a shift operator on the orthogonal complement of $\cM_{-1}$ and hence there can be no fixed point or tail of $\alpha$ outside $\cM_{-1}$. 
\end{proof}
The last assertion of Theorem \ref{thm:endo-pr} is very useful in applications.
Suppose we start with a Bernoulli shift $\beta$ over $\cN$ with generator $\cB_0$. We know from Lemma \ref{lem:bernoulli} that $\cN$ equals $\cB^\beta = \cB^\tail$; but how can one effectively determine $\cB^\beta$ or $\cB^\tail$ in applications? As soon as we succeed to realize $\beta$ as the restriction of
a shift $\alpha$ in the way of Theorem \ref{thm:endo-pr}, i.e. with $\cN \subset \cM^\alpha$ and $\cB_0 \subset \cM_0$, then $\cN = \cB^\beta = \cB^\tail$ must be equal to the left lower corner $\cM_{-1}$ of the first commuting square. In general this commuting square is easier to access in applications than tail or fixed points of $\beta$.  So, besides a nice general structure, Theorem \ref{thm:endo-pr} provides a convenient way to identify $\cN$. We demonstrate this idea for some examples in Section \ref{section:examples}.
  
An additional strong point of Theorem \ref{thm:endo-pr} is that it identifies
$\cM_{-1}$, the first algebra in the given tower, as the fixed point algebra of the endomorphism $\alpha = \lim_{m\to \infty}\alpha_{1}\alpha_{2} \alpha_{3}\cdots \alpha_{m}$. Dropping the first factors in this product representation, we can also identify the other algebras in the tower as fixed point algebras of certain partial shifts:
\begin{Corollary} \label{cor:endo-pr}
With assumptions as in Theorem \ref{thm:endo-pr}, for all $n \ge -1$ the 
algebra $\cM_n$ is equal to the fixed point algebra of $\alpha_{n+2,\infty} 
:= \lim_{m\to \infty}\alpha_{n+2}\alpha_{n+3} \alpha_{n+4}\cdots \alpha_{n+m}$ 
and equal to the tail algebra for the random sequence produced by 
$\alpha_{n+2,\infty}$ from $\cM_{n+1}$ as range of the time $0$-random 
variable. (Compare with Section \ref{section:sequences-braid}.)
\end{Corollary}
\begin{proof}
Apply Theorem \ref{thm:endo-pr} for $\alpha_{n+2,\infty}$ and the
tower $\cM_n \subset \cM_{n+1} \subset \cdots$.
\end{proof}
Corollary \ref{cor:endo-pr} shows that, given the commuting squares in 
Theorem \ref{thm:endo-pr}, the whole tower
$\cM_{-1} \subset \cM_0 \subset \cdots$ is determined in terms of fixed point algebras by the factors $(\alpha_n)_{n\in\Nset}$ of the product representation. The following consequence for braid group representations can be used to determine these fixed point algebras explicitly, compare with Section \ref{section:examples}.
\begin{Theorem} \label{thm:tower-reconstruction}
Let the setting be as in Theorem \ref{thm:endo-pr}. In addition assume that
for all $k \in \Nset$ we have $\alpha_k = \rho(\sigma_k)$ for a
representation $\rho \colon \Bset_\infty \to \Aut{\cM,\psi}$.
Then for all $n \ge -1$
\[
\cM_n = \cM^\rho_n,
\]
using the notation $\cM^\rho_n = \cM^{\rho(\Bset_{n+2,\infty})}$
introduced in Section \ref{section:sequences-braid}.
As a special case, if $\rho(\sigma_k) = Ad\, u_k$ for all $k \ge n+2$ then
$\cM_n$ is equal to the relative commutant 
$\cM \cap \set{u_k}{k \ge n+2}^\prime$.
\end{Theorem}
\begin{proof}
Combine Corollary \ref{cor:endo-pr} with Theorem \ref{thm:braided-HS}
for the shifted representation.
\end{proof}
In other words, under the given assumptions the tower we started from is 
automatically the tower of fixed point algebras. 
Note that this is a kind of converse to Theorem \ref{thm:braid-tower}
where we proved that starting from a braid group representation we can 
always construct commuting squares for the tower of fixed point algebras 
and establish the situation obtained in Theorem \ref{thm:tower-reconstruction}
from another direction. This is why we think of this structure as a kind of 
Galois type tower (compare with Remark \ref{rem:galois}).

Let us finally rewrite Theorem \ref{thm:tower-reconstruction} to obtain 
a more explicit form of the assumptions.
\begin{Corollary} \label{cor:tower-reconstruction}
Given a representation $\rho \colon \Bset_\infty \to \Aut{\cM,\psi}$
and a tower $(\cM_n)_{n \ge -1}$ of $\psi$-conditioned subalgebras such 
that $\cM = \bigvee_{n \ge -1} \cM_n$. Suppose further that, for all
$n \ge -1$, we have 
\[
\cM_n \subset \cM^\rho_n
\]
and we have a commuting square
\begin{eqnarray*}
\begin{matrix}
\cM_{n+1} &\subset &  \cM_{n+2} \\
  \cup    &         &  \cup\\
\cM_n & \subset &   \rho(\sigma_{n+2})(\cM_{n+1})  
\end{matrix}.
\end{eqnarray*}
Then for all $n \ge -1$
\[
\cM_n = \cM^\rho_n
\]
\end{Corollary}
\begin{proof}
Note that the corresponding endomorphism 
$\alpha = \lim_{n\to\infty} \rho(\sigma_1 \cdots \sigma_n)$ 
has already been defined and studied in Section \ref{section:endomorphisms-braid}. 
The assumption $\cM_n \subset \cM^\rho_n$ for all $n \ge -1$ gives 
$\cM_{-1} \subset \cM^\alpha$ for $n=-1$ and (\ref{item:pr-ii}) for 
$n\in\Nset_0$. Now (\ref{item:pr-ii}) together with the fact that 
$\alpha_j = \rho(\sigma_j)$ and $\alpha_k = \rho(\sigma_k)$ commute for 
$|j-k|\ge 2$ imply (\ref{item:pr-i}). Hence we have verified all 
assumptions of Theorem \ref{thm:tower-reconstruction}. 
\end{proof}
\end{appendix}
\bibliographystyle{alpha}                 
\label{section:bibliography}
\bibliography{ref-braid}                   
\end{document}